\documentclass{article}

\usepackage{graphicx} % Required for inserting images
\usepackage[utf8]{inputenc}
\usepackage{authblk} %for afiliations
\usepackage[dvipsnames]{xcolor}
\usepackage[colorlinks=true, 
urlcolor=blue,
citecolor=red]{hyperref}
\usepackage{booktabs}
\usepackage{placeins} % For float barriers

\usepackage{amsmath,amsthm,amssymb}
\usepackage{bm} %for \bm 
\usepackage{tikz}
\usepackage[top=1.2in, left=0.9in, bottom=1.2in, right=0.9in]{geometry}
\usepackage{color,soul}
\sethlcolor{pink}
\usepackage[colorinlistoftodos,prependcaption,textsize=small]{todonotes}

\newtheorem{prop}{Proposition}

\title{A mixed interpolation-regression method for function approximation on certain planar domains}
\author[1]{Ruymán Cruz-Barroso}
\author[2]{Lidia Fernández}
\author[3]{Francisco Marcellán}
\author[2]{Juan Antonio Villegas}

\affil[1]{Departamento de Análisis Matemático and Instituto de Matemáticas y sus Aplicaciones (IMAULL), La Laguna University, Spain}
\affil[2]{IMAG and Departamento de Matem\'atica Aplicada, Universidad de Granada, Spain}
\affil[3]{Departamento de Matemáticas, Universidad Carlos III de Madrid, Spain}
\date{\today}

\makeatletter
\newcommand{\keywords}[1]{%
  \begingroup
  \par\vspace{0.5em}
  \noindent\textbf{Keywords: }#1\par
  \endgroup
}
\newcommand{\MSC}[1]{%
  \begingroup
  \par\vspace{0.25em}
  \noindent\textbf{MSC (2010): }#1\par
  \endgroup
}
\makeatother

\begin{document}

\maketitle

\begin{abstract}
In this contribution, we introduce a mixed interpolation-regression operator for functions defined on certain planar domains. We focus  on an ellipse, an annulus and a polygon. An upper bound for the operator is obtained. Cubature formulas for weight functions defined on these domains are studied. The performance of the interpolation-regression methods is illustrated by numerical examples. 
\end{abstract}


\MSC{33C45; 42C05; 65D32; 62J05}

\section{Introduction}

In \cite{DDM15}, the authors introduce a method that combines interpolation and regression of a function in one variable whose values are known at a given set of points on the real line. The central idea is to avoid the Runge phenomenon when selecting the nodes used for interpolation, while still retaining the remaining data rather than discarding it. This additional information is then incorporated to improve the quality of the approximation obtained from the interpolation step.

In the case of two or more variables, the problem of reconstructing a function by a polynomial approximation that makes full use of data sampled on a regular grid has been explored through the development of the constrained mock‑Chebyshev least‑squares operator for square domains \cite{DDN22_1,DMN24}, and later extended to triangular settings \cite{DDN24}. The underlying principle of this methodology is the same as in the one variable case, so the authors  carefully select a subset of nodes arranged to approximate an optimal interpolation pattern—such as Padua points in the square or Leja sequences in the triangle. The remaining nodes are then incorporated into a simultaneous regression step that serves to refine and improve the accuracy of the overall approximation.
This strategy has been successfully applied to implement quadrature formulas \cite{DDN22_2}, to solve Fredholm integral equations \cite{DMNO24} and to approximate derivatives \cite{DN24}, and has recently been extended to Hermite interpolation \cite{DMN26} and to the unit circle \cite{CFM26} by some of the authors by considering points close to the zeros of para-orthogonal polynomials, among other problems.

In \cite{DMN25}, the authors extended this method to functions defined on a disk. They used the basis of Zernike polynomials, so the selection of the interpolation nodes became crucial in order to guarantee numerical stability. The function approximation technique proposed therein is also used to define accurate cubature formulas based on an interpolation-regression operator. An upper bound of the sup-norm of the operator is obtained  as well as an estimate of the error function. Numerical tests for different examples of sampling points are implemented.  

Several works, such as \cite{DMR26, FLNP15, NLDP14}, address the problem of constructing orthogonal bases analogous to Zernike polynomials on regions other than  a unit disk, including an ellipse, a polygon or an annulus. The main objective of the present work is to study the approximation of a function defined on any of these three domains by means of an interpolation–regression method of the type described above. To that end, we employ a mapping of the Zernike polynomial basis in such a way that the resulting functions are not polynomials but functions orthogonal  with respect to measures defined on the domain.   

Two principal strategies are commonly used to construct orthogonal bases on compact planar domains. The first one consists of generating the basis through Gram–Schmidt orthogonalization, but this approach can suffer from severe numerical instability when high‑degree polynomials are involved.
A second strategy addresses these difficulties by transplanting an already established orthogonal basis to a different domain by means of a diffeomorphism. This idea was explored in \cite{NLDP14}, where the Zernike polynomials defined on the unit disk were taken as a reference and subsequently mapped onto more general geometries.\\

The structure of the paper is as follows. In Section~\ref{sec:approximation-disk} we provide basic background on polynomials on the disk and emphasize that the choice of interpolation nodes  plays a key role. In a next step, if the value of a function is known at the optimal nodes, then it is possible to construct an accurate Zernike-based polynomial interpolant. In Section~\ref{sec:extension-zernike}, the starting point is Zernike polynomials as an orthonormal system on the disk. We apply a diffeomorphism $\varphi$ from the disk to the target region (ellipse, annulus, among others) in order to get a new basis. In this way, using the Jacobian of $\varphi^{-1}$, the basis remains orthogonal with respect to a perturbed measure in the domain and constitutes a complete system. This construction is illustrated when the domain is an ellipse, a circular annulus, and a regular polygon. In Section~\ref{sec:Opt_nodes}, the sampling nodes in such domains are given using  the above diffeomorphism. In Section~\ref{sec:mixed-IR-method} the mixed interpolation-regression operator is described for these three domains. Cubature formulas based on Gaussian nodes in the disk by using a change of variable are stated in the above three domains in Section~\ref{sec:cubature}. Finally, in Section~\ref{sec:experiments} the performance of the interpolation-regression method is illustrated for certain choices of the ellipse semiaxes, the inner and outer radii of the circular annulus and a 12-gon, respectively. The mean squared error, the maximum absolute error and the mean (resp. maximum) relative error for a function defined in the above domains are estimated for calculations based on the corresponding interpolation-regression operator.

\section{Approximation on the unit disk}
\label{sec:approximation-disk}

The unit disk, denoted by $D=\{\mathbf x\in\mathbb R^2: \|\mathbf x\|\leq 1\},$ where $\| \cdot \|$ denotes the Euclidean norm in $\mathbb{R}^2$, has been a widely studied domain in approximation theory. Given the variety of scientific fields in which functions defined on circular domains are employed, the efficient and accurate approximation of functions defined on the disk becomes a highly relevant problem. 

In  \cite[Section 8.1]{DX14}, generalized classical orthogonal polynomials on the $d$-dimensional ball are defined and, in \cite[Section 9.3]{DX14},  the convergence of Fourier expansions in terms of sequences of  polynomials orthogonal with respect to admissible weight functions supported on the $d$-dimensional ball is studied for functions belonging to weighted $L^{p}$-spaces.

In this section, we address the polynomial approximation of real functions on the disk. Denote by $\mathbb P_{n}(D)$ the linear space of polynomials in two variables of degree at most $n$. Within the orthogonal polynomial setting, there are several orthogonal bases of $\mathbb P_{n}(D)$. These include representations in terms of Gegenbauer polynomials, Jacobi polynomials in combination with spherical harmonics, and those computable via Rodrigues formulas; see \cite[Section 2.3]{DX14}. Although these bases are particularly effective depending on the application, we will focus on the Zernike polynomials. These are particularly significant in disk-based approximation and are also used in physical optics, ophthalmology, or astronomy; see \cite{BDF25, MDCZA09}. 

In \cite{NT22} an updated review on Zernike polynomials is presented. It provides a comprehensive treatment of their history, definitions, mathematical properties, roles in wavefront fitting, relations to optical aberrations, optical design, ophthalmic optics, and image analysis, among others, as well as a survey of their scientific and industrial applications.

We formally introduce Zernike polynomials on the disk in the following subsection.

\subsection{Zernike polynomials on the disk}

Zernike polynomials constitute a complete orthogonal basis on the unit disk~\cite{ZS34}. Originally introduced by Frits Zernike for the analysis of phase-contrast microscopy images, they quickly became fundamental in optical science because most optical aberrations are defined on circular pupils. Although they are not the only family of orthogonal polynomials on the disk~\cite[Section 2]{DX14}, their structure aligns closely with classical optical aberrations. They are naturally expressed in polar coordinates ($\rho$, $\phi$), with $\rho\in [0,1]$, $\phi\in[0,2\pi)$. Each Zernike polynomial is given as the product of a radial term and an angular term, together with appropriate normalization factors:
\begin{equation}\label{eq:zernike}
\begin{aligned}
&Z_{m}^{l}(\rho, \phi) = N_{m}^{l} R_{m}^{l} (\rho) \cos(m\phi),\quad  m\geq0,
\\
&Z_{m}^{-l}(\rho, \phi) = N_{m}^{l} R_{m}^{l} (\rho) \sin(m\phi), \quad m\geq1,
\end{aligned}
\end{equation}
for $0\leq l \leq m$ and $m-l$ even.
The normalization factors satisfy $N_{m}^{l}=\sqrt{2(m+1)}$ if $l\neq0$ and $N_{m}^{l}=\sqrt{m+1}$ if $l=0$. The functions $R_{m}^{l}$ are polynomials of degree $m$ in the variable $\rho$, defined by:
\begin{equation}\label{eq:radial}
R_{m}^{l}(\rho) = 
\sum_{i=0}^{(m-l)/2} \frac{(-1)^i(m-i)!}{i!\left((m+l)/2-i\right)!\left((m-l)/2-i \right)!}\rho^{m-2i},
\ \ m,l\ge 0, 
\ \ m-l \text{ even}.
\end{equation}

Here, it is important to emphasize that these polynomials are related to the classical Jacobi polynomials with the standard normalization \cite{Sze75} as
$$
R_{m}^{l}(\rho) = (-1)^{(m-l)/2}\rho^l P_{(m-l)/2}^{(l,0)}(1-2\rho^2),
$$
for $m,l \ge 0$ and $m-l$ even. With a change of variables to Cartesian coordinates $x=\rho \cos(\phi)$, $y=\rho \sin(\phi)$, we see that Zernike polynomials are, indeed, polynomials in $x$ and $y$ of total degree $m$, and are orthonormal on the unit disk $D=\{(x,y)\in\mathbb{R}^2 : x^2+y^2\leq1\}$; that is, 
\[
\int_{D} Z_m^l(x,y) Z_{m'}^{l'}(x,y) \, dx \, dy = \delta_{m,m'}\, \delta_{l,l'}. 
\]
It is possible to index Zernike polynomials with a single index $j$,
$$
j=\frac{m(m+2)+l}{2},
$$
so we can write
$$
Z_j(\rho, \phi)=Z_{m}^{l}(\rho, \phi).
$$

With Zernike polynomials as the central element, in \cite{DMN25} the authors developed a method to compute an approximation of a function defined on $D$ from a finite discrete set of points where the function is known, using interpolation and regression techniques simultaneously. This methodology is introduced in the following sections.

\subsection{Optimal nodes for interpolation on the disk}
To maintain numerical stability when interpolating functions with the Zernike polynomial basis, the choice of interpolation nodes $\{\mathbf p_1,\dots,\mathbf p_M\}$ on the disk plays a decisive role. A standard criterion for evaluating a node distribution is the condition number of the collocation matrix formed by sampling the Zernike basis at the selected nodes. Among the many unisolvent point sets available on the disk, one of the most effective is the Bos array \cite{BX03,DMR26, DMN25}. This configuration places the nodes on concentric circles determined by a prescribed maximal radial order $m$.

More precisely, 
$$
m_\nu=2m-4\nu+5
$$
equispaced nodes are placed on each circle of radius $\rho_\nu$, where 
$$
1\ge \rho_1 > \rho_2> \dots > \rho_K \ge 0, \quad K=K(m)=\left\lfloor \frac{m}{2}\right\rfloor+1.
$$

Regarding the choice of these radii, the objective is to minimize the conditioning of the Gram matrix $(Z_{j}(\mathbf p_i))_{i,j=1}^M$ and, in turn, the conditioning of the interpolation problem. Unfortunately, to date, an efficient way of obtaining the optimal configuration does not exist, so we use the quasi-optimal radii $\rho_\nu$, which can be computed by the following formula (see \cite{FLNP15, RSFM16}):
\begin{equation}\label{eq:rho_nu}
\rho_\nu=\rho_\nu(m)=1.1565\xi_{\nu,m}-0.76535\xi_{\nu,m}^2+0.60517\xi_{\nu,m}^3,
\end{equation}
where $\xi_{\nu,m}$ are the zeros of the $(m+1)$-th Chebyshev polynomial of the first kind, that is,
$$
\xi_{\nu,m}=\cos\left(\frac{(2\nu-1)\pi}{2(m+1)} \right), \qquad \nu=1,\dots,K.
$$

\begin{figure}
    \centering
    \begin{tabular}{cc}
       \includegraphics[width=0.3\linewidth]{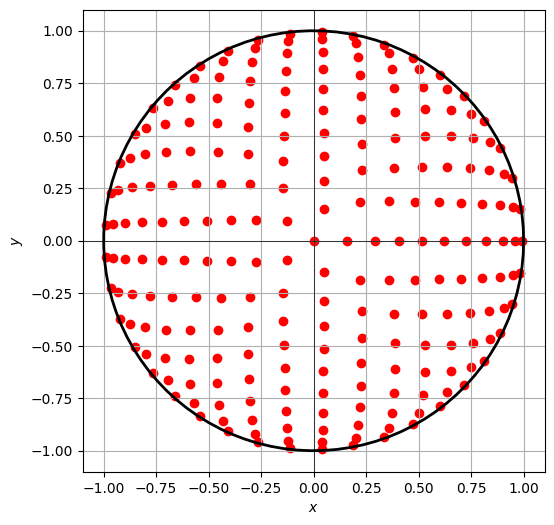}  & \includegraphics[width=0.3\linewidth]{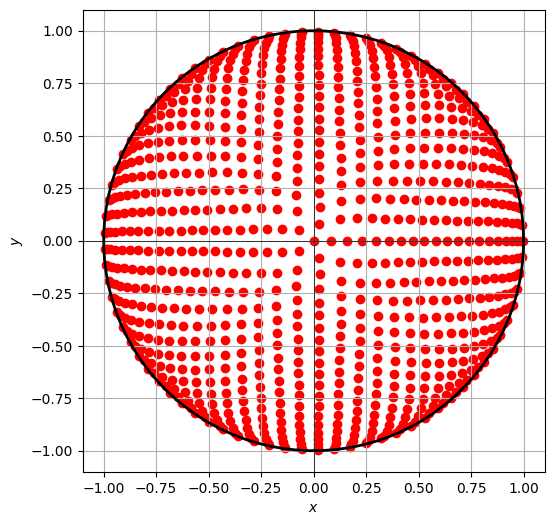}  \\
        (a) $m=20$ & (b) $m=40$ 
    \end{tabular}
    \caption{Optimal Concentric Samplings on the Disk for $m=20,40$.}
    \label{fig:OCS-disk}
\end{figure}

On each circle of radius $\rho_\nu$, the nodes are defined by
\begin{equation}\label{eq:p_opt}
\mathbf{p}_{\nu,\sigma}=\left(\rho_\nu \cos\left(\frac{2\pi (\sigma-1)}{m_\nu} \right), \, \rho_\nu \sin\left(\frac{2\pi (\sigma-1)}{m_\nu}  \right)\right),
\end{equation}
where $\sigma=1,\dots,m_\nu$, $\nu=1,\dots, K$. This configuration is referred to as the \emph{Optimal Concentric Sampling} (OCS); see Figure~\ref{fig:OCS-disk} for a graphical representation of the distribution of these points.
The set of optimal nodes is denoted by
$$
\Sigma_m^{opt}=\left\{\mathbf{p}_{\nu,\sigma} \,:\,\sigma=1,\dots,m_\nu, \, \nu=1,\dots, K  \right\},
$$
and 
$$M:=\sum_{\nu=0}^K m_\nu=\frac{(m+1)(m+2)}{2}
$$
is the dimension of $\mathbb{P}_m(\mathbb{R}^2)$.

In situations where the available data do not coincide with the optimal interpolation nodes, a different strategy is required. In the next subsection, we introduce an interpolation–regression approach that allows one to construct accurate Zernike-based approximations from general point sets on the disk.

\subsection{Interpolation-regression approach on the disk}\label{Sec:IntReg_Disk}
If the function $f$ is known at these optimal nodes, then it is possible to construct an accurate Zernike polynomial interpolant using them. However, if this is not the case and the function is known on a general grid, a mixed interpolation-regression method can be employed. The content of this section is developed in \cite{DMN25}.

Let us consider the set
\begin{equation}
\label{eq:Snn}
S_{n,n}=\left\{\mathbf{x}_{\eta,\kappa}=(\rho_\eta \cos(\theta_\kappa), \rho_\eta \sin(\theta_\kappa))\, : \, \eta=0,\dots, n, \,\, \kappa=0,\dots,n \right\} \subset D,
\end{equation}
where $0<\rho_0<\rho_1<\dots<\rho_n$ and $0\le \theta_0<\theta_1<\dots<\theta_{n} <2\pi$ are increasing sequences of radii and angles, respectively. Then 
$$
\# (S_{n,n})=(n+1)^2:=N.
$$
The set $S_{n,n}$ consists of points expressed in polar coordinates, but the authors argue that this is only for simplicity and that all the reasoning remains valid if we consider randomly sampled points uniformly distributed.

In order to define the interpolation-regression operator, let us fix $m, \tilde{r}\in\mathbb{N}$ such that
$$
m<\tilde{r}, \quad \tilde{R}=\frac{(\tilde{r}+1)(\tilde{r}+2)}{2}<N, \quad M=\frac{(m+1)(m+2)}{2}.
$$
The idea is to construct a polynomial of degree $\tilde{r}$ that interpolates the $M$ nodes closest to those of $\Sigma_m^{opt}$, while using the remaining nodes to improve the approximation accuracy through a simultaneous regression process.

The set of \textit{mock-optimal nodes} is 
\begin{equation}
\label{eq:mock-optimal-points}
\Sigma'_m :=\left\{\mathbf{x}'_{\nu,\sigma} \, : \, \sigma=1,\dots,m_\nu, \, \nu=1,\dots, K  \right\} \subset S_{n,n},
\end{equation}
where $\mathbf{x}'_{\nu,\sigma}$ is the solution to the minimization problem
\begin{equation}
\label{eq:minimization-problem}
\min_{\mathbf{x}_{\eta,\kappa}\in S_{n,n}} \|\mathbf{x}_{\eta,\kappa} - \mathbf{p}_{\nu,\sigma} \|.
\end{equation}
In \cite{DMN25}, an algorithm to select these points is presented, and it is guaranteed that
$\# (\Sigma'_m)=M$.

To simplify the notation,  a single-index notation is used for the nodes and the Zernike polynomials; that is,
\begin{align}
\mathbf{x}_{\eta,\kappa} &\rightarrow \mathbf{x}_i,  \quad i=1,\dots,N,
\\
\mathbf{x}'_{\nu,\sigma} &\rightarrow \mathbf{x}'_i,  \quad i=1,\dots,M,
\\
\mathbf{p}_{\nu,\sigma} &\rightarrow \mathbf{p}_j,  \quad j=1,\dots,M,
\\
Z_n^m &\rightarrow Z_j,  \quad j=1,\dots,M.
\end{align}

Let $\left\{u_1,\dots,u_{\tilde{R}}\right\}$ be a basis of $\mathbb{P}_{\tilde{r}}(D)$, the space of polynomials in two variables with total degree less than or equal to $\tilde{r}$, and suppose that the set $S_{n,n}$ is such that the matrix
$$
\mathcal{M}:=\left[
\begin{array}{cccc}
u_1(\mathbf{x}_1) & u_2(\mathbf{x}_1) & \dots & u_{\tilde{R}}(\mathbf{x}_1)
\\
u_1(\mathbf{x}_2) & u_2(\mathbf{x}_2) & \dots &u_{\tilde{R}}(\mathbf{x}_2)
\\
\vdots & \vdots & \ddots & \vdots
\\
u_1(\mathbf{x}_N) & u_2(\mathbf{x}_N) & \dots & u_{\tilde{R}}(\mathbf{x}_N)
\end{array}
\right] \in \mathbb{R}^{N\times \tilde{R}}
$$
has rank  $\tilde{R}$.

Let us also assume that the nodes of $\Sigma'_m$ are sufficiently close to those of the set $\Sigma^{opt}_m$ and, as a consequence, that the  matrix
$$
\mathcal{C}:=\left[
\begin{array}{ccccc}
u_1(\mathbf{x}'_1) & u_2(\mathbf{x}'_1) & \dots & u_{\tilde{R}}(\mathbf{x}'_1)
\\
u_1(\mathbf{x}'_2) & u_2(\mathbf{x}'_2) & \dots & u_{\tilde{R}}(\mathbf{x}'_2)
\\
\vdots & \vdots & \ddots & \vdots
\\
u_1(\mathbf{x}'_M) & u_2(\mathbf{x}'_M) & \dots & u_{\tilde{R}}(\mathbf{x}'_M)
\end{array}
\right] \in \mathbb{R}^{M\times \tilde{R}}
$$
has rank $M$. We can reorder the points in $S_{n,n}$ so that the mock-optimal nodes are the first $M$ ones, and thus $\mathcal{C}$ is the matrix formed by the first $M$ rows of the matrix $\mathcal{M}$. 

Let
$$
\bm{b}:=\big[f(\mathbf{x}_1), \dots, f(\mathbf{x}_N) \big]^T, \quad  \bm{d}:=\big[f(\mathbf{x}'_1), \dots, f(\mathbf{x}'_M) \big]^T.
$$
Then, an interpolation-regression operator is defined by
\begin{equation}
\hat{\Pi}_{\tilde{r},n}: f\in C(D,\mathbb{R}) \rightarrow \hat{\Pi}_{\tilde{r},n}[f]:=\sum_{i=1}^{\tilde{R}} \hat{a}_i u_i \in \mathbb{P}_{\tilde{r}}(\mathbb{R}^2),
\end{equation}
where the vector $\hat{\boldsymbol{a}}:=\big[\hat{a}_1, \hat{a}_2, \dots, \hat{a}_{\tilde{R}} \big]^T$
is the solution to the following least-squares problem: 
$$
\min_{\bm{a}\in\mathbb{R}^{\tilde{R}}} \|\mathcal{M}\bm{a}-\bm{b}\|_2,\quad \text{subject to } \quad \mathcal{C}\bm{a}=\bm{d}.
$$
It can be computed by solving the Karush-Kuhn-Tucker linear system
$$
\left[\begin{array}{cc}
   \mathcal{W}  &  \mathcal{C}^T
   \\
    \mathcal{C} & 0
\end{array}
\right] 
\left[\begin{array}{c}
   \hat{\bm{a}}  
   \\
    \hat{\bm{z}} 
\end{array}
\right]
=
\left[\begin{array}{c}
   2 \mathcal{M}^T\bm{b}  
   \\
    \bm{d} 
\end{array}
\right], 
\qquad
\mathcal{W}= 2 \mathcal{M}^T \mathcal{M},
$$
where $\hat{\bm{z}}$ is the vector of Lagrange multipliers.

To solve the system, we compute the QR factorization of $\mathcal{C}$, 
$$
\mathcal{C}=\mathcal{Q_C R_C}, \quad \mathcal{Q_C}\in\mathbb{R}^{M\times M}, \quad \mathcal{Q_C} \mathcal{Q_C}^T= \mathcal{I},
$$
and the matrix $\mathcal{R_C}$ is decomposed as
$$
\mathcal{R_C}=\big[\mathcal{R}_{11}, \mathcal{R}_{12} \big],  \quad \mathcal{R}_{11}\in\mathbb{R}^{M\times M}, \quad \mathcal{R}_{12} \in \mathbb{R}^{M\times (\tilde{R}-M)}.
$$
If we write $\hat{\boldsymbol{a}}:=
\left[\begin{array}{c}
\hat{\boldsymbol{a}}_1 
\\ \hat{\boldsymbol{a}}_2
\end{array}\right]$, with $\hat{\boldsymbol{a}}_1\in\mathbb{R}^M$ and  $\hat{\boldsymbol{a}}_2\in\mathbb{R}^{\tilde{R}-M}$ and
$$
\mathcal{M}=\big[\mathcal{M}_1, \mathcal{M}_2 \big], \qquad \mathcal{M}_1\in\mathbb{R}^{N\times M}, \quad\mathcal{M}_2 \in\mathbb{R}^{N\times (\tilde{R}-M)},
$$
then
\begin{align}
\hat{\bm{a}}_1 & = \mathcal{R}_{11}^{-1}\Big(\mathcal{Q_C}^T \bm{d}-\mathcal{R}_{12}\hat{\bm{a}}_2 \Big),
\\
\hat{\bm{a}}_2 & = \big(\mathcal{V}_1^T \mathcal{V}_1 \big)^{-1} \mathcal{V}_1^T \bm{b}_1,
\end{align}
where
$$
\mathcal{V}_1=\mathcal{M}_2-\mathcal{M}_1 \mathcal{R}_{11}^{-1} \mathcal{R}_{12}\in \mathbb{R}^{N\times (\tilde{R}-M)}, \qquad \bm{b}_1=\bm{b}-\mathcal{M}_1\mathcal{R}_{11}^{-1} \mathcal{Q_C}^T \bm{d}.
$$
For more details, see \cite{Bjo96,DMN25}. \\

The framework developed in this section provides the necessary tools for stable and accurate approximation on the unit disk, combining Zernike polynomials with suitable sampling and reconstruction strategies. In many applications, however, the domain of interest is not the disk but a more general planar region. A natural way to extend these techniques is to transport them through appropriate mappings. In the next section, we introduce a methodology to extend Zernike polynomials to more general domains by means of diffeomorphisms (almost everywhere), preserving their main structural properties.

\section{Extension of Zernike polynomials}
\label{sec:extension-zernike}

In \cite{NLDP14}, the authors present a unified approach for constructing complete and orthonormal systems of functions on optical apertures with a wide range of shapes. They develop a general method based on applying a diffeomorphism to transform the unit disk into the target aperture. In this way, well‑known orthonormal systems defined on the disk, such as  Zernike polynomials, can be transported to more complex geometries while preserving most of their mathematical structure and physical interpretation.

They define the appropriate mapping from the unit disk to an angular sector of an elliptical annulus. This is the most general case considered, and all other standard aperture shapes arise as special cases obtained by fixing specific parameter values. The authors derive explicit formulas for the mapping, its inverse, and its Jacobian, and show that the transformed functions remain orthonormal with respect to the metric induced by the change of variables.

The general idea of the method is as follows:
\begin{enumerate}
    \item Start with an orthonormal system on the disk (e.g., Zernike polynomials).
    \item Apply a diffeomorphism $\varphi : D \to \Omega$ that transforms the unit disk into the target region $\Omega$ (ellipse, annulus, sector, etc.).
    \item Define new basis functions on $\Omega$ by
    $$
        K_j(x,y) = Z_j\bigl(\varphi^{-1}(x,y)\bigr).
    $$
    \item Using the Jacobian of the mapping $\varphi^{-1}$, these new functions remain orthonormal with respect to the weight $|J_{\varphi^{-1}}(x,y)|\,dx\, dy$, and they constitute a complete system on $\Omega$.
\end{enumerate}

Two cases of particular interest in optics are analyzed in detail: elliptical apertures and annular apertures. These are studied here.

In \cite{FLNP15}, the authors extend the previous results to the case where the domain is a diffeomorphism almost everywhere and include the case where $\Omega$ is a polygon.

\bigskip

\subsubsection*{Elliptical Zernike basis}
In the case of an ellipse, the mapping from the unit disk to the ellipse is an affine transformation involving a scaling along the coordinate axes and a rotation by  an angle $\alpha$ (see \cite{NLDP14}).
An ellipse with semi-axes $A \ge B>0$ and orientation $\alpha$ is obtained from $(u,v)\in D$ through
$$
\varphi(u,v)= (Au \cos\alpha -B v \sin\alpha,  Au \sin\alpha + B v\cos\alpha).   
$$
The inverse transformation is also linear:
$$
\varphi^{-1}(x,y)=
\left(
\frac{x\cos\alpha + y\sin\alpha}{A},\
\frac{-x\sin\alpha + y\cos\alpha}{B}
\right).
$$
We will focus  on ellipses aligned with the axes. In this case 
\begin{equation}
    \label{eq:map-ellipse}
    \varphi(u,v)= (Au,  B v), \qquad \varphi^{-1}(x,y)=
\left(
\frac{x}{A},\
\frac{y}{B}
\right). 
\end{equation}
The Jacobian of the transformation $\varphi^{-1}$ is constant, $J(x,y)=1/(AB)$, so the family of polynomials $\{E_j \}_{j\ge 0}$ given by
$$
E_j(x,y):=\frac{1}{\sqrt{AB}} \,Z_j\left(\frac{x}{A},\
\frac{y}{B} \right).
$$
is orthonormal on the ellipse  \cite{NLDP14,DMR26}.

\bigskip

\subsubsection*{Annular Zernike basis} 
In the case of a circular annulus, the unit disk is mapped to the domain
$$
O=\{(x,y)\in\mathbb{R}^2 \, : \, a^2\le x^2+y^2\le A^2 \}
= \{(r,\theta)\,:\, a \le r \le A,\ 0\le \theta <2\pi\},
$$
where $A$ is the outer radius and $a = hA$ is the radius of the central obscuration, with $0<h<1$. In this case, the diffeomorphism \cite{NLDP14} is the nonlinear radial mapping:
\begin{equation}
    \label{eq:map-annulus}
    (x,y)=\varphi(\rho,\phi)
=
\left( A ((1-h)\rho + h)\cos\phi,\;
A ((1-h)\rho + h)\sin\phi \right),
\end{equation}
which sends the radius $\rho\in[0,1]$ of the disk to
$$
r=A [(1-h)\rho + h] \in[a,A]
$$
while the angle remains unchanged $\theta=\phi$. 
Its inverse is given by
$$
(u,v)=\varphi^{-1}(x,y)=\left( \frac{r-hA}{A(1-h)}\cos\theta,\;
\frac{r-hA}{A(1-h)}\sin\theta \right).
$$
This transformation is nonlinear in the radial coordinate.
As a consequence, the Jacobian is not constant and is given  by
$$
J(r,\theta)=\frac{r-hA}{rA^2 (1-h)^2},
$$
which induces a non-Euclidean metric on the annulus.
This nonuniform Jacobian explains why the corresponding basis functions cease to be polynomials. There are two orthogonal bases defined by:
$$
O_j(r,\theta)
:=
Z_j\!\left( \frac{r-hA}{A(1-h)}\cos\theta,\;
\frac{r-hA}{A(1-h)}\sin\theta \right),
$$
$$
\tilde{O}_j(r,\theta)
:=
\sqrt{\frac{r-hA}{rA^2 (1-h)^2}}\, Z_j\!\left( \frac{r-hA}{A(1-h)}\cos\theta,\;
\frac{r-hA}{A(1-h)}\sin\theta \right).
$$
In the numerical experiments, we use the first one, see Section \ref{sec:experiments}. Both are orthonormal with respect to different weights.

\subsubsection*{Polygonal Zernike basis}
\begin{figure}
    \centering
    \includegraphics[width=0.4\linewidth]{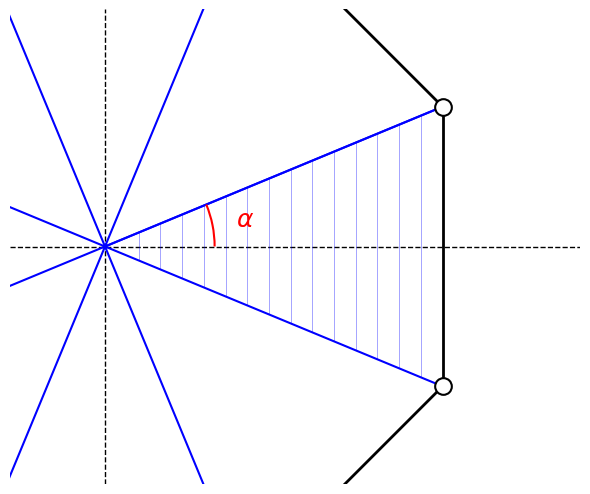}
    \caption{Graphic interpretation of $\alpha =\pi/p$ with $p=8$ sides.}
    \label{fig:alpha}
\end{figure}

Let us denote by $\Omega_p$ a polygon with $p$ sides and radius $1$ centered at the origin. Define $\alpha=\pi/p$ and assume that one of the radii of the polygon is oriented at an angle $\alpha$ from the positive $x$-axis, as
indicated in Figure~\ref{fig:alpha}.

To construct a mapping between the disk and the polygon, the authors in \cite{FLNP15} define the functions
\begin{align}
U_\alpha(\phi)&:=\phi-\left\lfloor\frac{\phi+\alpha}{2\alpha} \right\rfloor 2\alpha, \quad 0\le\phi<2\pi,
\\
\label{eq:function-R}
R_\alpha(\phi)&:=\frac{\cos\alpha}{\cos(U_\alpha(\phi))},
\end{align}
and the mapping $\varphi$ that transforms $D$  into the polygon $\Omega_p$ is:
\begin{equation}\label{eq:map_pol}
(x,y)=\varphi_\alpha(u,v)=(u \, R_\alpha(\arctan(v/u)), v \, R_\alpha(\arctan(v/u))).
\end{equation}
Alternatively, in polar coordinates,
\begin{equation}\label{eq:map_pol-polar}
    r(\rho,\phi)=\rho R_\alpha(\phi), \quad \theta(\rho, \phi)=\phi, \qquad 0\le \rho\le 1,\,\, 0\le \phi<2\pi.
\end{equation}

The inverse mapping is given by
$$
(u,v)=\varphi_\alpha^{-1}(x,y)=\left(\frac{x}{R_\alpha(\arctan(y/x))}, \frac{y}{R_\alpha(\arctan(y/x))} \right),
$$
or, in polar coordinates,
$$
\rho(r,\theta)= \frac{r}{R_\alpha(\theta)},  \quad \phi(r, \theta)=\theta.
$$

The Jacobian of the inverse transformation is $J_\alpha(\theta)=\frac{1}{R_\alpha(\theta)^2}$, and two orthogonal families are given by \cite[Section 3]{DMR26}
$$
K_j(r,\theta)=  Z_j\left(\frac{r}{R_\alpha(\theta)},\theta\right),
$$
$$
\tilde K_j(r,\theta)= \frac{1}{R_\alpha(\theta)} Z_j\left(\frac{r}{R_\alpha(\theta)},\theta\right).
$$
As  in the annulus, these new functions are no longer polynomials since the mapping is not linear.\\

The constructions presented above show how Zernike-based orthonormal systems can be transported from the unit disk to more general planar domains by means of suitable mappings. However, for approximation purposes, it is not sufficient to transfer only the basis functions; an appropriate choice of sampling nodes on the target domain is also required to ensure good numerical behavior. In the next section, we address this issue by introducing optimal node configurations obtained by mapping well-conditioned point sets from the disk.

\section{Optimal nodes}\label{sec:Opt_nodes}

In \cite{DMR26}, the authors use the mappings introduced above to construct sampling patterns that preserve favorable numerical conditioning. The central idea is to map the sampling points from the disk to the target domain $\Omega$ using the same diffeomorphism $\varphi$  employed to construct the orthogonal functions. That is, if
$$
S=\{s_1,s_2,\dots,s_M \}\subset D
$$
is a node configuration on the unit disk, where $M$ is the number of Zernike polynomials of order $\le m$, then 
$$
\varphi(S)= \{\varphi(s_1), \varphi(s_2), \dots, \varphi(s_M) \}
$$
is a node configuration on $\Omega$. 

\subsubsection*{Optimal nodes on the ellipse} In the ellipse case, the diffeomorphism is very simple, so the optimal set of interpolation nodes in the ellipse is given by
$$
\Sigma_m^{opt}(E)=\left\{\mathbf{p}_{\nu,\sigma}^{E} \,:\,\sigma=1,\dots,m_\nu, \, \nu=1,\dots, K  \right\},
$$
$$
\mathbf{p}_{\nu,\sigma}^{E}=\left(A\, \rho_\nu \cos\left(\frac{2\pi (\sigma-1)}{m_\nu} \right), B \, \rho_\nu \sin\left(\frac{2\pi (\sigma-1)}{m_\nu} \right)\right), 
$$
with $\rho_\nu$ given in \eqref{eq:rho_nu}; see the first column of Figure~\ref{fig:OCS-planar-surfaces}. 

\subsubsection*{Optimal nodes on the annulus}
The optimal set of interpolation nodes in the annulus is given by 
$$
\Sigma_m^{opt}(O)=\left\{\mathbf{p}_{\nu,\sigma}^{O} \,:\,\sigma=1,\dots,m_\nu, \, \nu=1,\dots, K  \right\},
$$
$$
\mathbf{p}_{\nu,\sigma}^{O}=\left(A ((1-h)\rho_\nu + h)  \cos\left(\frac{2\pi (\sigma-1)}{m_\nu} \right), \, A ((1-h)\rho_\nu + h)   \sin\left(\frac{2\pi (\sigma-1)}{m_\nu} \right)\right), 
$$
with $\rho_\nu$ given in \eqref{eq:rho_nu}; see the central column of Figure~\ref{fig:OCS-planar-surfaces}.

\subsubsection*{Optimal nodes on a polygon}
\begin{figure}
    \centering\begin{tabular}{ccc}
        \includegraphics[height=4cm]{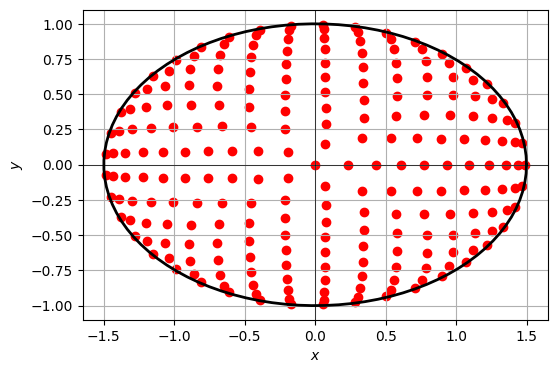} &
        \includegraphics[height=4cm]{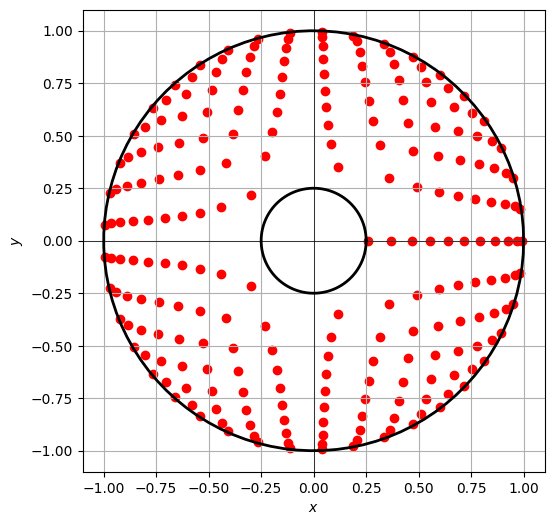} &
        \includegraphics[height=4cm]{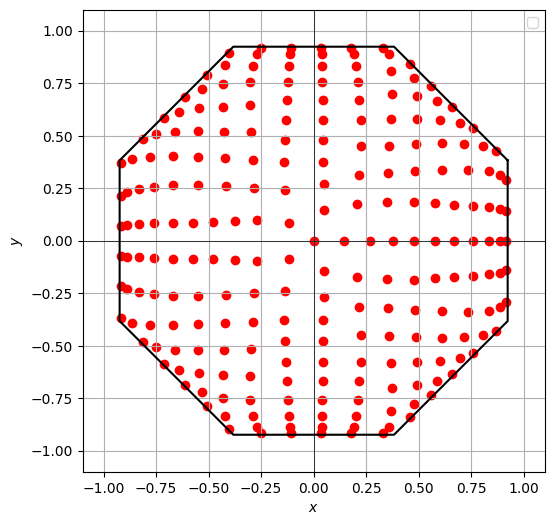} \\
        \includegraphics[height=4cm]{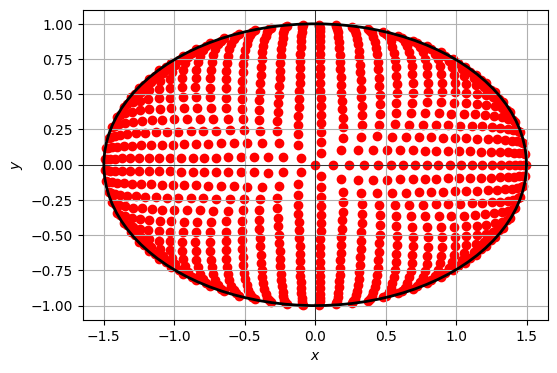} &
        \includegraphics[height=4cm]{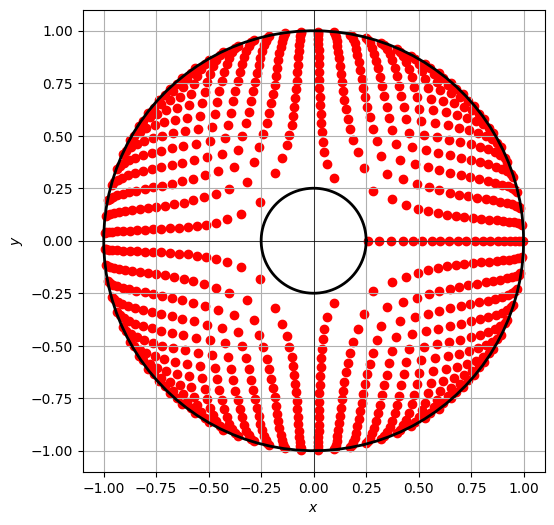} &
        \includegraphics[height=4cm]{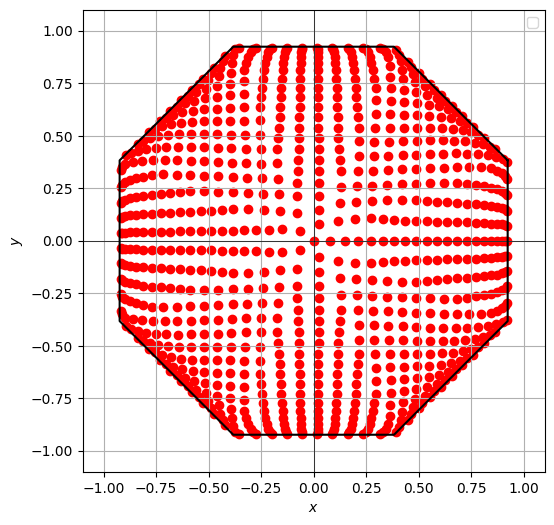} \\
    \end{tabular}
    \caption{From left to right: OCS nodes on the ellipse with semiaxes $1.5$ and $1$; the annulus with radii $1/4$ and $1$; and the regular polygon with $p=8$ sides; using $m=20$ (first row) and $m=40$ (second row).}
    \label{fig:OCS-planar-surfaces}
\end{figure}

As stated in \cite{DMR26},  the optimal nodes in the polygon $\Sigma_m^{opt}(\Omega_p)$ are obtained as the images, under the mapping $\varphi_\alpha$ defined in \eqref{eq:map_pol}, of the optimal points in the disk given in \eqref{eq:p_opt}; see the right column of Figure~\ref{fig:OCS-planar-surfaces}.\\

~\\
The construction of suitable node configurations on general domains provides the necessary framework for extending approximation techniques beyond the unit disk. In particular, once both the basis functions and the associated sampling points have been appropriately defined on $\Omega$, it becomes possible to generalize the interpolation–regression strategy introduced on the disk. This extension is developed in the next section, where a mixed interpolation–regression method is formulated for these transformed domains.

\section{Mixed interpolation-regression method}
\label{sec:mixed-IR-method}

Let us consider a function $f$ defined on one of the previous domains $\Omega$. We follow the same ideas as in Section \ref{Sec:IntReg_Disk}. 

Let $S_{n,n}\subset \Omega$ be a set of  $N=(n+1)^2$ points uniformly distributed over the domain. Although the authors of \cite{DMN25} consider a specific structured choice for this set (see \eqref{eq:Snn}), the proposed method remains applicable to any collection of points for which the minimization problem \eqref{eq:minimization-problem} is solvable and the set $\Sigma_m'$ defined in \eqref{eq:mock-optimal-points} contains exactly $M$ points. Accordingly, throughout this work,  $S_{n,n}$ is taken to be a realization of a sample of size $N$ from a uniform distribution supported on $\Omega$, and we assume that the values of the function $f$ are known at these sampling points.

In order to define the interpolation-regression operator, let us fix $m, \tilde{r}\in\mathbb{N}$ such that
$$
m<\tilde{r}, \quad \tilde{R}=\frac{(\tilde{r}+1)(\tilde{r}+2)}{2}<N, \quad M=\frac{(m+1)(m+2)}{2}.
$$

Furthermore, we consider the set of optimal points $\Sigma_m^{opt}(\Omega)$ defined in Section \ref{sec:Opt_nodes}. We explain how the interpolation-regression operator is constructed, following the approach introduced in Section~\ref{Sec:IntReg_Disk}, but adapted to these new domains.

\subsubsection*{Interpolation-regression on the ellipse}
First, we want to construct a polynomial of degree $\tilde{r}$ that interpolates the $M$ nodes closest to those in $\Sigma_m^{opt}(E)$, while the remaining nodes are used to further refine the approximation through a simultaneous regression.

The set of \textit{mock-optimal nodes} is 
\begin{equation}
\Sigma'_m(E) :=\left\{\mathbf{x}'_{\nu,\sigma} \, : \, \sigma=1,\dots,m_\nu, \, \nu=1,\dots, K  \right\} \subset S_{n,n},
\end{equation}
where $\mathbf{x}'_{\nu,\sigma}$ is the solution of the minimization problem
\begin{equation}
\min_{\mathbf{x}_{\eta,\kappa}\in S_{n,n}} \|\mathbf{x}_{\eta,\kappa} - \mathbf{p}_{\nu,\sigma} \|.
\end{equation}
We select these points as in Section \ref{Sec:IntReg_Disk}, and it is guaranteed that
$\# (\Sigma'_m(E))=M$.
Again, we use a single-index notation to simplify the exposition
\begin{align*}
\mathbf{x}_{\eta,\kappa} &\rightarrow \mathbf{x}_i,  \quad i=1,\dots,N,
\\
\mathbf{x}'_{\nu,\sigma} &\rightarrow \mathbf{x}'_j,  \quad j=1,\dots,M.
\end{align*}

Let $\left\{u_1,\dots,u_{\tilde{R}}\right\}=\left\{E_1,\dots,E_{\tilde{R}}\right\}$ be the basis of $\varphi(\mathbb{P}_{\tilde r} (D))$ given in Section~\ref{sec:extension-zernike}, and suppose that the set $S_{n,n}$ is such that the matrix $\mathcal{M}$  has rank $\tilde{R}$, where
$$
\mathcal{M}:=\left[
\begin{array}{cccc}
u_1(\mathbf{x}_1) & u_2(\mathbf{x}_1) & \dots & u_{\tilde{R}}(\mathbf{x}_1)
\\
u_1(\mathbf{x}_2) & u_2(\mathbf{x}_2) & \dots &u_{\tilde{R}}(\mathbf{x}_2)
\\
\vdots & \vdots & \ddots & \vdots
\\
u_1(\mathbf{x}_N) & u_2(\mathbf{x}_N) & \dots & u_{\tilde{R}}(\mathbf{x}_N)
\end{array}
\right] \in \mathbb{R}^{N\times \tilde{R}}.
$$

Suppose that the nodes of $\Sigma'_m(E)$ are sufficiently close to those of the set $\Sigma^{opt}_m(E)$ and, as a consequence, that the matrix
$$
\mathcal{C}:=\left[
\begin{array}{ccccc}
u_1(\mathbf{x}'_1) & u_2(\mathbf{x}'_1) & \dots & u_{\tilde{R}}(\mathbf{x}'_1)
\\
u_1(\mathbf{x}'_2) & u_2(\mathbf{x}'_2) & \dots & u_{\tilde{R}}(\mathbf{x}'_2)
\\
\vdots & \vdots & \ddots & \vdots
\\
u_1(\mathbf{x}'_M) & u_2(\mathbf{x}'_M) & \dots & u_{\tilde{R}}(\mathbf{x}'_M)
\end{array}
\right] \in \mathbb{R}^{M\times \tilde{R}},
$$
has rank  $M$.  We reorder the points in $S_{n,n}$ so that the mock-optimal nodes are the first $M$ ones, and thus $\mathcal{C}$ becomes the matrix formed by the first $M$ rows of $\mathcal{M}$. 

Let
$$
\bm{b}:=\big[f(\mathbf{x}_1), \dots, f(\mathbf{x}_N) \big]^T, \quad  \bm{d}:=\big[f(\mathbf{x}'_1), \dots, f(\mathbf{x}'_M) \big]^T.
$$
Then, an interpolation-regression operator is defined by
\begin{equation}\label{eq:Pi_el}
\hat{\Pi}_{\tilde{r},n}: f \longmapsto \hat{\Pi}_{\tilde{r},n}[f]:=\sum_{i=1}^R \hat{a}_i u_i \in \mathbb{P}_{\tilde{r}}(\mathbb{R}^2),
\end{equation}
where the vector $\hat{\bm{a}}:=\big[\hat{a}_1, \hat{a}_2, \dots, \hat{a}_{\tilde{R}} \big]^T$
is the solution to the least-squares problem 
$$
\min_{\bm{a}\in\mathbb{R}^R} \|\mathcal{M}\bm{a}-\bm{b}\|_2,\quad \text{subject to } \quad \mathcal{C}\bm{a}=\bm{d},
$$
and can be computed in the same way as in Section \ref{Sec:IntReg_Disk}, yielding analogous definitions for $\mathcal{R}_{11}$, $\mathcal{R}_{12}$, $\mathcal{Q}_{\mathcal{C}}$, $\mathcal{M}_1$ and $\mathcal{M}_2$, which will be used in Proposition \ref{pr:norm}.

\subsubsection*{Interpolation-regression on the annulus and the polygon}

In all other cases, when the Zernike polynomial basis is transported through the mapping $\varphi$, the resulting functions are no longer polynomials. Consequently, the operator $\hat{\Pi}_{\tilde{r},n}$ does not yield a polynomial of degree $\tilde{r}$, but rather a function given as a linear combination of the $\tilde{R}$ functions obtained by mapping the Zernike basis through $\varphi$.  

Let us denote by $\mathbb{F}_{\tilde{r}}(\mathbb{R}^2)$ the space generated by 
$$
\{u_j=\tilde{O}_j\,:\,j=1,\dots, \tilde{R}\}
$$
in the case of the annulus, and
$$
\{u_j=K_j\,:\,j=1,\dots, \tilde{R}\}
$$
in the case of the polygon. We can proceed as in the case of the ellipse, except that the interpolation-regression operator takes values in the space $\mathbb{F}_{\tilde{r}}$:
\begin{equation}\label{eq:Pi_An_Pol}
\hat{\Pi}_{\tilde{r},n}: f \longmapsto \hat{\Pi}_{\tilde{r},n}[f]:=\sum_{i=1}^{\tilde{R}} \hat{a}_i u_i \in \mathbb{F}_{\tilde{r}}(\mathbb{R}^2). 
\end{equation}

Theorems 1 and 2 in \cite{DMN25} admit extensions to the cases studied in this article, thanks to the construction of the interpolation and regression operators introduced here and their connection with the corresponding construction on the unit disk. Since the proofs are similar, we present only the statements, which are referred to as Proposition~\ref{pr:norm} and Proposition~\ref{pr:error}.

\begin{prop}\label{pr:norm}
The norms of the operators $\hat{\Pi}_{\tilde{r},n}$ given in \eqref{eq:Pi_el} and \eqref{eq:Pi_An_Pol} satisfy 
$$
\|\hat{\Pi}_{\tilde{r},n}\|_{\infty} \le \max_{i=1,\dots,\tilde{R}} \|u_i\|_{\infty} (K_1(n,m)+K_2(n,m)),
$$
where
\begin{align*}
    K_1(n,m)&=\|\mathcal{R}_{11}\|_1 \Big(M\|\mathcal{Q}_{\mathcal{C}}^T\|_1 + \|\mathcal{R}_{12}\|_1 K_2(n,m)\Big),
    \\
    K_2(n,m)&= \|(\mathcal{V}_1^T\mathcal{V}_1)^{-1} \mathcal{V}_1^T \|_1 \Big(N +M \|\mathcal{M}_1 \mathcal{R}_{11}^{-1} \mathcal{Q}_{\mathcal{C}}^T \|_1 \Big).
\end{align*}
\end{prop}

For the second result, let us denote by
$$
\varepsilon_{\tilde{r}}(f) =\|f- p^{\star}_{\tilde{r}}[f] \|_{\infty},
$$
where $p^{\star}_{\tilde{r}}[f]$ is a best uniform approximation in the space $\mathbb{P}_{\tilde{r}}(\mathbb{R}^2)$ or $\mathbb{F}_{\tilde{r}}(\mathbb{R}^2)$ depending on whether the domain is an  ellipse, an annulus, or a polygon. Observe that, in the case where the space consists of polynomials, this best approximation is unique; however, this may not be the case when the space  is $\mathbb{F}_{\tilde{r}}(\mathbb{R}^2)$. 

\begin{prop}\label{pr:error}
Let $f$ be a function defined on an ellipse, an annulus, or a polygon. Then,
$$
\|f- \hat{\Pi}_{\tilde{r},n}[f] \|_{\infty} \le \Big( 1+\max_{i=1,\dots,\tilde{R}} \|u_i\|_{\infty} \big(K_1(n,m)+K_2(n,m)\big)\Big) \varepsilon_{\tilde{r}}(f).
$$
\end{prop}

\begin{figure}
    \centering\begin{tabular}{ccc}
        \includegraphics[height=4.5cm]{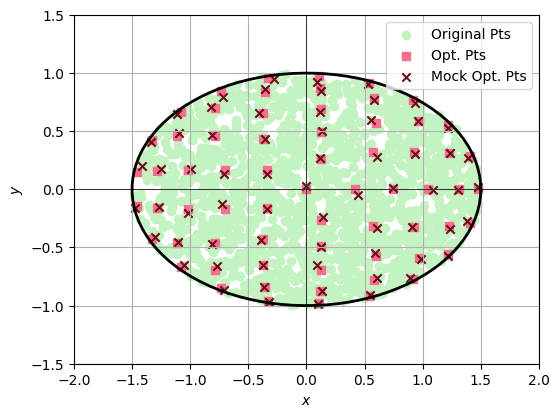} &
        \includegraphics[height=4.5cm]{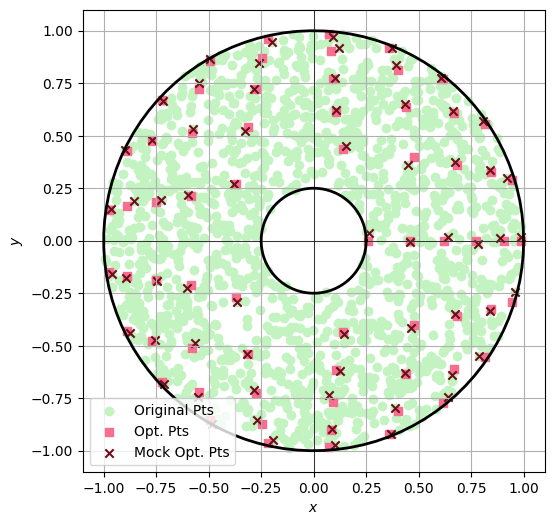} &
        \includegraphics[height=4.5cm]{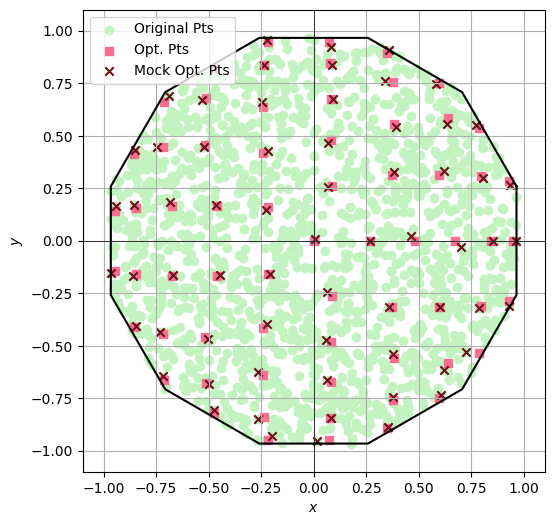}
    \end{tabular}
    \caption{Points from $S_{40,40}$ (green circles), OCS optimal points for $m=10$ (pink squares) and mock-optimal points in the ellipse with semiaxes $1.5$ and $1$; the annulus with radii $1/4$, $1$ and the regular polygon with $p=12$ sides.}
    \label{fig:mock-optimal-planar-surfaces}
\end{figure}

\noindent\textbf{Remark.}
An alternative approach would consist of transporting the approximation problem to the unit disk via the inverse of the diffeomorphism. More precisely, given \(f:\Omega \to \mathbb{R}\) and a diffeomorphism (a.e.) \(\varphi:D \to \Omega\), one may consider the composition \(f \circ \varphi\), approximate it on the disk by \(\hat{\Pi}^D_{\tilde{r},n}[f \circ \varphi]\), and then map the result back to \(\Omega\). 

However, our numerical experiments indicate that the resulting approximation is significantly less accurate than the direct construction \(\hat{\Pi}^{\Omega}_{\tilde{r},n}[f]\). This suggests that the interpolation–regression procedure is not invariant under such transformations and that performing the approximation directly on the target domain yields substantially better results.

This behavior can be attributed to the distortion introduced by the mapping \(\varphi\), which affects the distribution of nodes and the stability of the approximation.\\

The interpolation–regression framework developed in the previous section provides accurate approximations of functions defined on general domains $\Omega$ from discrete data. Beyond approximation, this construction can also be exploited to address related numerical tasks. In particular, it naturally leads to efficient strategies for numerical integration when function values are only available at scattered points. In the next section, we show how to derive cubature formulas based on this approach.

\section{Cubature formulas through an interpolation-regression method}
\label{sec:cubature}

Let us consider a Riemann-Stieltjes integrable function $f:\Omega \longrightarrow \mathbb{R}$ defined on one of the domains $\Omega$ considered in this paper. If we seek to approximate the integral 
$$
\int_{\Omega} f(\mathbf x) d\mathbf x,
$$
but the values of $f$ are available only on a randomly distributed set of $N=(n+1)^2$ points, denoted by $S_{n,n}$, then this integral can be approximated by integrating $\hat{\Pi}_{\tilde r,n}[f]$. That is,
\begin{equation}
    \label{eq:approx-integral}
    \int_{\Omega} f(\mathbf x)  d\mathbf{x} \approx \int_{\Omega} \hat{\Pi}_{\tilde r,n}[f](\mathbf x)  d\mathbf x.
\end{equation}

In this context, two potential approaches can be considered for evaluating the right-hand side integral in \eqref{eq:approx-integral}. The most straightforward method involves the direct computation of the integral, since the integrand is explicitly known on $\Omega$. For the specific case of an ellipse, the approximation $\hat{\Pi}_{\tilde r,n}[f](\mathbf x)$ is a polynomial of degree $\tilde r$ on $E$. Consequently, the procedure reduces to the integration of a polynomial over an elliptical domain, a task that can be readily accomplished either symbolically or numerically. Conversely, since $\hat{\Pi}_{\tilde r,n}[f](\mathbf x)$ is not a polynomial  for polygonal or annular geometries, the integration process becomes significantly more complex and computationally expensive.

Nevertheless, several computational tests have been conducted, yielding highly accurate approximations within short execution times when the integral is evaluated over an ellipse. As anticipated, these favorable results do not extend to annular and polygonal domains; even for low values of $\tilde r$, execution times exceed 30 minutes, ultimately producing unsatisfactory results.

Another viable approach is to employ a cubature formula, thereby approximating the integral as a weighted sum of the evaluations of $\hat{\Pi}_{\tilde r,n}[f](\mathbf x)$ at a suitable set of integration nodes. The computational cost of this method is substantially lower than that of integrating the approximation itself. This cost depends primarily on the value of $\tilde r$, since both the structural complexity of the operator and its evaluation time increase with larger values of $\tilde r$.

In the remainder of this section, we explain how to construct cubature formulas over an ellipse, an annulus, and a polygon, using existing formulas on the unit disk as a starting point.

\subsection{Cubature formulas using Gaussian nodes}

In the unit disk $D$, let us consider a cubature formula with algebraic degree of exactness $q \geq 0$. This requires the existence of a set of $Q = \frac{1}{2}(q+1)(q+2)$ Gaussian cubature nodes $\{\xi_1, \dots, \xi_Q\} \subset D$ and their associated weights $\{\omega_1, \dots, \omega_Q\}$, such that an integral over the disk may be approximated via the cubature rule
\begin{equation}
    \label{eq:cubature-disk}
    \int_D f(x,y)\, dx\, dy \approx \sum_{i=1}^Q \omega_i f(\xi_i),
\end{equation}
where the approximation becomes an exact equality whenever $f(x, y)$ is a polynomial of degree at most $q$ on $D$. In polar coordinates, \eqref{eq:cubature-disk} reads as
\begin{equation}
    \label{eq:cubature-disk-polar}
    \int_D f(x,y)\, dx\, dy = \int_{0}^{2\pi}\int_{0}^{1} \tilde f (\rho,\phi) \rho \,d\rho \, d\phi \approx \sum_{i=1}^Q \omega_i \tilde f(\rho_i, \phi_i),
\end{equation}
where $\tilde f (\rho,\phi) = f(\rho \cos\phi,\rho\sin\phi)$ represents the function $f$ in polar coordinates, and $(\rho_i,\phi_i)$ denote the polar coordinates of the Gaussian nodes $\xi_i$ for $i=1,\dots,Q$.

By employing an appropriate change of variables, we can construct cubature formulas for an ellipse, an annulus, and a polygon, using the Gaussian nodes and their associated weights from the unit disk. To this end, let $\{\zeta_1,\dots,\zeta_Q\}$ denote the set of cubature nodes on the corresponding domain $\Omega$, and let $\{w_1,\dots,w_Q\}$ be their associated weights. However, since the original function $f$ may not be known at these specific nodes, we evaluate the interpolation-regression operator at these locations to approximate the function values. Consequently, the integral can be approximated as
\begin{equation}
    \label{eq:cubature}
    \int_\Omega f(x,y)\, dx\, dy \approx \sum_{i=1}^Q f(\zeta_i)w_i \approx \sum_{i=1}^Q \hat\Pi_{\tilde r, n}[f](\zeta_i) w_i .
\end{equation}

In what follows, we describe how to derive the cubature nodes and their associated weights for each specific geometry.

\subsubsection{Cubature formulas on an ellipse}

Let $f:E\longrightarrow \mathbb R$ be a continuous function defined on an ellipse $E$ with semiaxes $A$ and $B$. By employing the mapping $\varphi:D\longrightarrow E$ given by $\varphi(u,v)=(Au,Bv)$, as introduced in \eqref{eq:map-ellipse}, whose Jacobian determinant is $|J_\varphi|=AB$, the integral can be expressed over the unit disk as
$$
\int_E f(x,y) \, dx\, dy = \int_D f(\varphi(u,v))  AB \, du\, dv.
$$
Applying the disk cubature rule \eqref{eq:cubature-disk} to the right-hand side yields
$$
\int_E f(x,y)\,  dx\, dy \approx \sum_{i=1}^Q f(\varphi(\xi_i))\,AB \,\omega_i.
$$
Consequently, the integral over the elliptical domain is approximated via the cubature scheme \eqref{eq:cubature}, where the target nodes $\{\zeta_1,\dots,\zeta_Q\}$ and their corresponding weights $\{w_1,\dots,w_Q\}$ are explicitly defined by
\begin{align}
\label{eq:cubature-ellipse}
    \zeta_i &= (Au_i, B v_i), & w_i &= AB\, \omega_i, \qquad i=1,\dots,Q,
\end{align}
where $\xi_i=(u_i,v_i)$ are the Cartesian coordinates of the Gaussian nodes on $D$. This cubature formula remains exact for any function $f \in \mathbb{P}_{q}(E)$, the space of bivariate polynomials of degree at most $q$ on the ellipse.

\subsubsection{Cubature formulas on an annulus}

Let $f:O\longrightarrow\mathbb R$ be a continuous function defined on a circular annulus $O$ with outer radius $A$ and inner radius $a=hA$ ($0<h<1$). Let $\tilde f(r, \theta)=f(r\cos\theta,r\sin\theta)$ denote the evaluation of $f$ at a point $(x,y)\in O$ expressed in polar coordinates $(r,\theta)$. By recalling the mapping $\varphi:D\longrightarrow O$ given in \eqref{eq:map-annulus}, its action in polar coordinates for any $(\rho,\phi)\in D$  yields $\varphi(\rho,\phi)=(A((1-h)\rho+h),\phi)$. Since the Jacobian determinant of this transformation is $|J_{\varphi}|=A-a$, the integral over the annulus can be reformulated as follows:
\begin{align*}
    \int_O f(x,y)\, dx\, dy &= \int_{0}^{2\pi}\int_{a}^{A}\tilde f(r,\theta)\, r \, dr \, d\theta\\ 
    &= \int_{0}^{2\pi}\int_{0}^{1}\tilde f(\varphi(\rho,\phi)) \, A((1-h)\rho+h)\, (A-a) \, d\rho \, d\phi\\
    &= \int_{0}^{2\pi}\int_{0}^{1}\tilde f(\varphi(\rho,\phi))\, \dfrac{A((1-h)\rho+h)\, (A-a)}\rho \, \rho \, d\rho \, d\phi.
\end{align*}
Applying the polar-form disk cubature rule \eqref{eq:cubature-disk-polar} to this last expression yields
$$
\int_O f(x,y)\, dx\, dy \approx \sum_{i=1}^Q \tilde f(\varphi(\rho_i,\phi_i)) \, \dfrac{A((1-h)\rho_i+h)\, (A-a)}{\rho_i} \omega_i.
$$
Therefore, the integral can be approximated via the cubature scheme \eqref{eq:cubature}, where the target nodes $\{\zeta_1,\dots,\zeta_Q\}$ and their associated weights $\{w_1,\dots,w_Q\}$ are given by
\begin{align}
\label{eq:cubature-annulus}
    \zeta_i &= (A((1-h)\rho_i+h)\cos\phi_i, A((1-h)\rho_i+h)\sin\phi_i),  & w_i &= \dfrac{A((1-h)\rho_i+h)\, (A-a)}{\rho_i} \omega_i,
\end{align}
for $i=1,\dots,Q$, where $(\rho_i,\phi_i)$ are the polar coordinates of the original disk nodes. This integration rule is exact for any function belonging to $\mathbb F_{q}(O)$, the space of bivariate functions on the annulus obtained by transporting the polynomial space $\mathbb{P}_q(D)$ through the mapping $\varphi$ in \eqref{eq:map-annulus}.

\subsubsection{Cubature formulas on a polygon}

Given a continuous function $f:\Omega_p \longrightarrow \mathbb R$ defined on a regular polygon $\Omega_p$ with $p\geq 3$ sides, let $\tilde f(r, \theta)=f(r\cos\theta,r\sin\theta)$ denote its representation in polar coordinates. We recall the mapping $\varphi:D\longrightarrow \Omega_p$ expressed in polar coordinates in \eqref{eq:map_pol-polar}, whose Jacobian determinant is given by $|J_{\varphi}|=R_\alpha(\phi)$, where the function $R_\alpha(\phi)$ is defined in \eqref{eq:function-R}. By applying the change of variables defined by $\varphi$, the integral over the polygon can be written as
\begin{align*}
    \int_{\Omega_p} f(x,y)\, dx\, dy &= \int_{0}^{2\pi}\int_{0}^{R_\alpha(\theta)}\tilde f(r,\theta) \, r \, dr \, d\theta  \\
    &= \int_{0}^{2\pi} \int_{0}^{1} \tilde f(\varphi(\rho,\phi)) \,  (R_{\alpha}(\phi))^2 \rho\, d\rho \, d\phi .
\end{align*}
This latter integral can be approximated by employing the polar-form disk cubature rule \eqref{eq:cubature-disk-polar}, which yields
$$
\int_{\Omega_p} f(x,y)\, dx\, dy \approx \sum_{i=1}^Q \tilde f(\varphi(\rho_i,\phi_i)) (R_{\alpha}(\phi_i))^2\omega_i. 
$$
In this way, we obtain a cubature scheme \eqref{eq:cubature} in the polygonal domain, where the target nodes $\{\zeta_1,\dots,\zeta_Q\}$ and their associated weights $\{w_1,\dots,w_Q\}$ are given by 
\begin{align}
\label{eq:cubature-polygon}
    \zeta_i &= \varphi(\rho_i,\phi_i), & w_i &= (R_\alpha(\phi_i))^2 \omega_i, \qquad i=1,\dots,Q,
\end{align}
with $(\rho_i,\phi_i)$ denoting the polar coordinates of the original disk nodes. This integration rule is exact for any function belonging to $\mathbb F_{q}(\Omega_p)$, the space of bivariate functions on the polygon obtained by transporting the polynomial space $\mathbb{P}_q(D)$ through the diffeomorphism $\varphi$ in \eqref{eq:map_pol}.\\

In summary, we have established a systematic framework for constructing cubature formulas across various geometries by mapping the Gaussian nodes and weights of the unit disk onto the target domains. To evaluate the efficiency and accuracy of the proposed integration methods, extensive computational tests have been conducted. A detailed analysis of their performance, execution times, and accuracy is provided in Section~\ref{sec:num-exp-cubature}.

In the following section, we present numerical experiments  to validate our theoretical results. We outline the various parameters used during the tests, introduce the metrics used to evaluate the performance of the proposed methods, and discuss the results obtained in detail.

\section{Numerical experiments} \label{sec:experiments}

To illustrate the performance of these interpolation-regression methods, numerous experiments were conducted by varying the parameters $m$ and $\tilde r$ to approximate different test functions on the previously described planar surfaces. In these simulations, we considered an ellipse with semiaxes $1.5$ (horizontal) and $1$ (vertical), a circular annulus with inner radius $1/4$ and outer radius $1$, and a regular polygon with $p=12$ sides. These geometric parameters can be easily modified to accommodate other configurations.

For each domain $\Omega$, a randomly distributed sample of points $S_{n,n}$ was generated by fixing $n=100$, resulting in a data set of size $N = \#S_{100,100} = (101)^2 = 10201$. We assume that the values of the target function are available exclusively at these $N$ points, representing the finite and discrete information at our disposal. The choice of $n=100$ ensures that the functions are sampled over a sufficiently large data set, potentially yielding a higher approximation accuracy, albeit at the expense of significantly increased computational time and memory requirements. All computations were carried out on the \emph{ALBAICIN} supercomputer at the University of Granada (\url{https://supercomputacion.ugr.es/}).

We have employed the previously described interpolation-regression methodology for the following test functions
$$
\begin{array}{lll}
f_1(x,y)=\sin(xy), & f_2(x,y)=e^{-xy}, & f_3(x,y)= e^{-(x^2+y^2)},
\\[10pt]
f_4(x,y)=\frac{1}{x^2+y^2+1}, & f_5(x,y)=\cos(x)\sin(y), & f_6(x,y)=\ln(x^2+y^2+1).
\end{array}
$$

\begin{figure}
    \centering\begin{tabular}{ccc}
        \includegraphics[height=4cm]{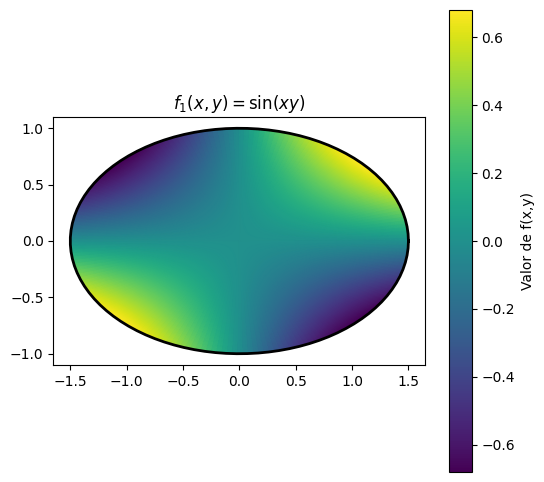} &
        \includegraphics[height=4cm]{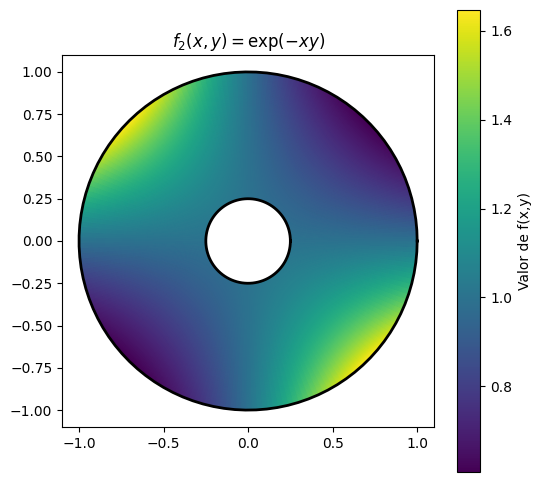} &
        \includegraphics[height=4cm]{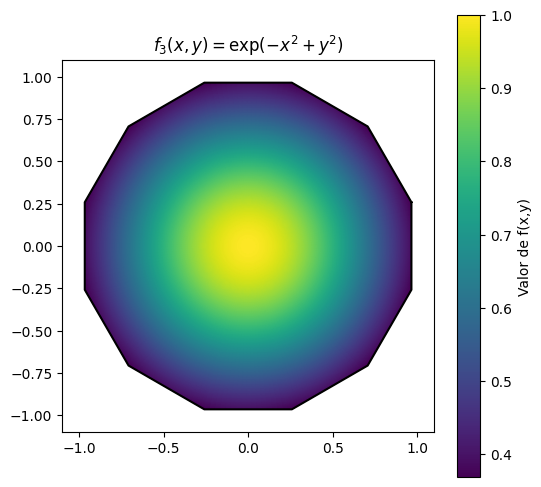} \\
        \includegraphics[height=4cm]{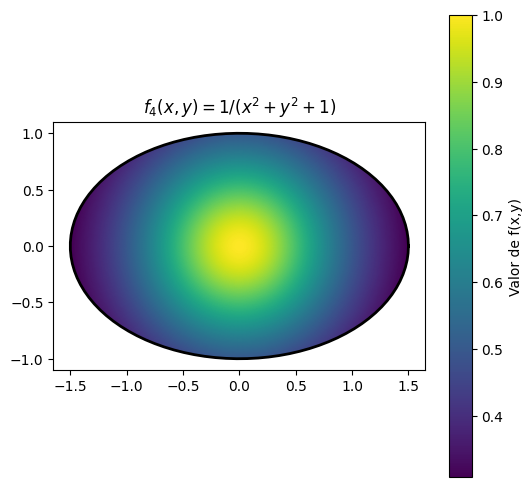} &
        \includegraphics[height=4cm]{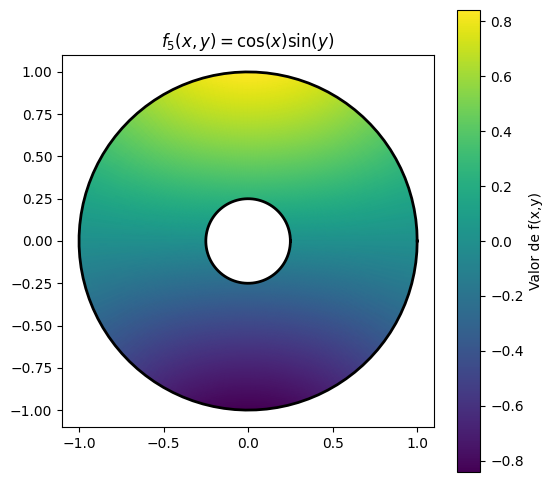} &
        \includegraphics[height=4cm]{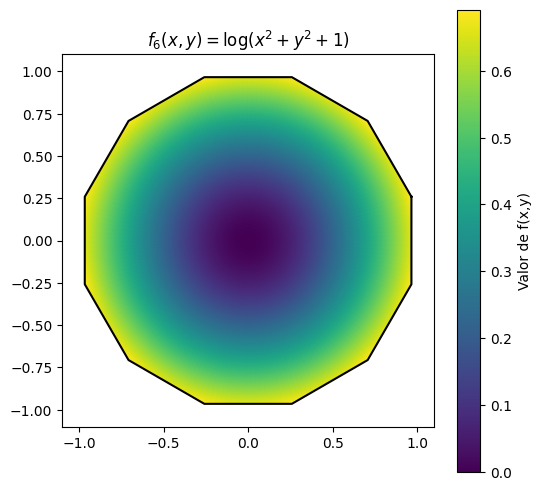} \\
    \end{tabular}
    \caption{Heatmap plots of the functions $f_1, \dots, f_6$. $f_1$ and $f_4$ are plotted on the ellipse, $f_2$, and $f_5$ on the annulus and $f_3$ and $f_6$ on the polygon.}
    \label{fig:plot-functions}
\end{figure}

See Figure~\ref{fig:plot-functions} for plots of these functions.

For each test function $f \in \{f_1, \dots, f_6\}$ defined on the corresponding domains, we construct the interpolation-regression operator for various combinations of the parameters $m$ and $\tilde{r}$. To evaluate the accuracy of the operator $\hat\Pi_{\tilde r,n}[f]$, a validation set consisting of $5000$ test points, denoted by $\mathcal{T}$, is independently generated. Consequently, for each considered value of $m$ and $\tilde{r}$, we compute the following performance metrics:
\begin{itemize}
    \item \textbf{Mean Squared Error (MSE):} The average of the squared local errors, defined as
    $$
    \mathrm{MSE} = \frac{1}{\#\mathcal{T}} \sum_{\mathbf{x} \in \mathcal{T}} \left(f(\mathbf{x}) - \hat\Pi_{\tilde r,n}[f](\mathbf{x})\right)^2. 
    $$
    \item \textbf{Maximum Absolute Error (MaxAE):} The maximum absolute local error, defined as
    $$
    \mathrm{MaxAE} = \max_{\mathbf{x} \in \mathcal{T}} \left\{ \left| f(\mathbf{x}) - \hat\Pi_{\tilde r,n}[f](\mathbf{x}) \right| \right\}.
    $$
    \item \textbf{Mean Relative Error (MRE):} The average of the relative local errors, defined as
    $$
    \mathrm{MRE} = \frac{1}{\#\mathcal{T}} \sum_{\mathbf{x} \in \mathcal{T}} \frac{\left| f(\mathbf{x}) - \hat\Pi_{\tilde r,n}[f](\mathbf{x}) \right|}{\left| f(\mathbf{x}) \right|}.
    $$
    \item \textbf{Maximum Relative Error (MaxRE):} The maximum relative local error, defined as
    $$
    \mathrm{MaxRE} = \max_{\mathbf{x} \in \mathcal{T}} \left\{ \frac{\left| f(\mathbf{x}) - \hat\Pi_{\tilde r,n}[f](\mathbf{x}) \right|}{\left| f(\mathbf{x}) \right|} \right\}.
    $$
\end{itemize}

While MSE and MRE provide a global measure of the approximation accuracy across the entire domain, MaxAE and MaxRE highlight the largest error. Consequently, a low MSE or MRE indicates a generally accurate approximation; however, local discrepancies in specific regions may still occur, leading to high MaxAE or MaxRE values.

Furthermore, evaluating the efficiency of the method is also essential. Although increasing the values of $m$ or $\tilde r$ may yield a more accurate approximation, larger values also introduce additional structural complexity, making the operator more computationally expensive to evaluate. To quantify this aspect, we introduce the execution time metric (\text{ExTime}), which measures the time (in seconds) required to evaluate the approximant $\hat\Pi_{\tilde r,n}[f]$ at the $5000$ test points within $\mathcal{T}$.\\

First, both parameters $m$ -- the dimension of the interpolation part -- and $\tilde r$ -- the dimension of the regression part -- are varied, with the aim of empirically identifying the optimal value of $m$, while $n$ remains fixed.

\subsection{Varying the dimension of the interpolation part}

First, once $n=100$ has been fixed, we perform calculations for several values of $m$ -- the degree of the interpolation part -- and, following the indications in \cite{DMN25}, set $\tilde{r} = m+\lfloor\sqrt m \rfloor$. In this way, we can empirically identify a potentially optimal value of $m$ in terms of both the accuracy of the approximant and its efficiency.

\begin{table}[!ht]
    \centering
    \begin{tabular}{ccccccccc}\toprule
        $m$ & $\tilde r$ & $M$ & $\tilde R$ & \textbf{MSE} & \textbf{MaxAE} & \textbf{MRE} & \textbf{MaxRE} & \textbf{ExTime} \\ \toprule
        5 & 7 & 21 & 36 & 0.000583713 & 0.00400818 & 0.000320182 & 0.00577733 & 2.85298 \\ \hline
        10 & 13 & 66 & 105 & 1.06712e-07 & 9.14e-07 & 4.97149e-08 & 1.27423e-06 & 13.2923 \\ \hline
        15 & 18 & 136 & 190 & 7.20522e-12 & 6.77112e-11 & 3.50548e-12 & 1.20842e-10 & 37.1076 \\ \hline
        20 & 24 & 231 & 325 & 9.66913e-16 & 1.23235e-14 & 7.38556e-16 & 1.86074e-14 & 93.8597 \\ \hline
        25 & 30 & 351 & 496 & 2.115e-15 & 5.4623e-14 & 1.52365e-15 & 3.34204e-14 & 202.61 \\ \hline
        30 & 35 & 496 & 666 & 3.53692e-14 & 9.7411e-13 & 1.37903e-14 & 1.51396e-12 & 352.077 \\ \hline
        35 & 40 & 666 & 861 & 7.01799e-14 & 3.1517e-12 & 1.03332e-14 & 1.94669e-12 & 571.362 \\ \hline
        40 & 46 & 861 & 1128 & 6.69038e-14 & 2.16138e-12 & 1.12521e-14 & 2.91347e-12 & 945.089 \\ \hline
        45 & 51 & 1081 & 1378 & 4.96089e-13 & 2.61772e-11 & 6.00625e-14 & 3.92902e-11 & 1392.89 \\ \bottomrule
    \end{tabular}
    \caption{Results for $f_2$ on the ellipse}
    \label{tab:results-f2-ellipse}
\end{table}

\begin{table}[!ht]
    \centering
    \begin{tabular}{ccccccccc}
    \toprule
    $m$ & $\tilde r$ & $M$ & $\tilde R$ & \textbf{MSE} & \textbf{MaxAE} & \textbf{MRE} & \textbf{MaxRE} & \textbf{ExTime} \\ \toprule
        5 & 7 & 21 & 36 & 0.0169161 & 0.0588409 & 0.0125955 & 0.0759125 & 2.6541 \\ \hline
        10 & 13 & 66 & 105 & 0.00762156 & 0.036897 & 0.00394671 & 0.0382813 & 12.7661 \\ \hline
        15 & 18 & 136 & 190 & 0.00696165 & 0.0423871 & 0.00444445 & 0.0406742 & 35.8995 \\ \hline
        20 & 24 & 231 & 325 & 0.00487832 & 0.0337166 & 0.0017697 & 0.034984 & 92.533 \\ \hline
        25 & 30 & 351 & 496 & 0.00472315 & 0.0385835 & 0.00253984 & 0.0400453 & 197.908 \\ \hline
        30 & 35 & 496 & 666 & 0.00380984 & 0.033021 & 0.00125955 & 0.0341393 & 346.274 \\ \hline
        35 & 40 & 666 & 861 & 0.00439778 & 0.122885 & 0.00210113 & 0.196257 & 566.606 \\ \hline
        40 & 46 & 861 & 1128 & 0.00339229 & 0.0323305 & 0.0009733 & 0.0334337 & 957.921 \\ \hline
        45 & 51 & 1081 & 1378 & 0.0242035 & 1.46888 & 0.00369996 & 1.73898 & 1412.98 \\ \bottomrule
    \end{tabular}
    \caption{Results for $f_2$ on the annulus}
    \label{tab:results-f2-annulus}
\end{table}

\begin{table}[!ht]
    \centering
    \begin{tabular}{ccccccccc}
    \toprule
        $m$ & $\tilde r$ & $M$ & $\tilde R$ & \textbf{MSE} & \textbf{MaxAE} & \textbf{MRE} & \textbf{MaxRE} & \textbf{ExTime} \\ \toprule
         5 & 7 & 21 & 36 & 0.0054574 & 0.0358493 & 0.00332016 & 0.0311524 & 3.2935 \\ \hline
        10 & 13 & 66 & 105 & 0.0044551 & 0.0280493 & 0.00293695 & 0.0272348 & 14.5202 \\ \hline
        15 & 18 & 136 & 190 & 0.00217561 & 0.0193766 & 0.00148307 & 0.0126075 & 37.2147 \\ \hline
        20 & 24 & 231 & 325 & 0.0016313 & 0.013676 & 0.00108684 & 0.012305 & 93.9419 \\ \hline
        25 & 30 & 351 & 496 & 0.00128992 & 0.0128795 & 0.000795406 & 0.00797253 & 202.5 \\ \hline
        30 & 35 & 496 & 666 & 0.00100655 & 0.0161114 & 0.000637465 & 0.0121597 & 351.226 \\ \hline
        35 & 40 & 666 & 861 & 0.000858515 & 0.00816177 & 0.000529015 & 0.010715 & 570.218 \\ \hline
        40 & 46 & 861 & 1128 & 0.001024 & 0.0268464 & 0.000544163 & 0.0433482 & 955.875 \\ \hline
        45 & 51 & 1081 & 1378 & 0.00166466 & 0.0430288 & 0.000739708 & 0.0388824 & 1403.85 \\ \bottomrule
    \end{tabular}
    \caption{Results for $f_2$ on the polygon}
    \label{tab:results-f2-polygon}
\end{table}

\begin{figure}
    \centering\begin{tabular}{ccc}
        \includegraphics[width=0.3\textwidth]{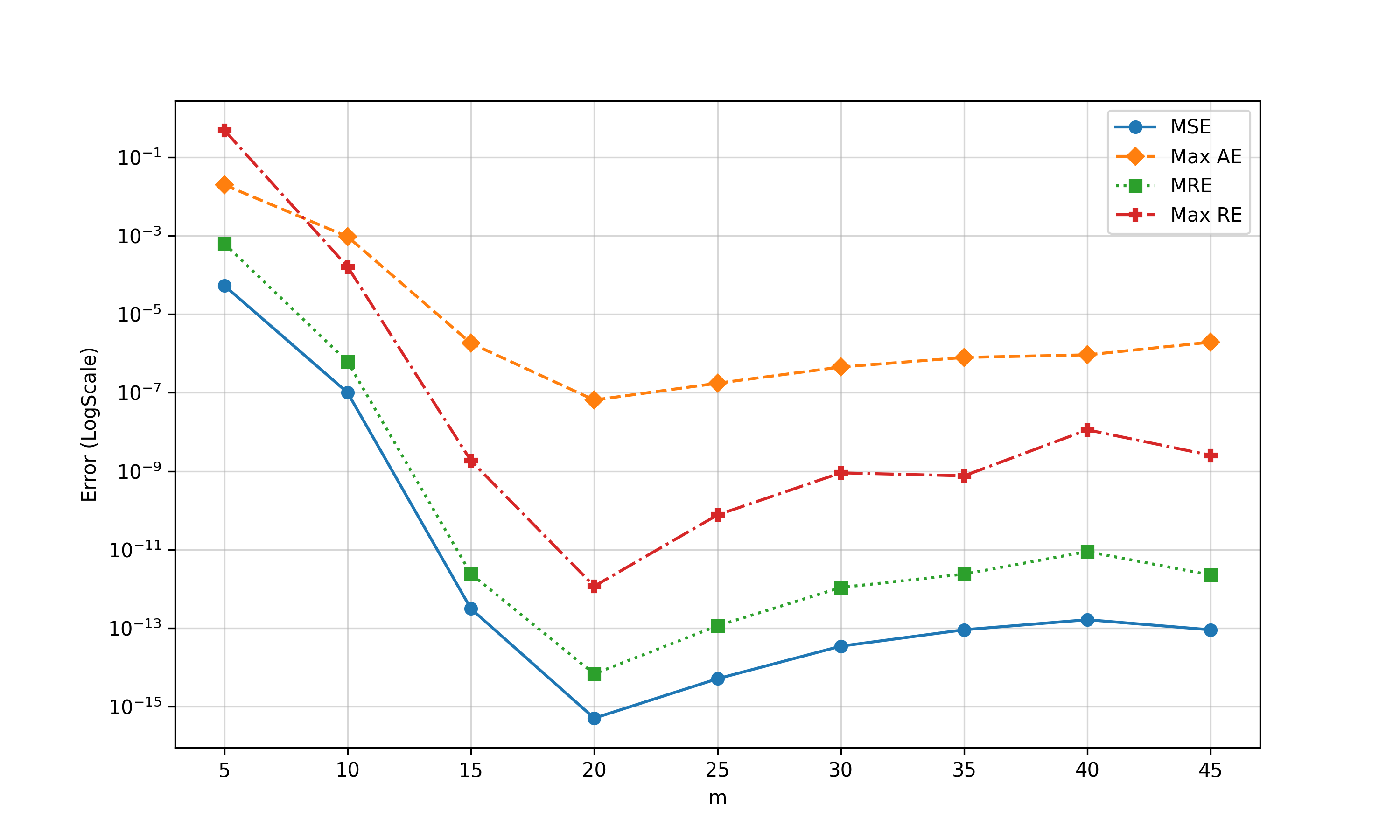} & \includegraphics[width=0.3\textwidth]{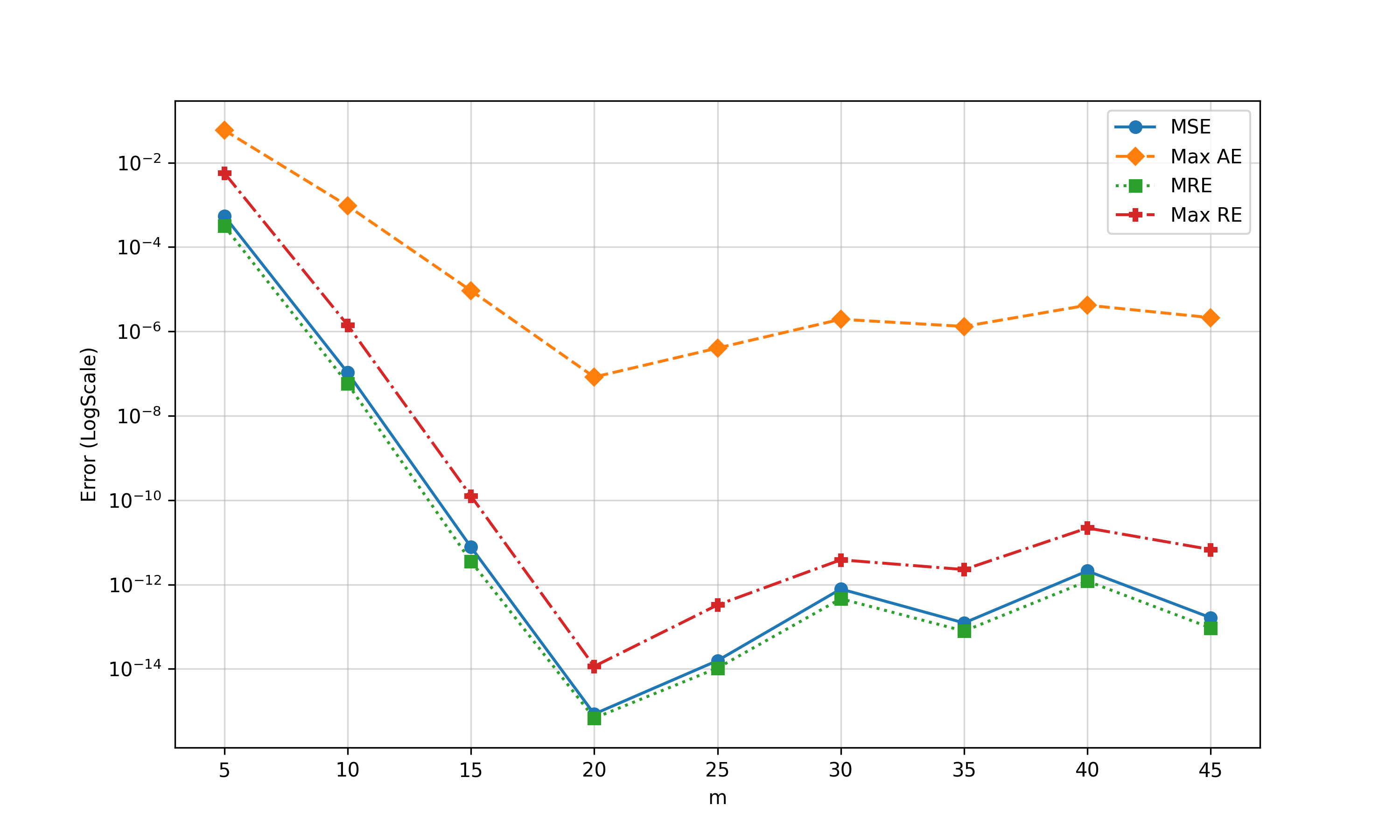} &
        \includegraphics[width=0.3\textwidth]{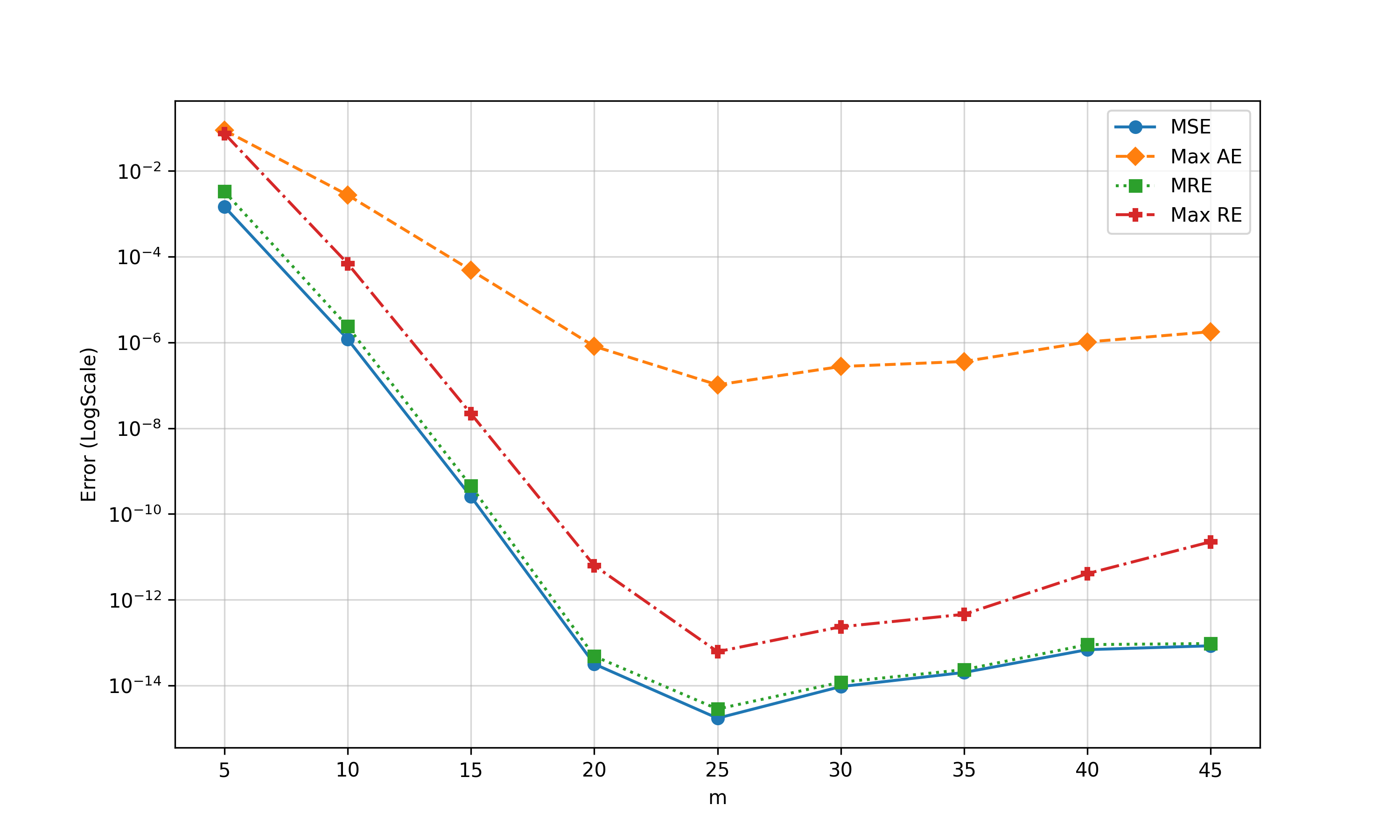} \\ \includegraphics[width=0.3\textwidth]{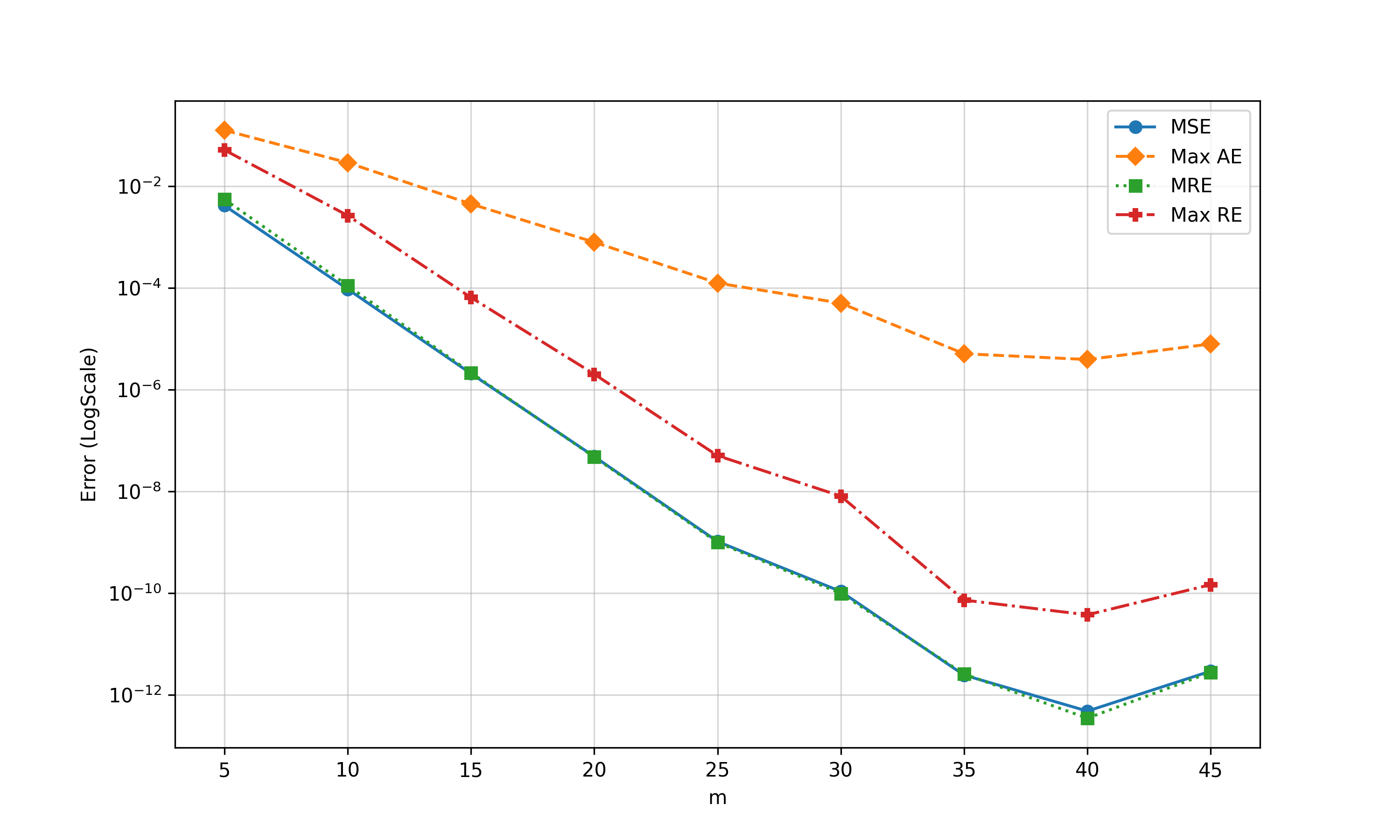} &
        \includegraphics[width=0.3\textwidth]{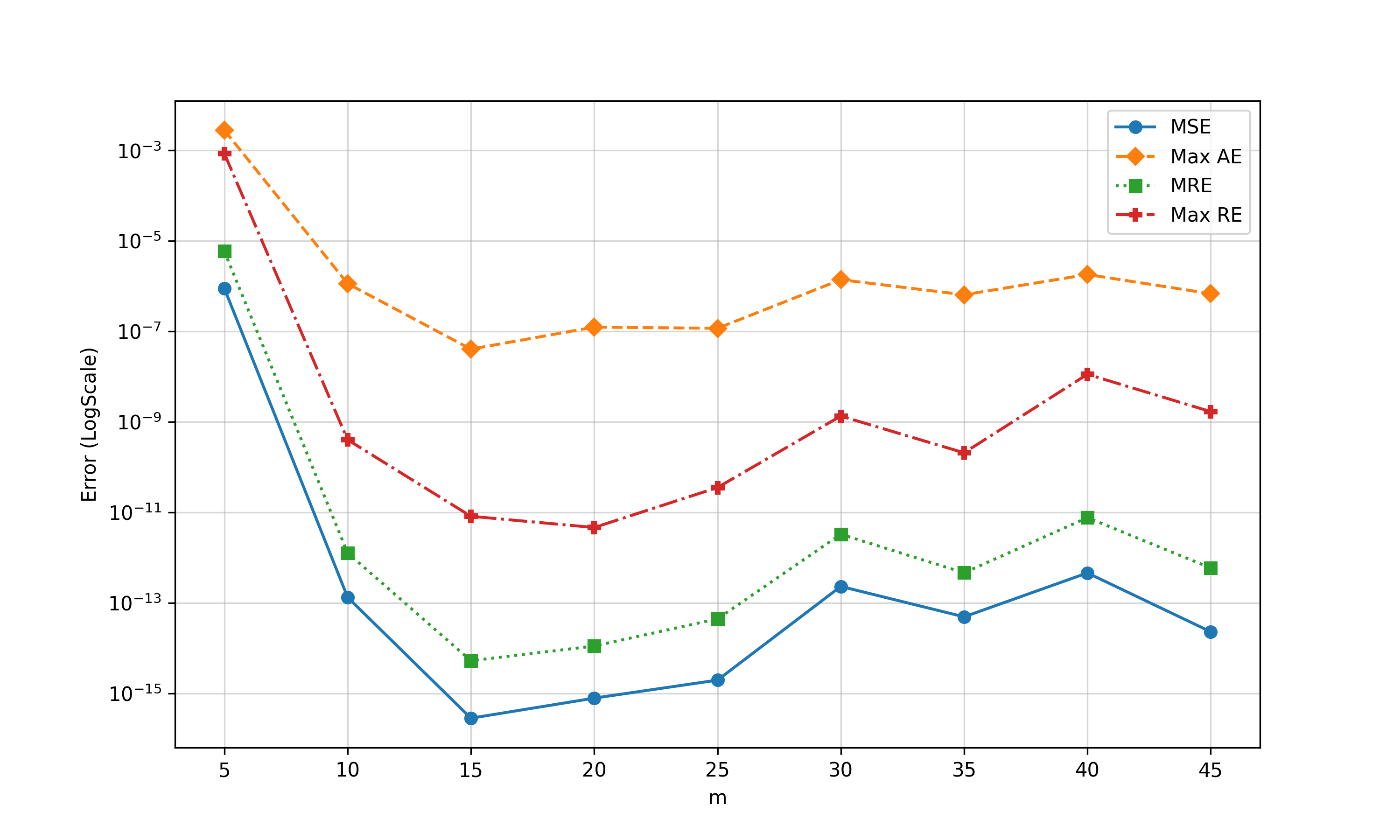} & \includegraphics[width=0.3\textwidth]{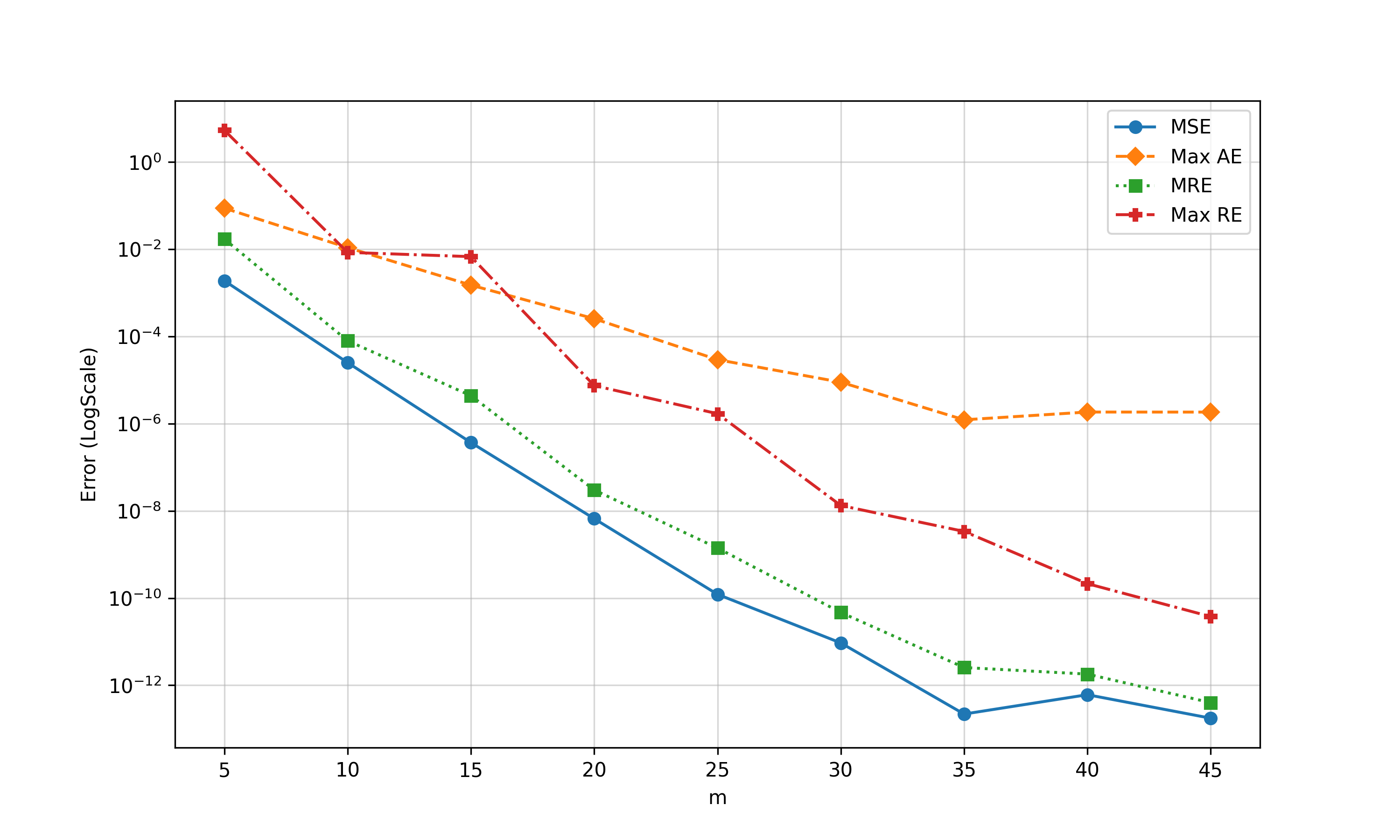} \\
        
    \end{tabular}
    \caption{Errors versus $m$ on the ellipse. First row: $f_1,f_2, f_3$, second row: $f_4,f_5, f_6$.}
    \label{fig:errors-ellipse}
\end{figure}

\begin{figure}
    \centering\begin{tabular}{ccc}
        \includegraphics[width=0.3\textwidth]{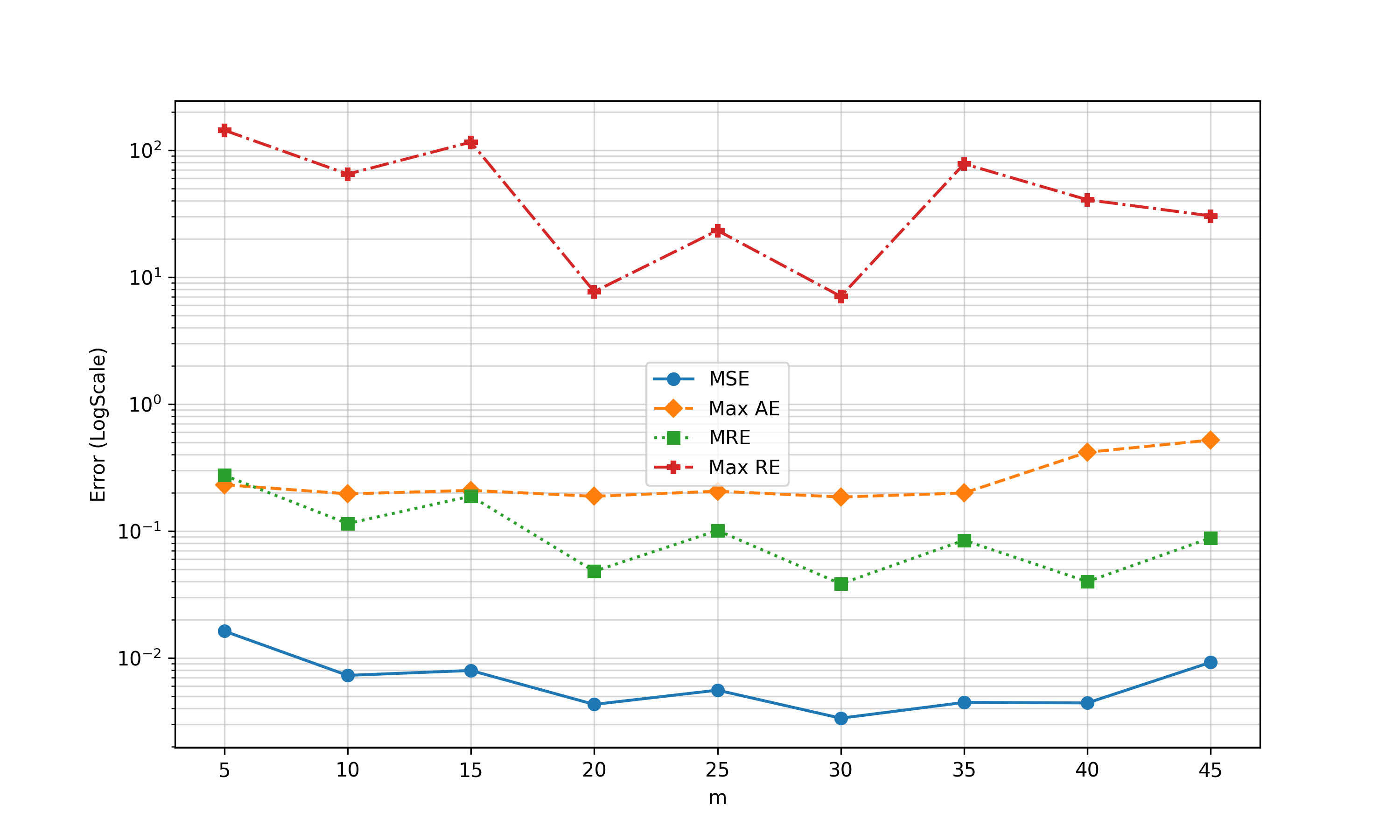} & \includegraphics[width=0.3\textwidth]{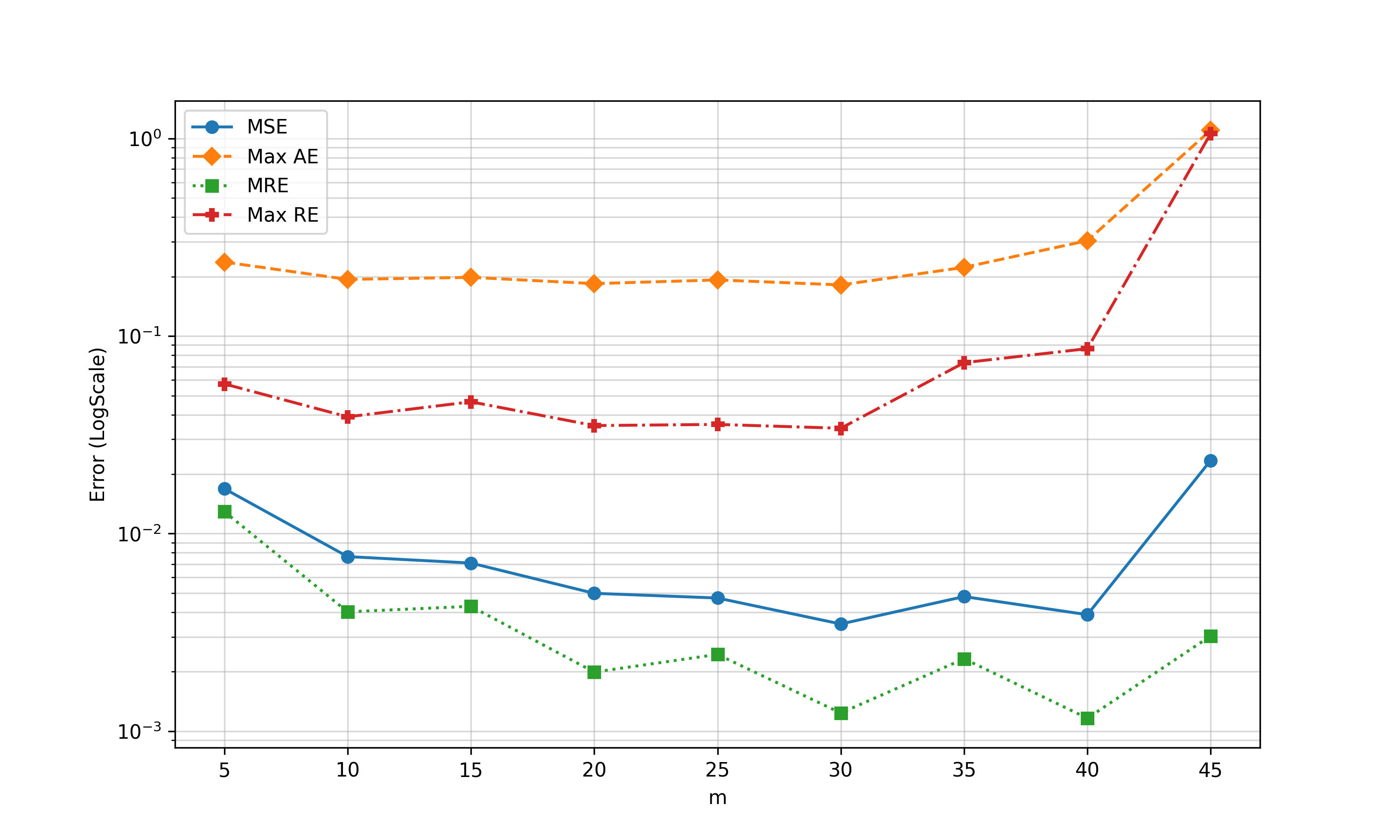} &
        \includegraphics[width=0.3\textwidth]{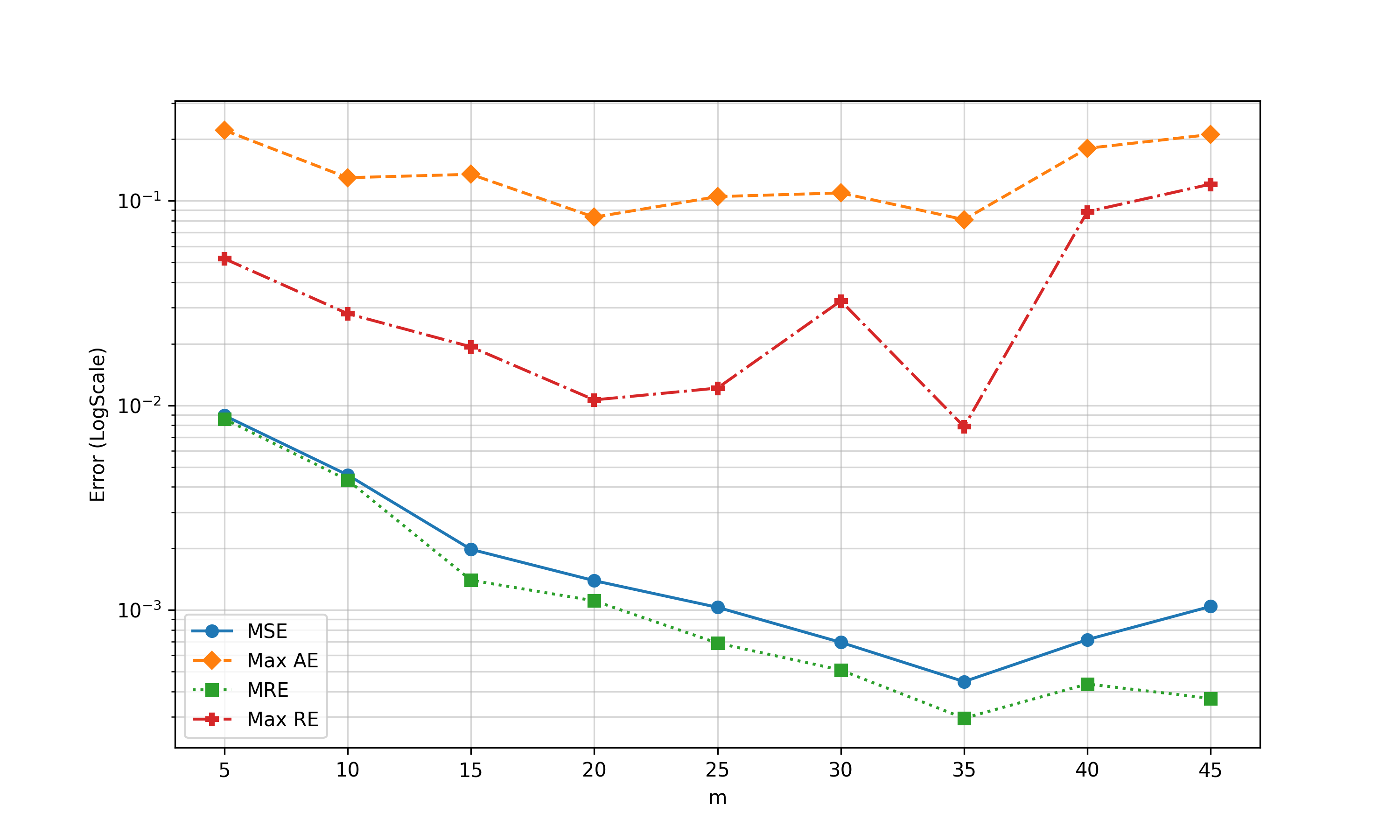} \\ \includegraphics[width=0.3\textwidth]{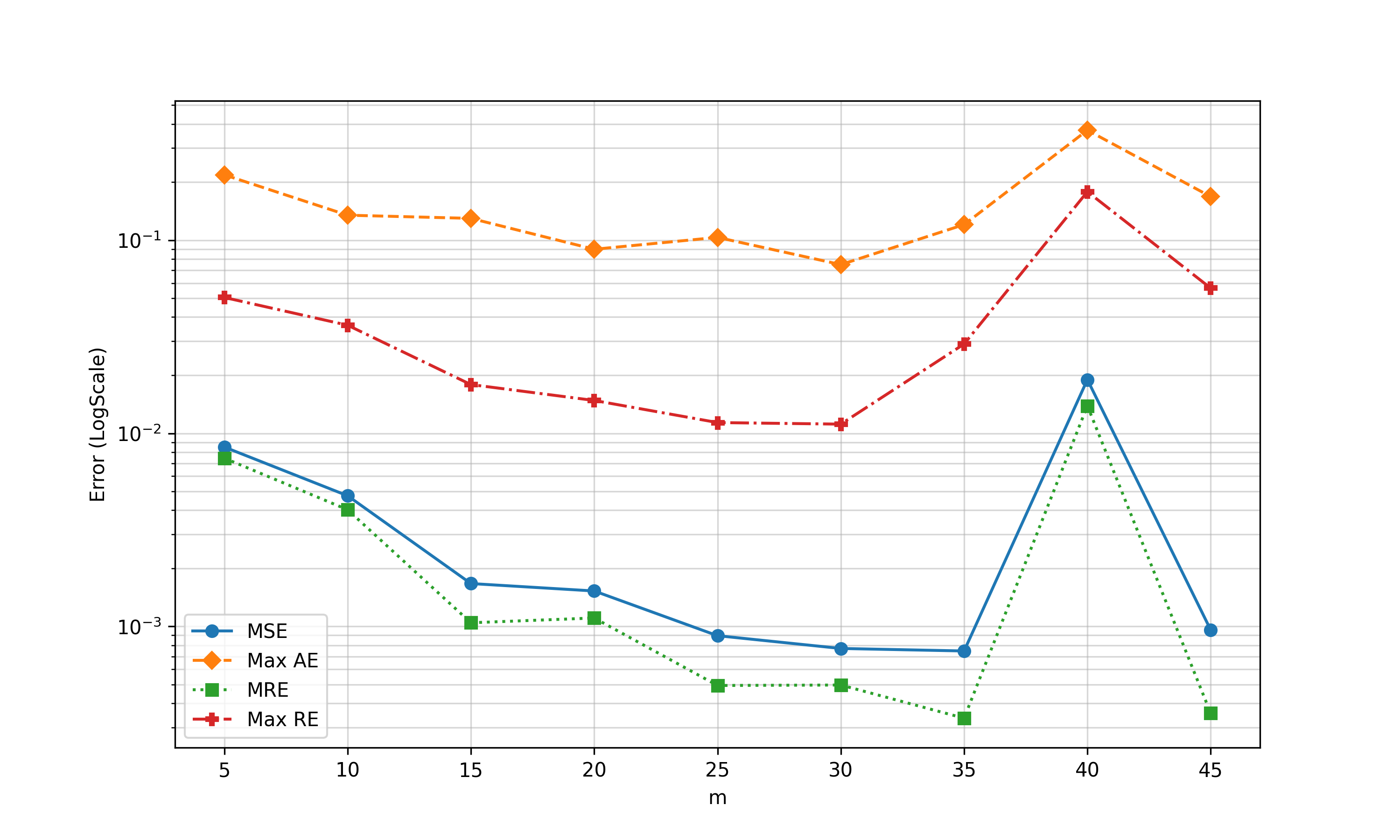} &
        \includegraphics[width=0.3\textwidth]{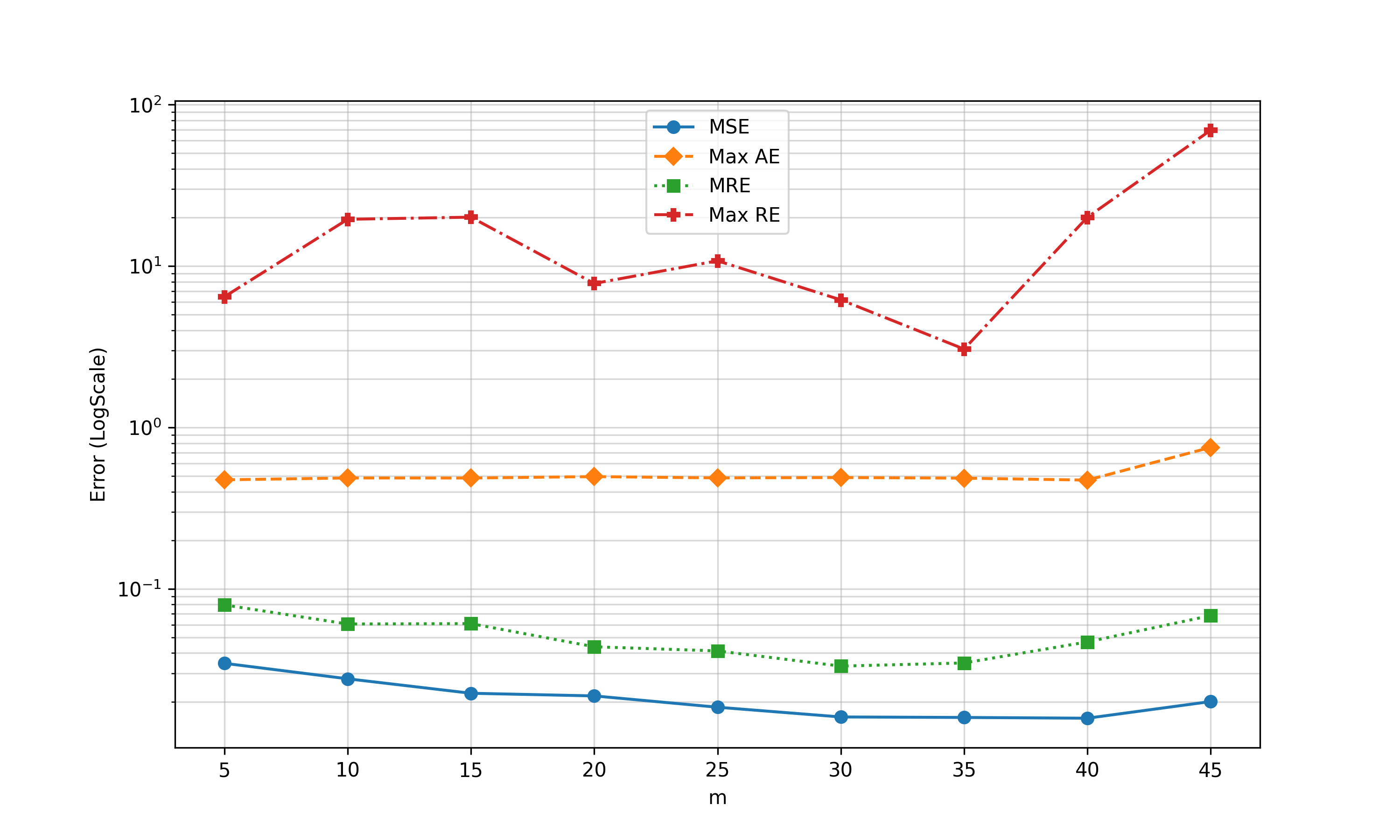} & \includegraphics[width=0.3\textwidth]{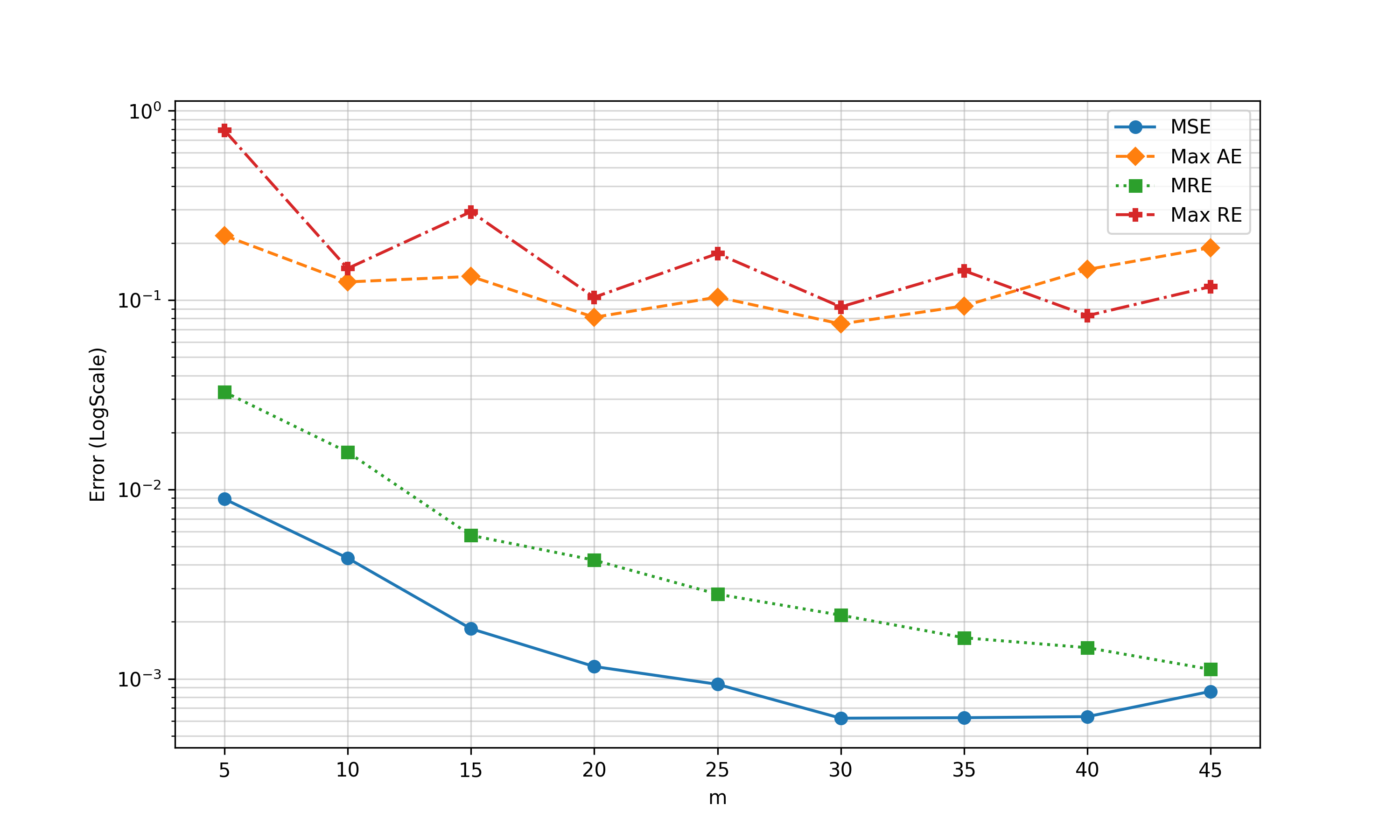} \\
        
    \end{tabular}
    \caption{Errors versus $m$ on the annulus. First row: $f_1,f_2, f_3$, second row: $f_4,f_5, f_6$.}
    \label{fig:errors-annulus}
\end{figure}

\begin{figure}
    \centering\begin{tabular}{ccc}
        \includegraphics[width=0.3\textwidth]{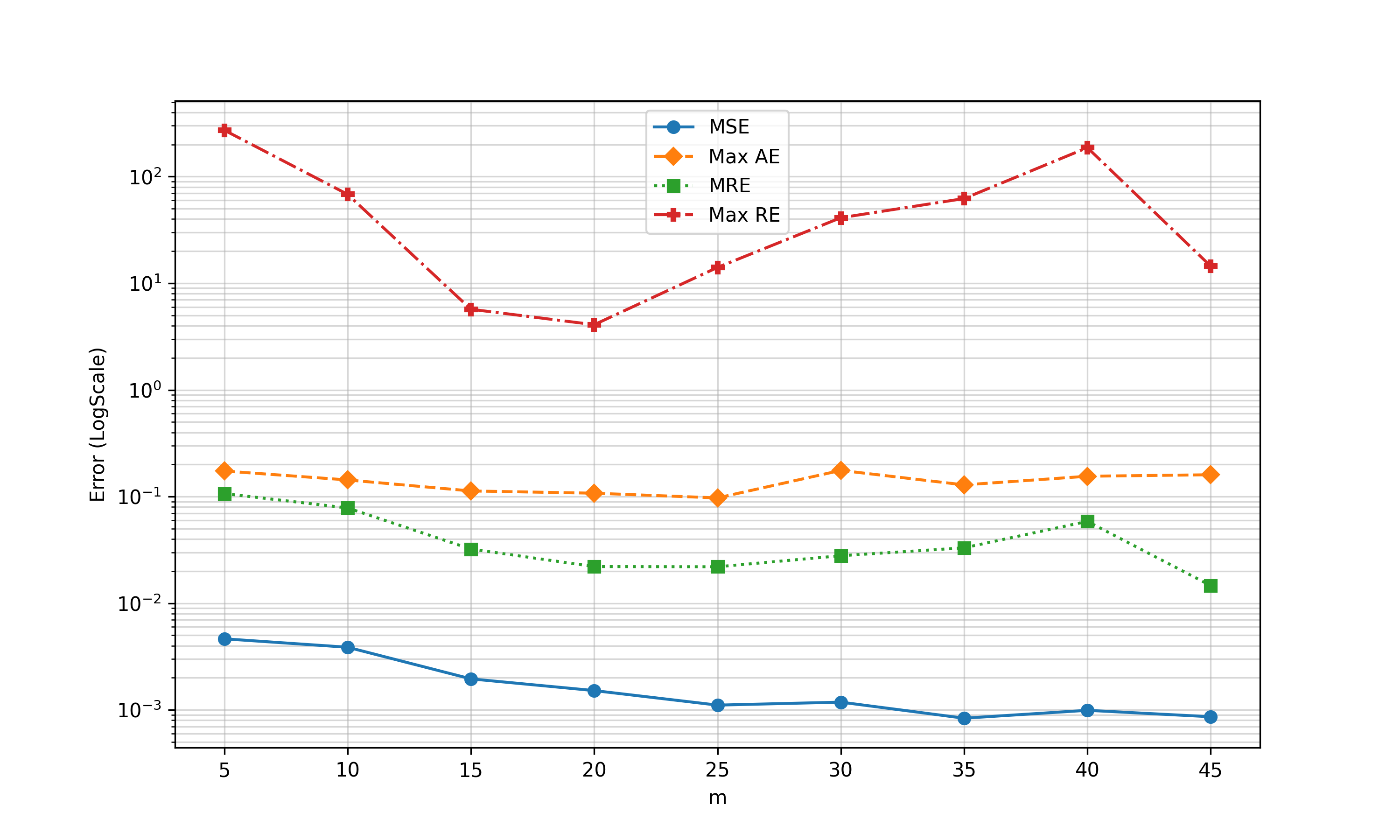} & \includegraphics[width=0.3\textwidth]{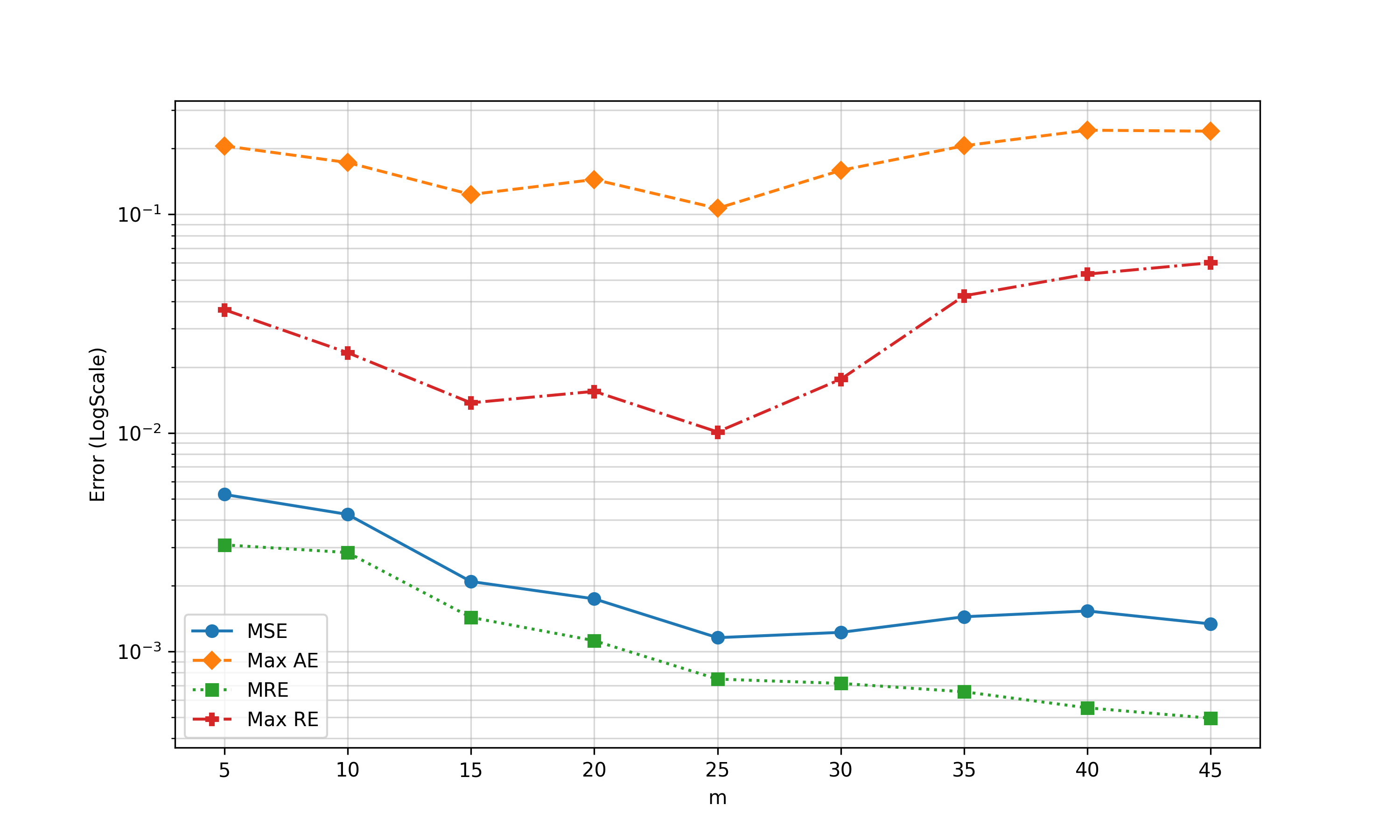} &
        \includegraphics[width=0.3\textwidth]{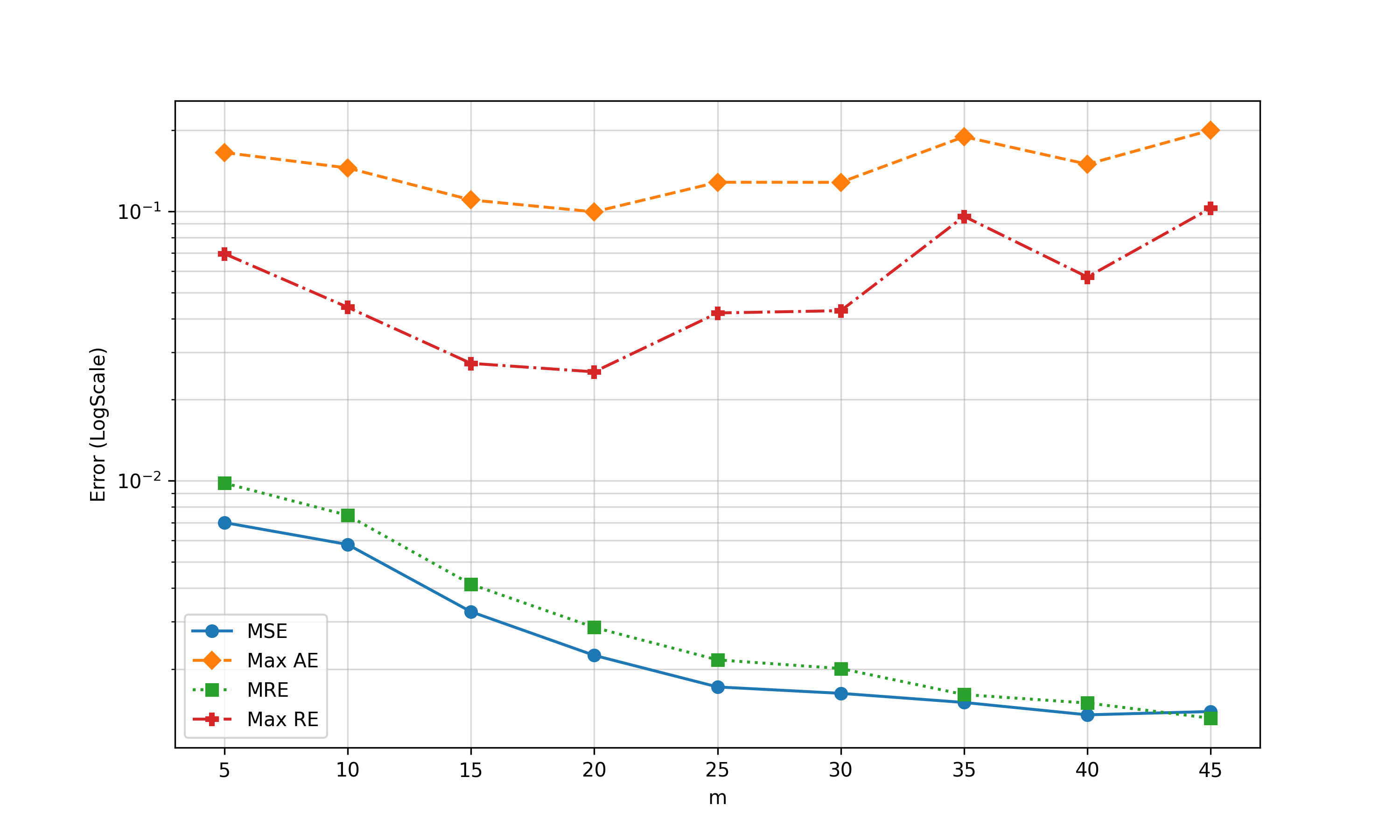} \\ \includegraphics[width=0.3\textwidth]{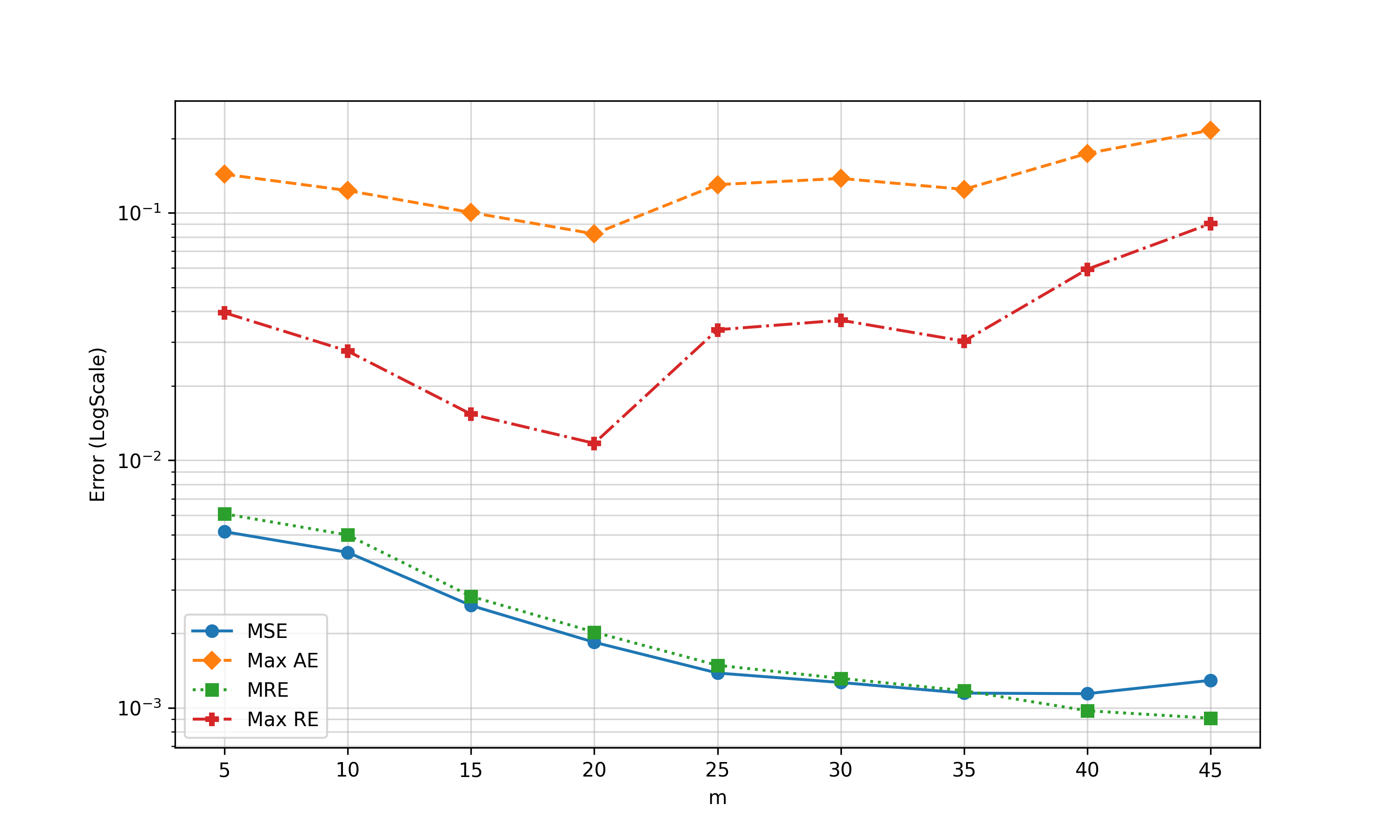} &
        \includegraphics[width=0.3\textwidth]{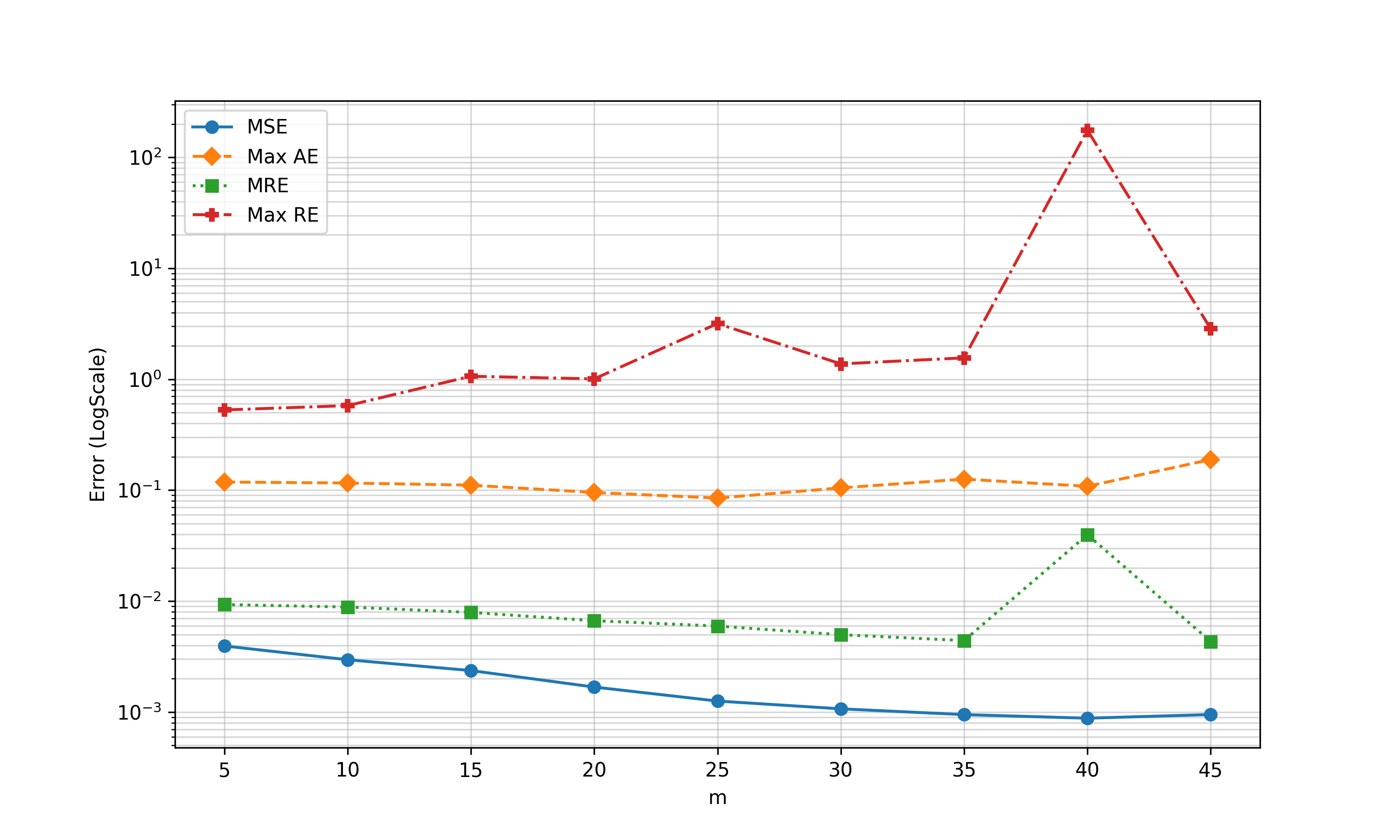} & \includegraphics[width=0.3\textwidth]{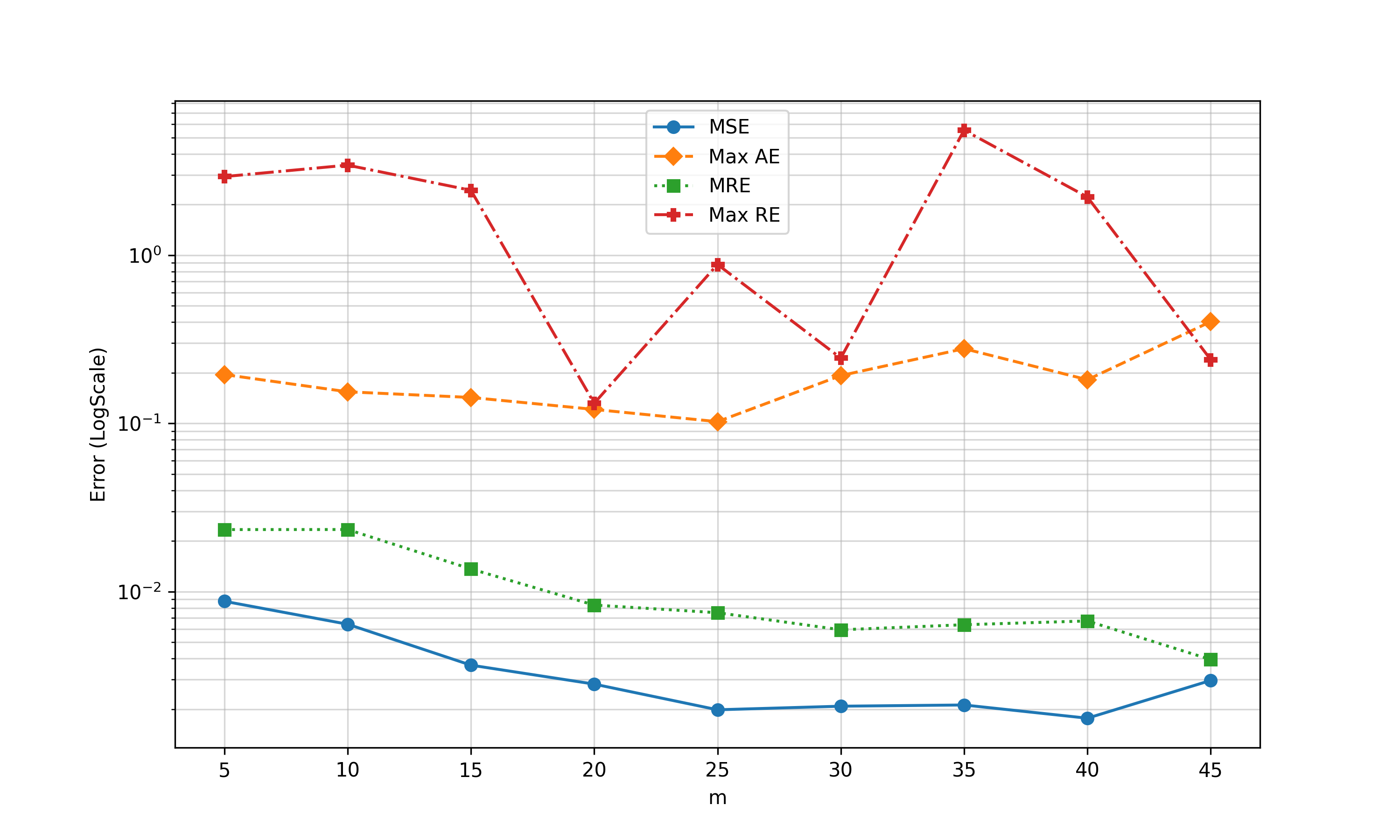} \\
        
    \end{tabular}
    \caption{Errors versus $m$ on the polygon. First row: $f_1,f_2, f_3$, second row: $f_4,f_5, f_6$.}
    \label{fig:errors-polygon}
\end{figure}

\pagebreak

The numerical results for $f_2$ are presented in Table~\ref{tab:results-f2-ellipse} for the ellipse, Table~\ref{tab:results-f2-annulus} for the annulus, and Table~\ref{tab:results-f2-polygon} for the polygon. Each entry in these tables reports, for a given interpolation dimension $m$, the corresponding regression dimension $\tilde r = m + \lfloor\sqrt m\rfloor$, the dimensions $M$ and $\tilde R$ of the polynomial spaces $\mathbb P_{m}(\mathbb R^2)$ and $\mathbb P_{\tilde r}(\mathbb R^2)$, the previously defined error metrics, and the execution time (\text{ExTime}).

The errors obtained for the ellipse are substantially lower than those observed for the annulus and the polygon. This behavior is expected, as the mapping from the unit disk to the ellipse is linear. For the annulus and the polygon, the errors predominantly range between $10^{-2}$ and $10^{-4}$, with noticeably higher values of MaxAE and MaxRE as $m$ increases.

Furthermore, the execution time increases rapidly with $m$, reaching approximately $1400$ seconds ($\approx 23$ minutes) when $m=45$.

Figures~\ref{fig:errors-ellipse}, \ref{fig:errors-annulus}, and \ref{fig:errors-polygon} illustrate the behavior of the different error metrics as a function of $m$. The solid blue line with circle markers represents the MSE, the dashed orange line with diamond markers corresponds to the MaxAE, the dotted green line with square markers denotes the MRE, and the dash-dotted red line with plus markers indicates the MaxRE.

Ultimately, it is essential to strike a balance between accuracy and computational efficiency, occasionally accepting slightly lower precision to obtain more efficient approximations. Based on the plots and tabulated data, setting $m = 20 = n/5$ appears to be an optimal choice: it yields high accuracy while keeping the execution time around $90$ seconds in most experiments, which represents a very reasonable computational cost.\\

Next, we fix $m$ and vary the dimension of the regression part $\tilde r$, to study the effect of the difference between the two parameters.

\FloatBarrier

\subsection{Varying the dimension of the regression part}

In the previous section, the proposed operator was computed and validated for each test function by fixing $n=100$ and varying the dimension of the interpolation part $m$, while keeping $\tilde r = m +\lfloor\sqrt m\rfloor$. In this section, we fix $m=\lfloor n/5\rfloor=20$ as a choice that balances accuracy and efficiency, and instead vary the dimension of the regression part, $\tilde r$.

Since $n=100$, the interpolation dimension is fixed at $m=20$ with $M=231$ throughout this section. Given that $\tilde r > m$, we begin by setting $\tilde r = 25$ and progressively increase this value. Tables~\ref{tab:results-f2-ellipse-r}, \ref{tab:results-f2-annulus-r}, and \ref{tab:results-f2-polygon-r} display the results obtained for $f_2$ on the ellipse, the annulus, and the polygon, respectively. As observed previously, the results on the ellipse are substantially better due to the simplicity of its underlying mapping. Notably, increasing the value of $\tilde r$ does not initially produce significant differences; however, when $\tilde r$ becomes too large, both accuracy and efficiency deteriorate, resulting in larger errors and increased execution times (\text{ExTime}).

\begin{figure}
    \centering\begin{tabular}{ccc}
        \includegraphics[width=0.3\textwidth]{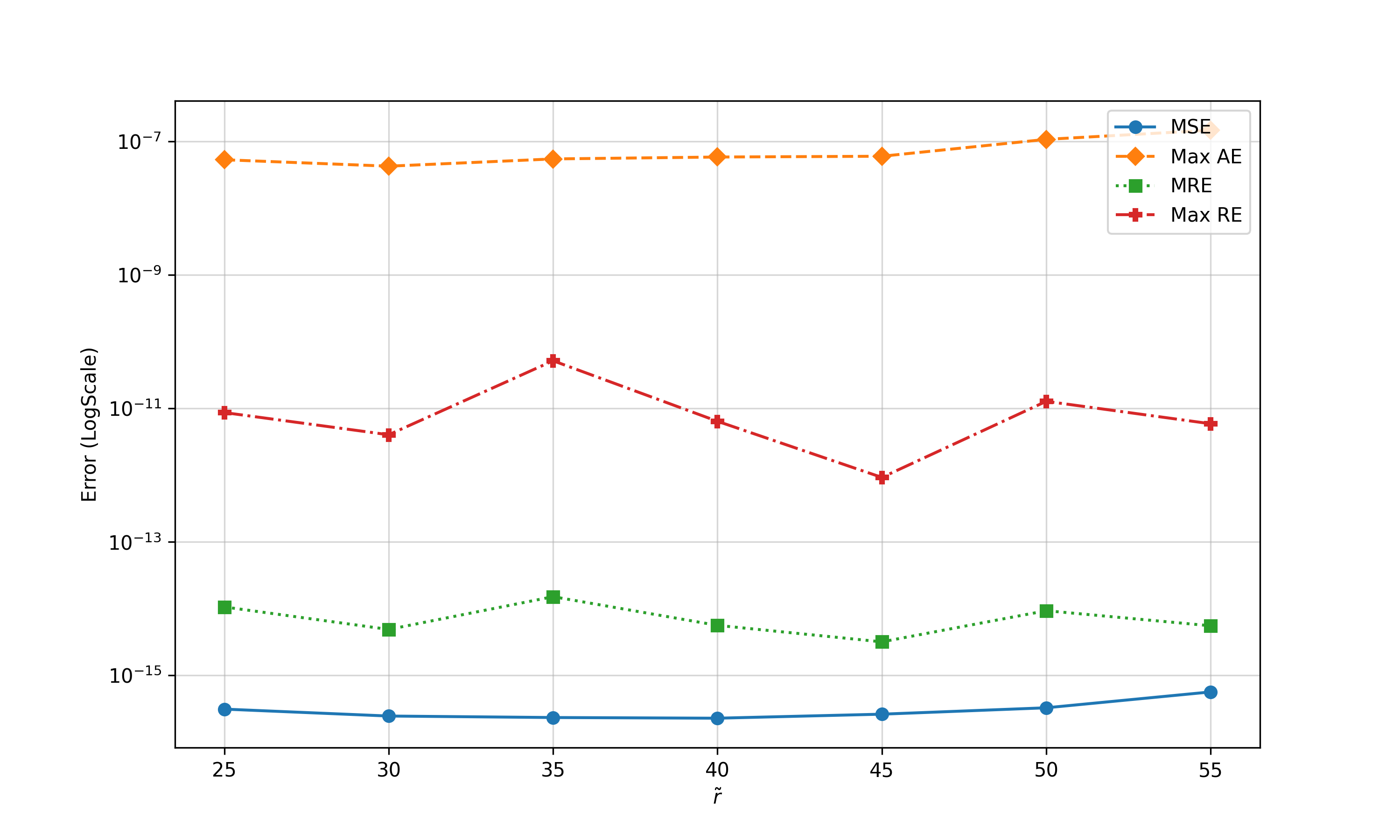} & \includegraphics[width=0.3\textwidth]{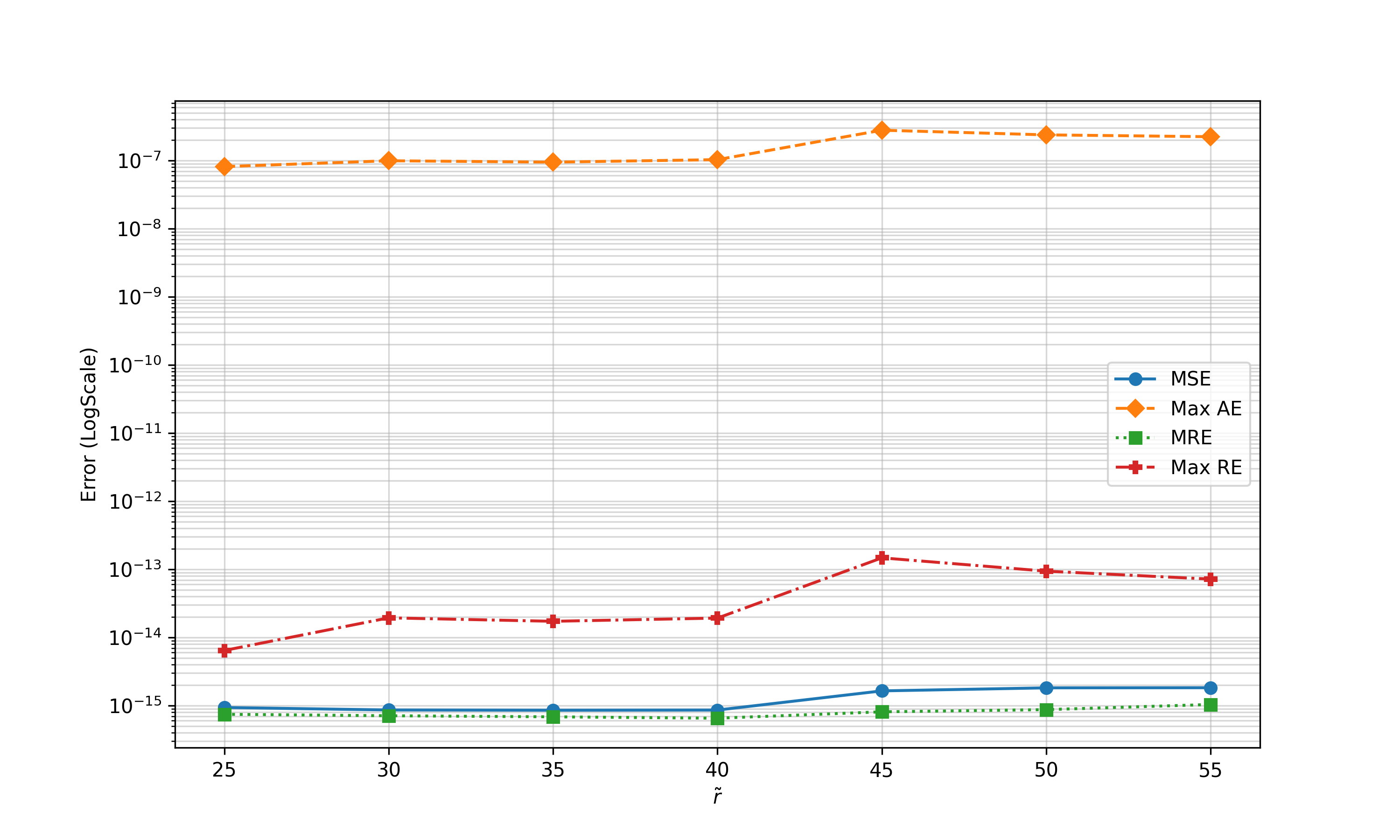} &
        \includegraphics[width=0.3\textwidth]{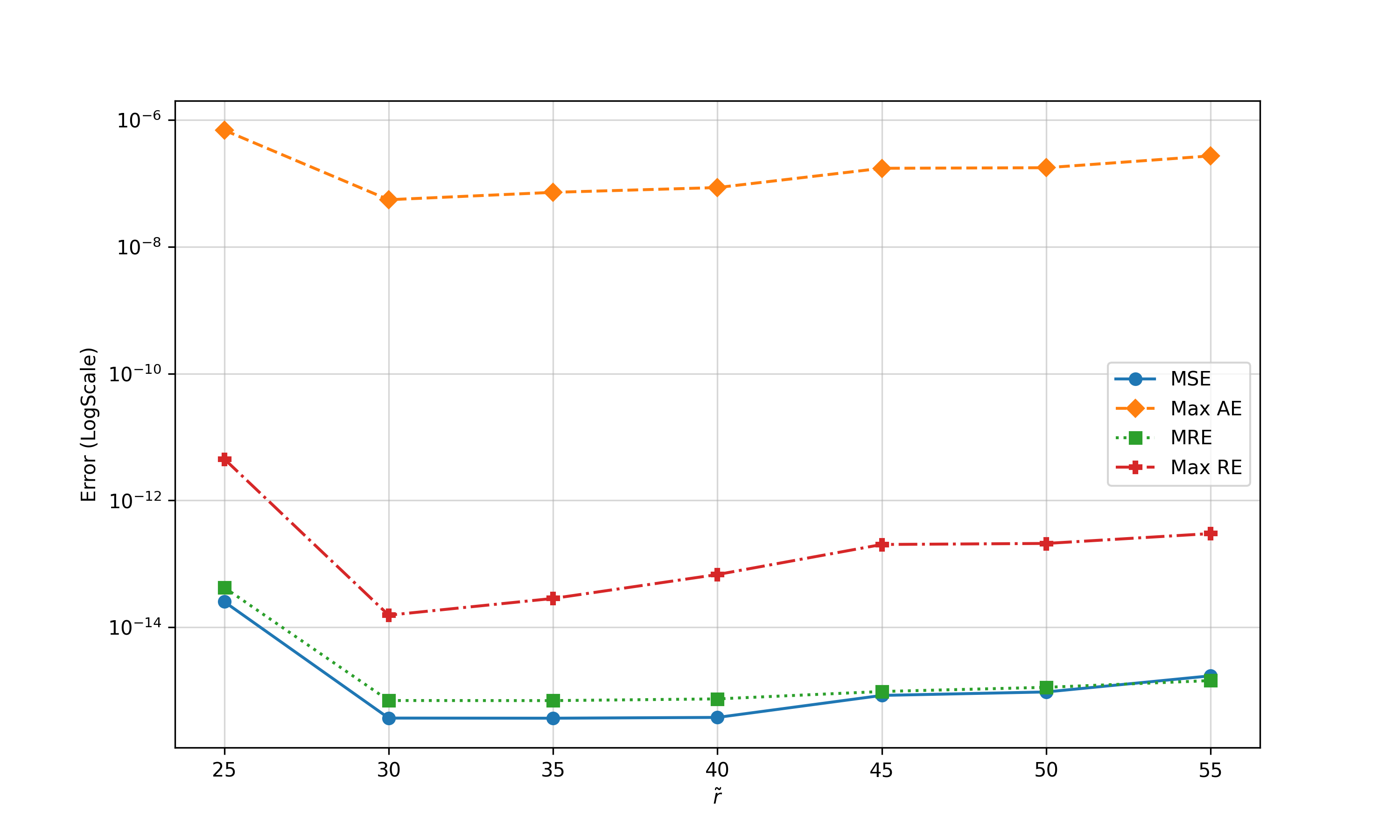} \\ \includegraphics[width=0.3\textwidth]{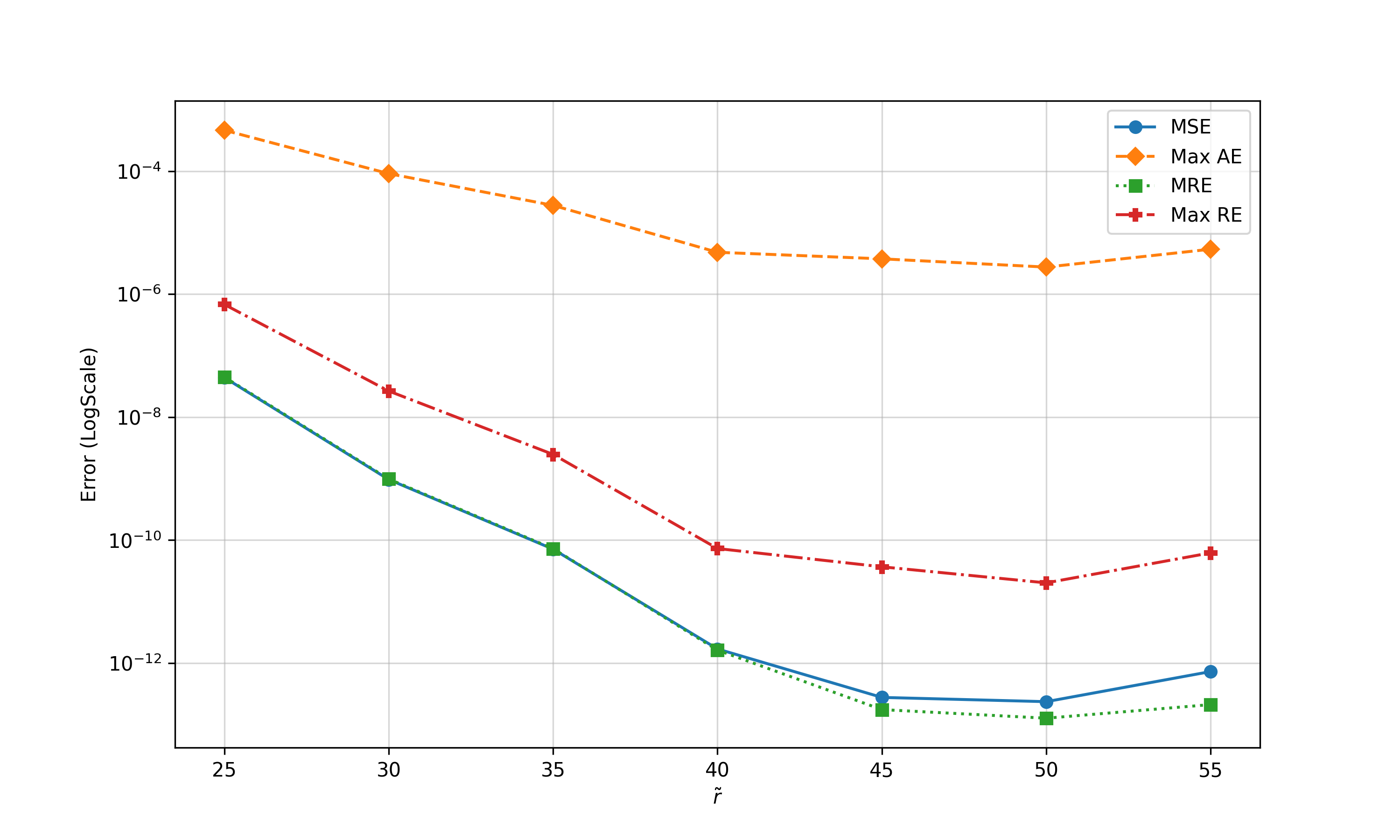} &
        \includegraphics[width=0.3\textwidth]{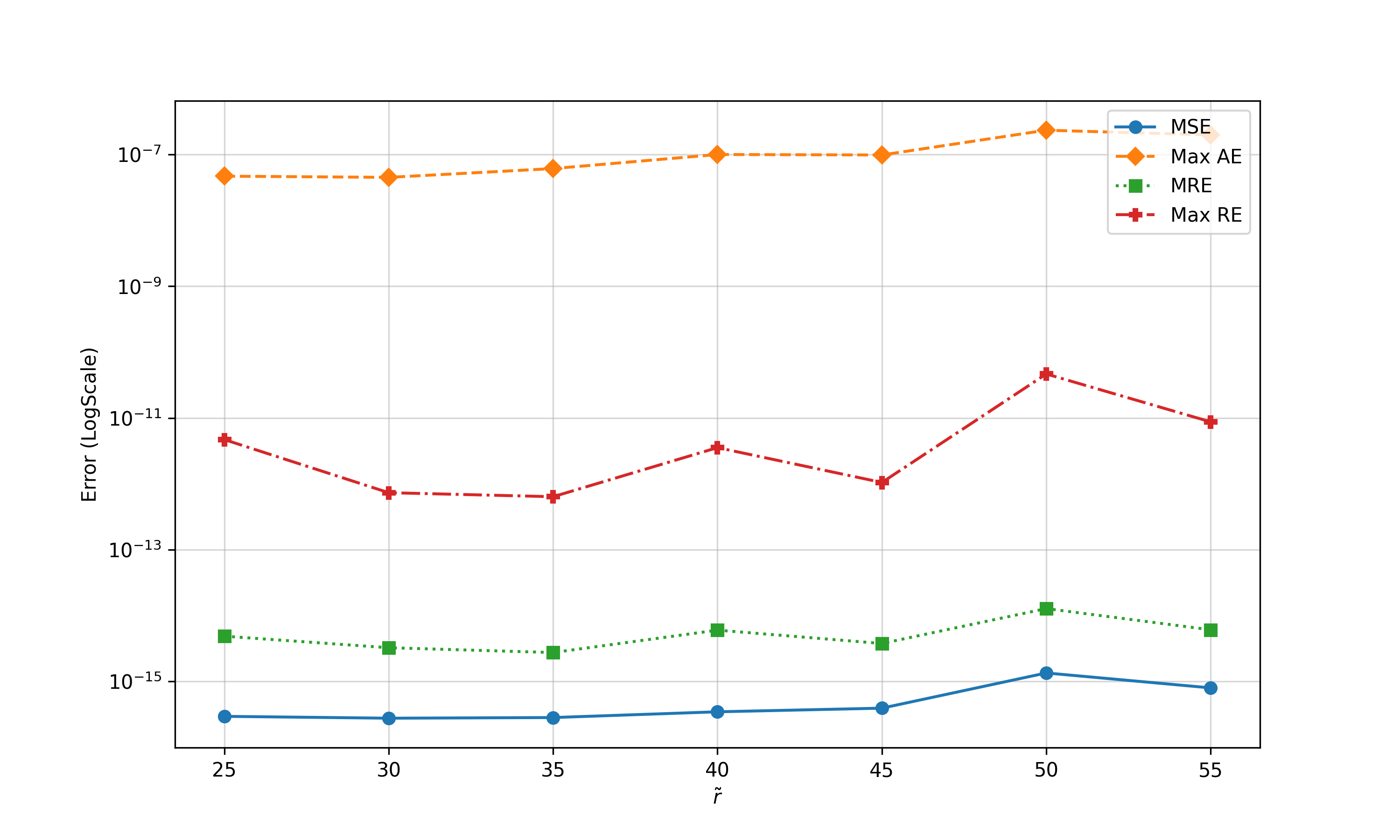} & \includegraphics[width=0.3\textwidth]{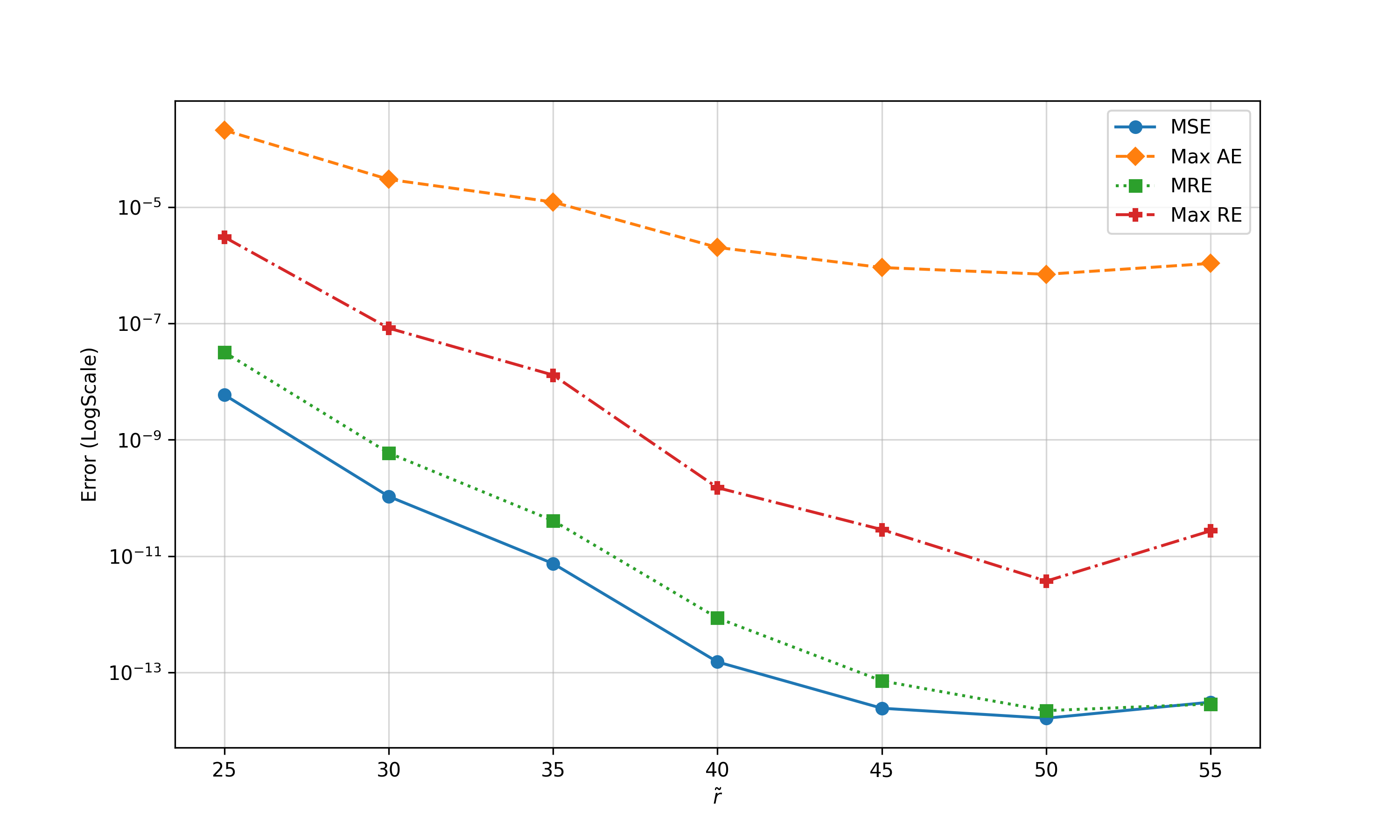} \\ 
    \end{tabular}
    \caption{Errors versus $\tilde r$ on the ellipse. First row: $f_1,f_2, f_3$, second row: $f_4,f_5, f_6$.}
    \label{fig:errors-ellipse-r}
\end{figure}

\begin{figure}
    \centering\begin{tabular}{ccc}
        \includegraphics[width=0.3\textwidth]{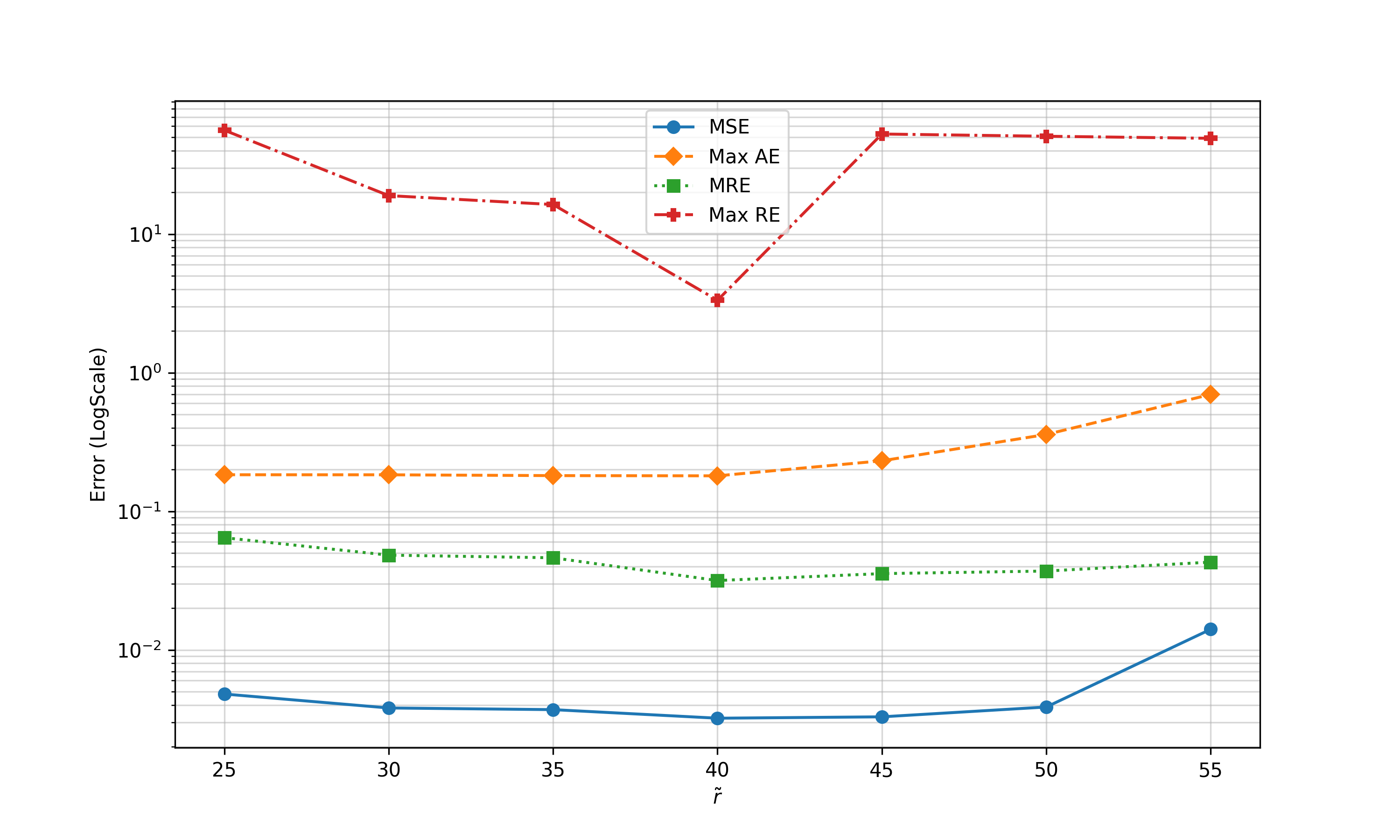} & \includegraphics[width=0.3\textwidth]{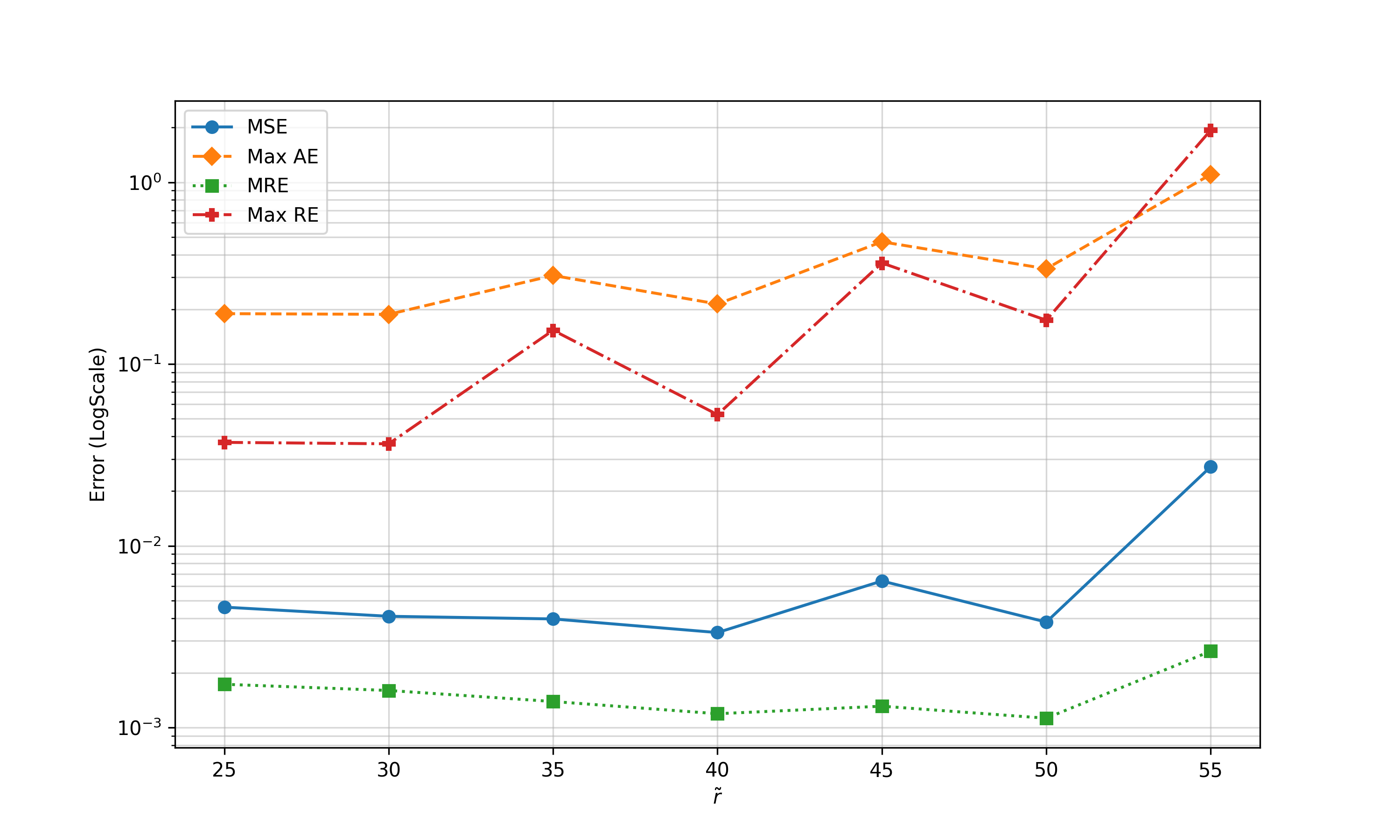} &
        \includegraphics[width=0.3\textwidth]{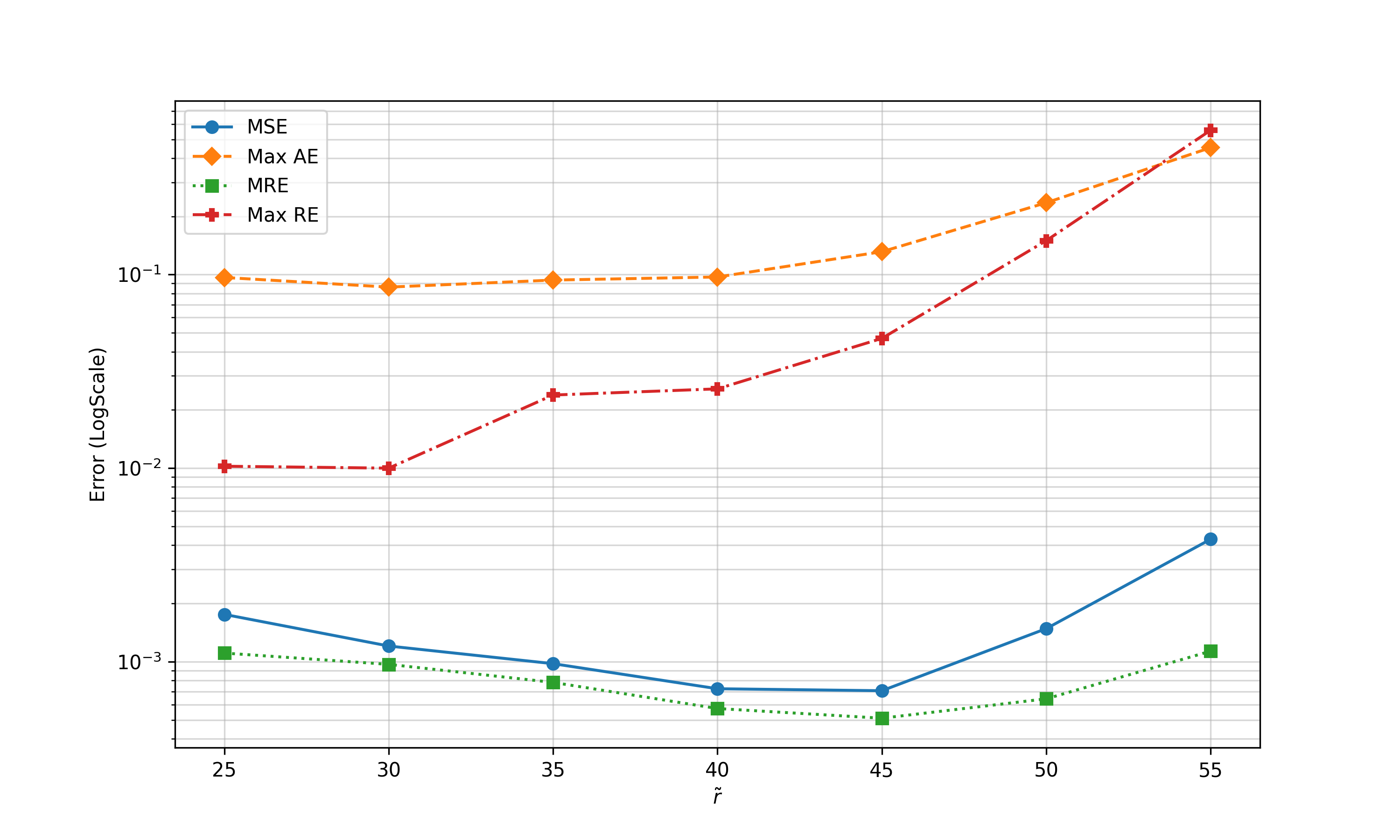} \\ \includegraphics[width=0.3\textwidth]{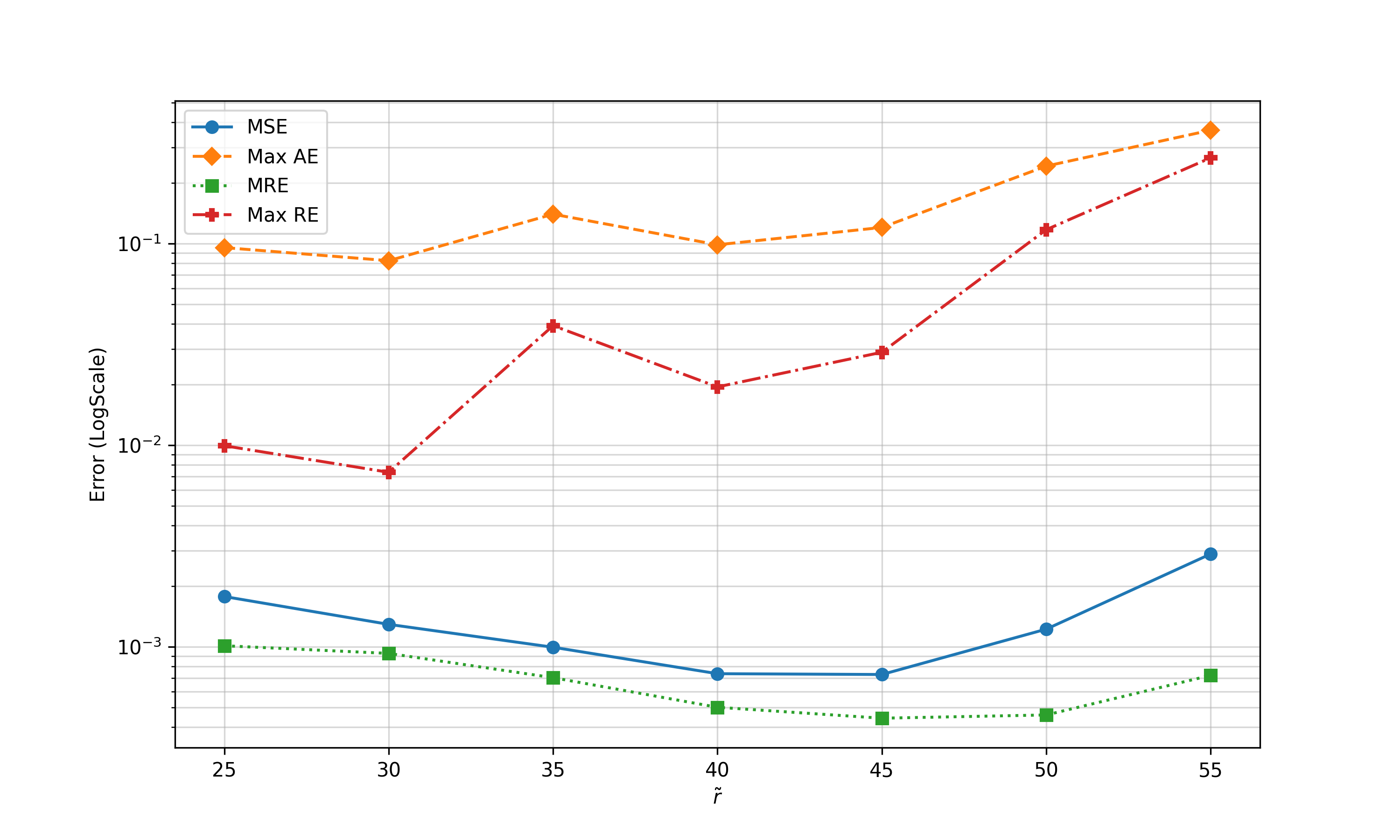} &
        \includegraphics[width=0.3\textwidth]{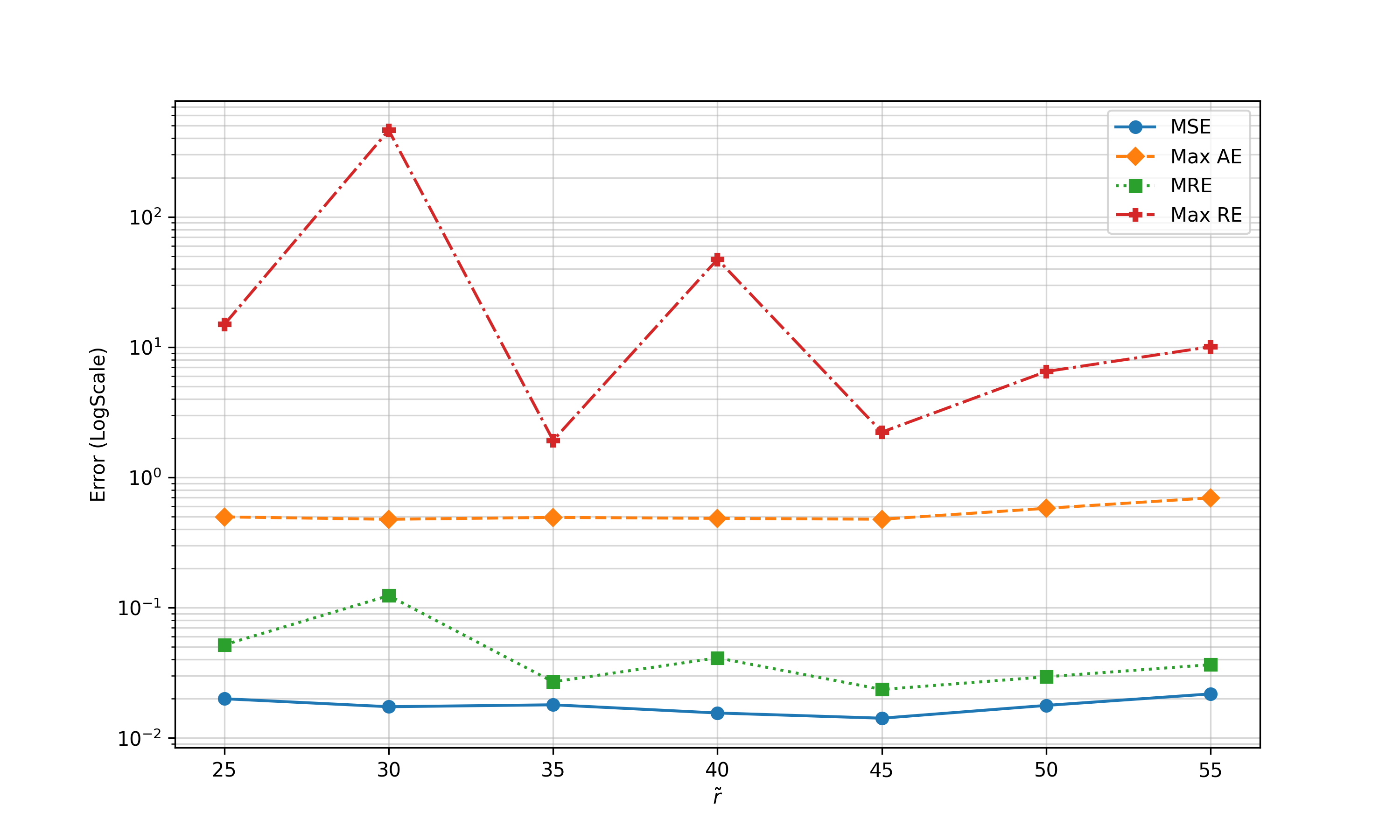} & \includegraphics[width=0.3\textwidth]{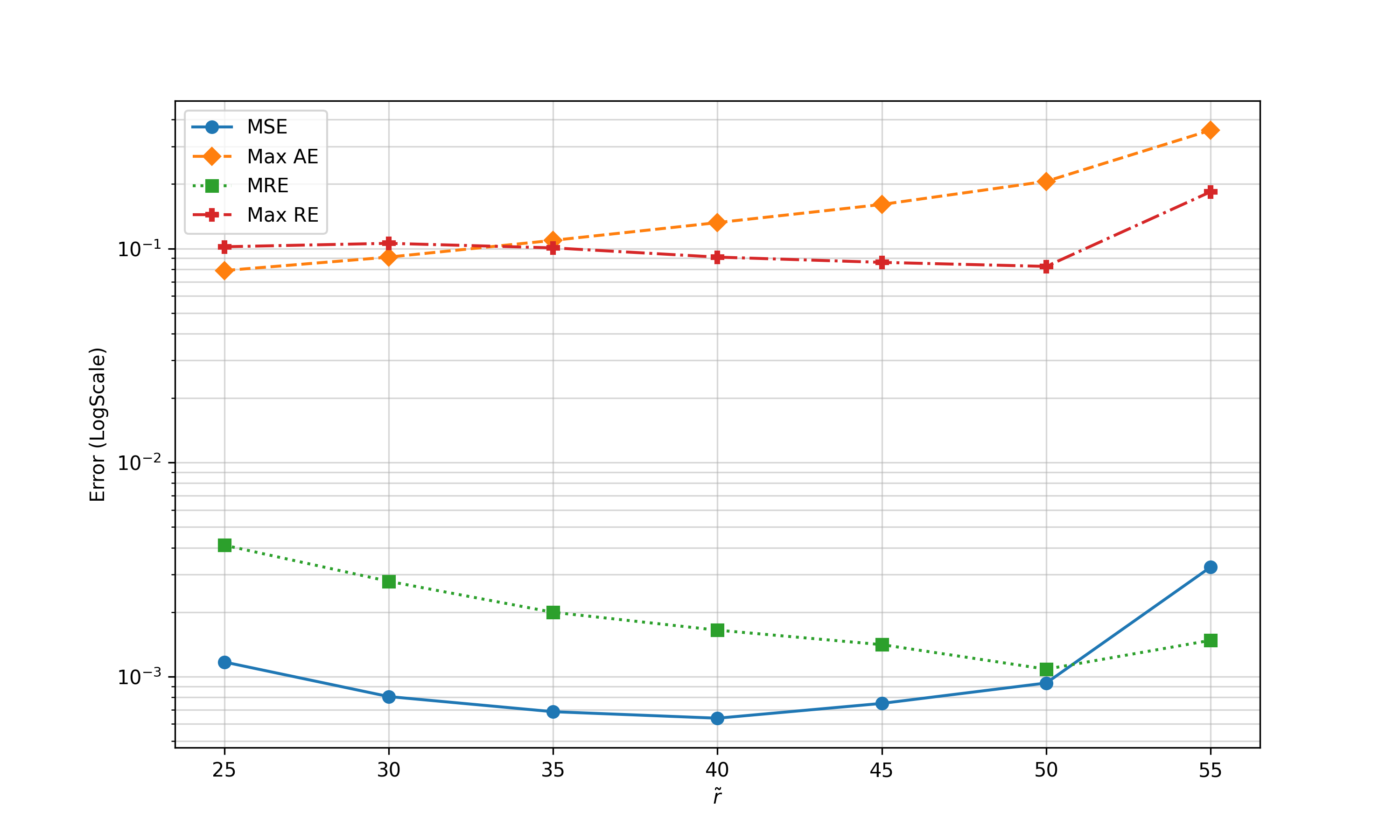} \\ 
    \end{tabular}
    \caption{Errors versus $\tilde r$ on the annulus. First row: $f_1,f_2, f_3$, second row: $f_4,f_5, f_6$.}
    \label{fig:errors-annulus-r}
\end{figure}

\begin{figure}
    \centering\begin{tabular}{ccc}
        \includegraphics[width=0.3\textwidth]{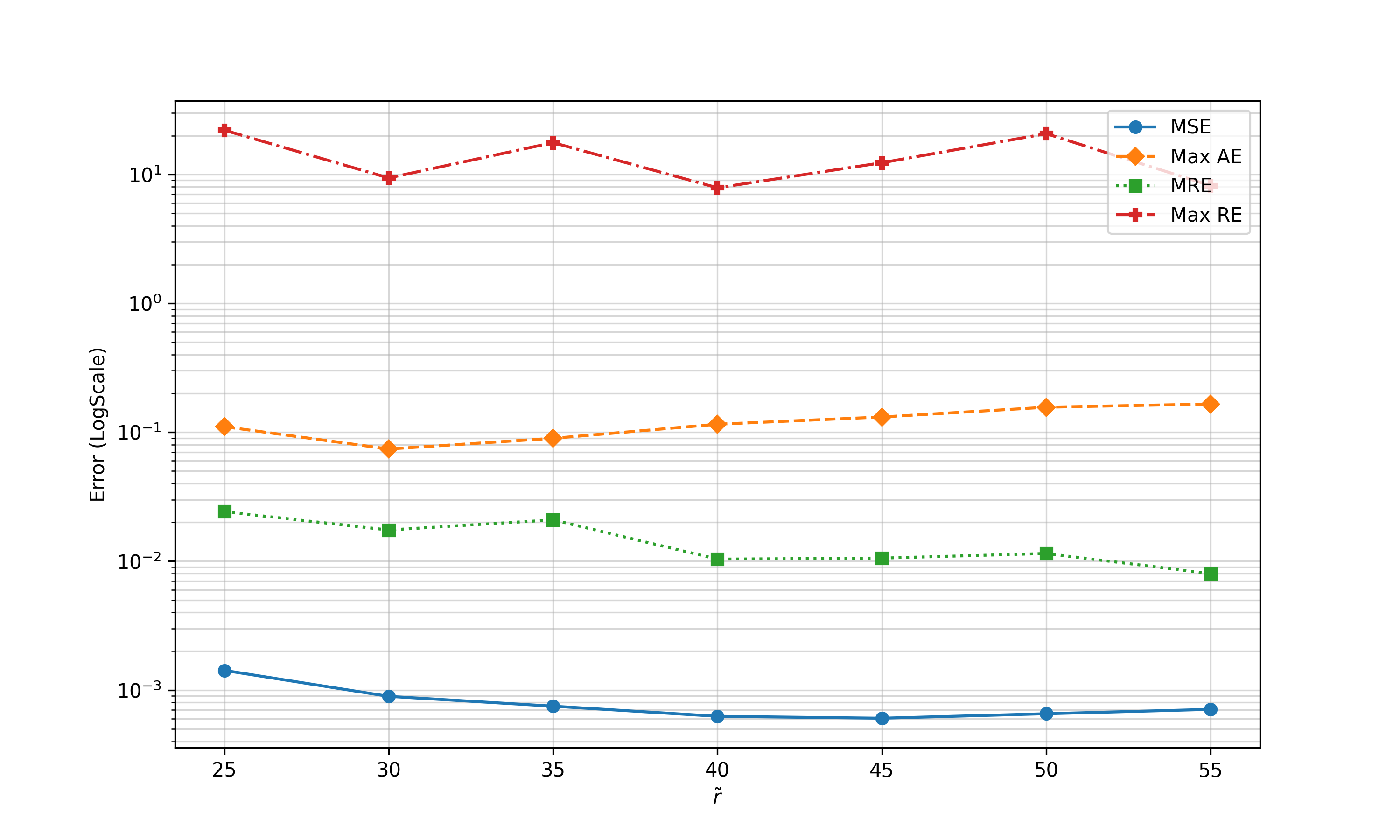} & \includegraphics[width=0.3\textwidth]{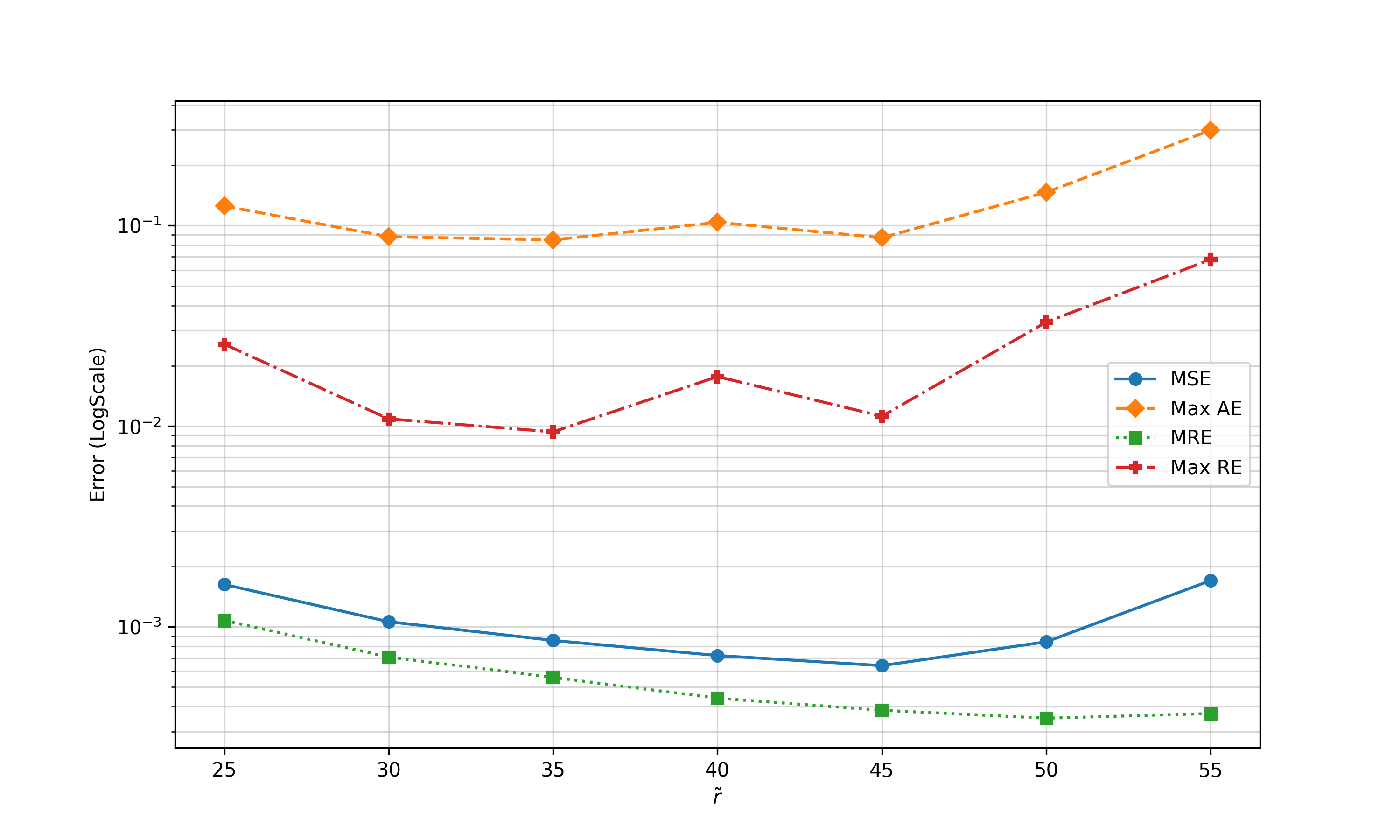} &
        \includegraphics[width=0.3\textwidth]{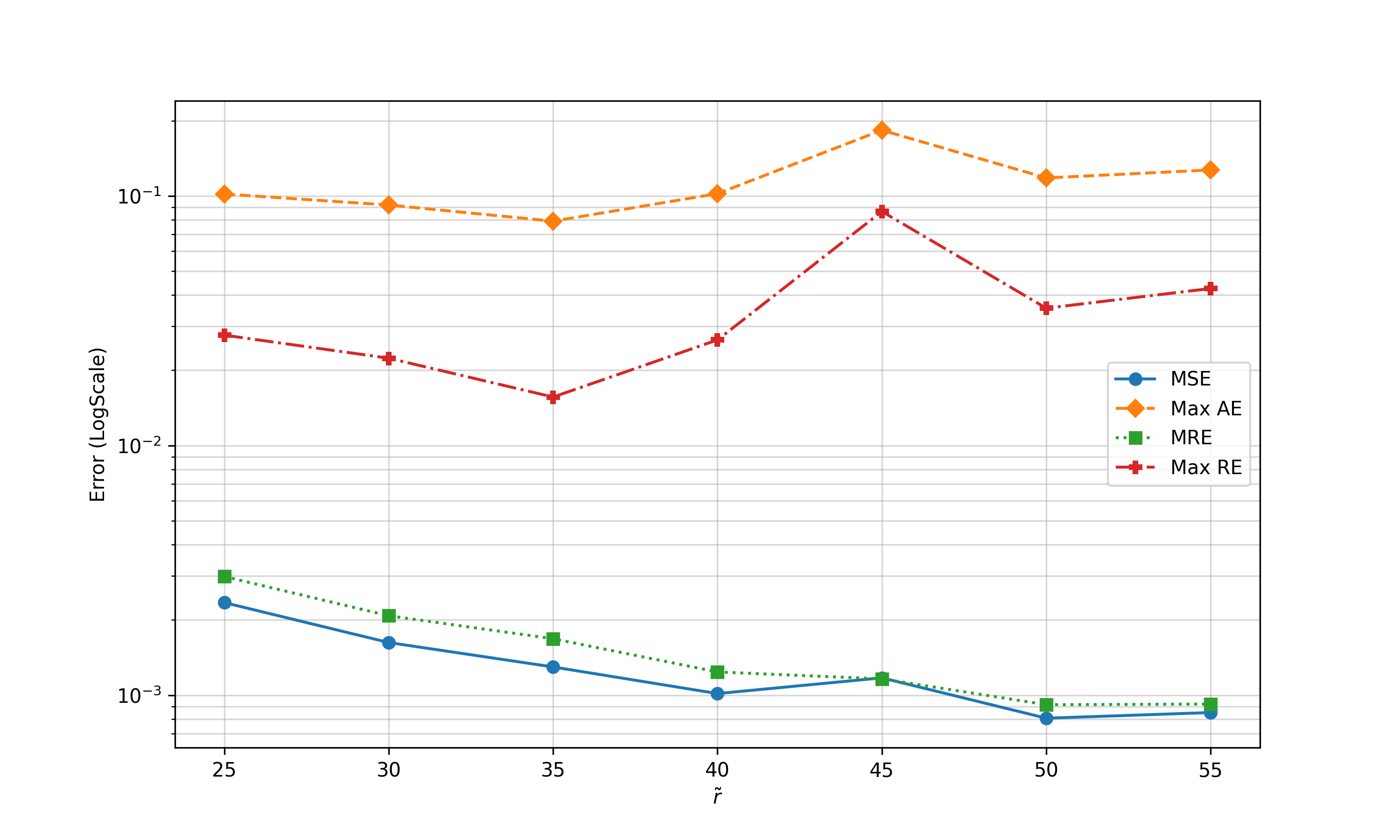} \\ \includegraphics[width=0.3\textwidth]{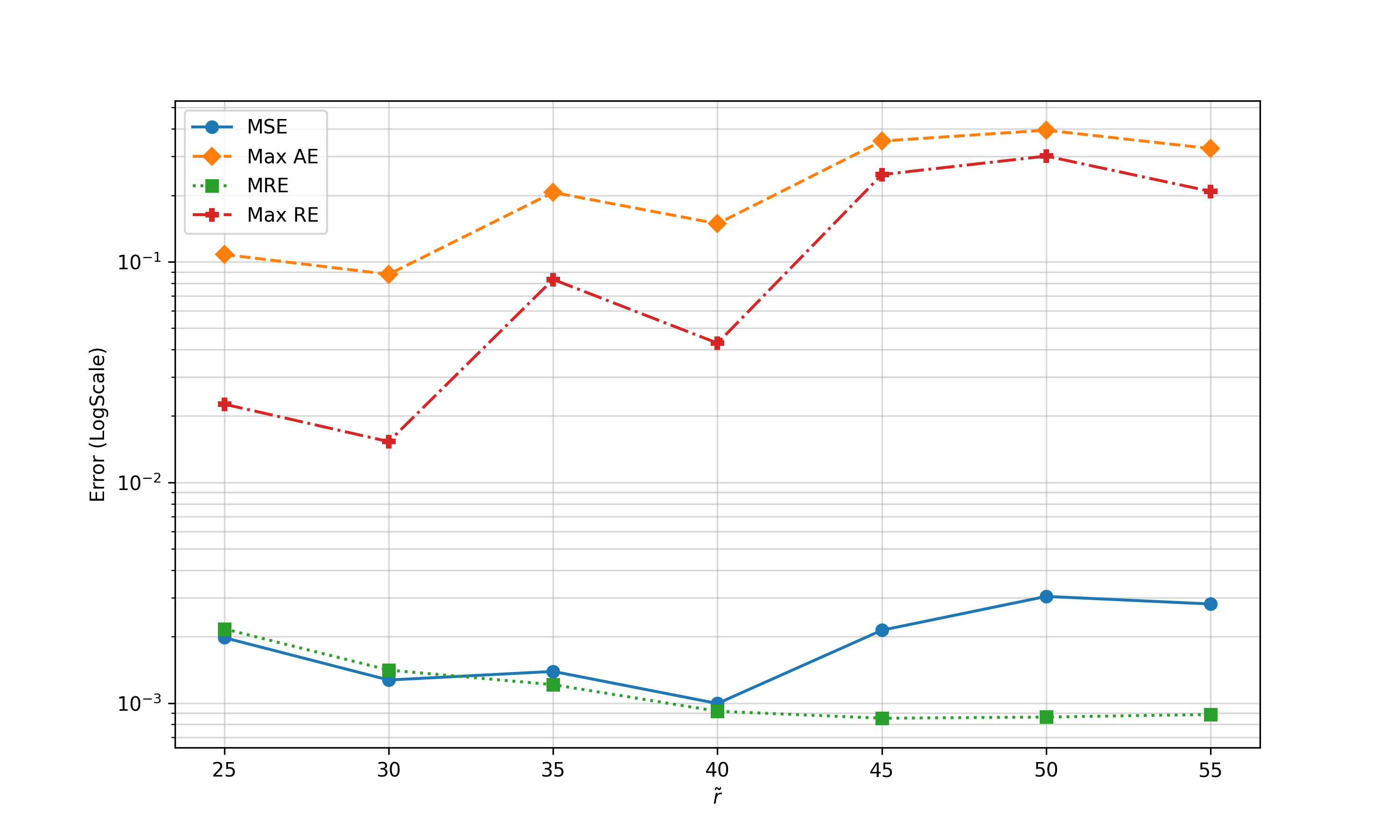} &
        \includegraphics[width=0.3\textwidth]{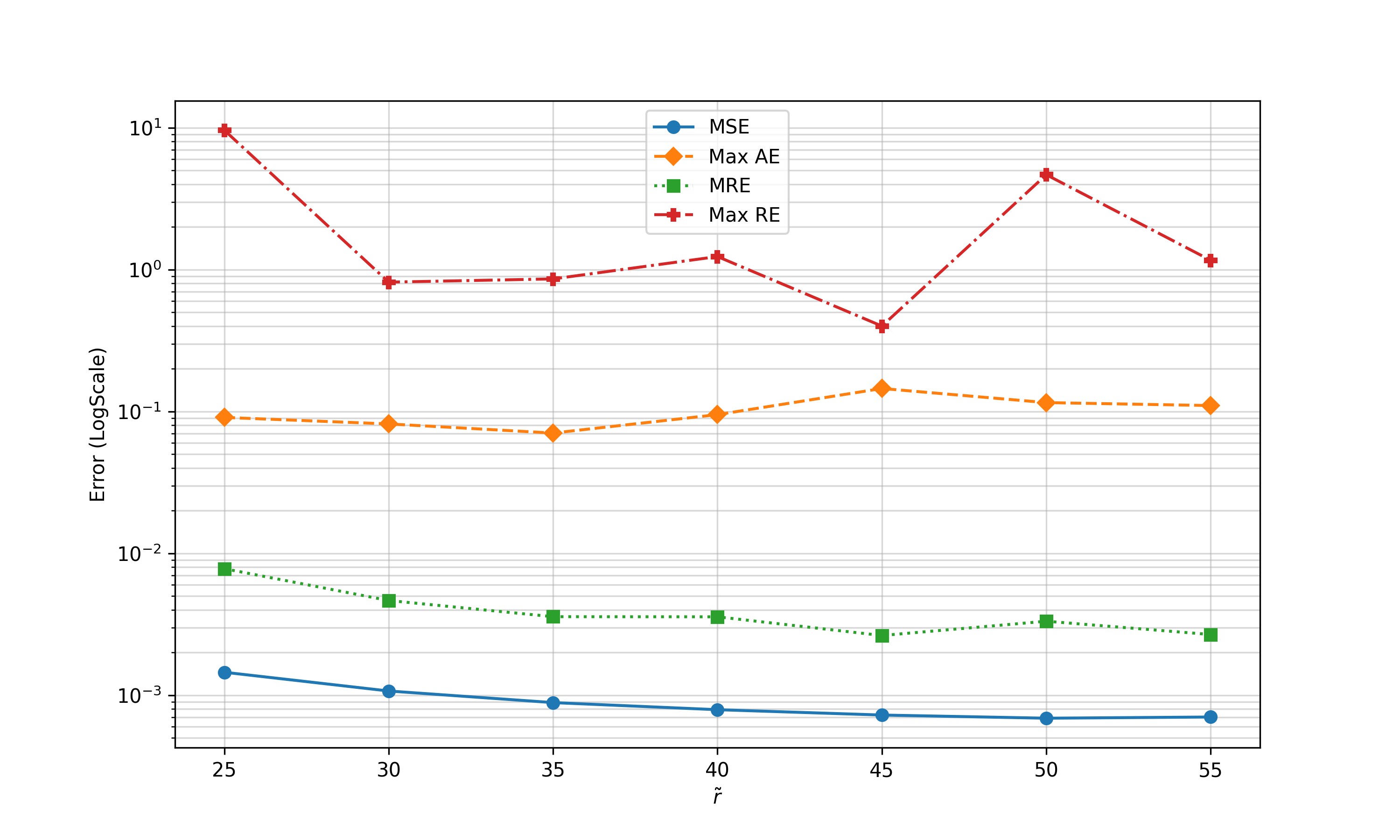} & \includegraphics[width=0.3\textwidth]{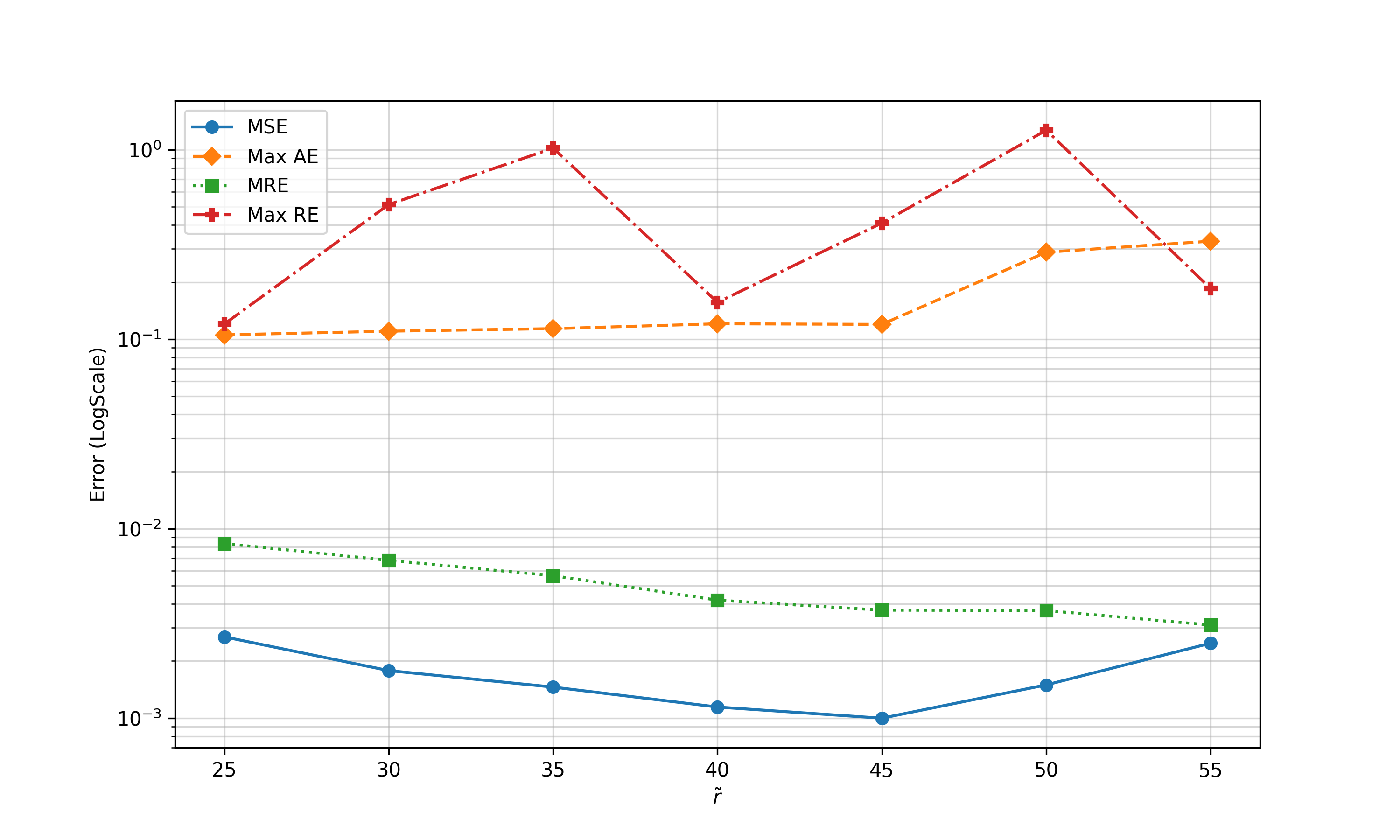} \\ 
    \end{tabular}
    \caption{Errors versus $\tilde r$ on the polygon. First row: $f_1,f_2, f_3$, second row: $f_4,f_5, f_6$.}
    \label{fig:errors-polygon-r}
\end{figure}

The remaining results concerning the approximation of the functions $f_1, \dots, f_6$ are illustrated in Figures~\ref{fig:errors-ellipse-r}, \ref{fig:errors-annulus-r}, and \ref{fig:errors-polygon-r}. While Figure~\ref{fig:errors-ellipse-r} reveals noticeable performance variations for $f_3$, $f_4$, and $f_6$ on the ellipse, the approximations on the annulus and the polygon are far less sensitive to variations in $\tilde r$. As shown in Figures~\ref{fig:errors-annulus-r} and \ref{fig:errors-polygon-r}, the error curves decrease by barely an order of magnitude, and most of them eventually increase for larger values of $\tilde r$ (as expected). These findings suggest that computational efficiency should be prioritized, as larger values of $\tilde r$ do not necessarily improve accuracy but invariably increase the computational cost of the operator.

%\pagebreak

\begin{table}[!ht]
    \centering
    
    \begin{tabular}{ccccccc}
    \toprule
        $\tilde r$ & $\tilde R$ & \textbf{MSE} & \textbf{MaxAE} & \textbf{MRE} & \textbf{MaxRE} & \textbf{ExTime} \\ \toprule
        25 & 351 & 9.38722e-16 & 8.1617e-08 & 7.46852e-16 & 6.44899e-15 & 103.661 \\ \hline
        30 & 496 & 8.63509e-16 & 9.94031e-08 & 7.10322e-16 & 1.94233e-14 & 196.919 \\ \hline
        35 & 666 & 8.56258e-16 & 9.48304e-08 & 6.85007e-16 & 1.72929e-14 & 343.577 \\ \hline
        40 & 861 & 8.61288e-16 & 1.03238e-07 & 6.53768e-16 & 1.92757e-14 & 560.487 \\ \hline
        45 & 1081 & 1.64574e-15 & 2.79173e-07 & 8.11188e-16 & 1.47529e-13 & 869.137 \\ \hline
        50 & 1326 & 1.82353e-15 & 2.38419e-07 & 8.73527e-16 & 9.41106e-14 & 1290.9 \\ \hline
        55 & 1596 & 1.83079e-15 & 2.24014e-07 & 1.03918e-15 & 7.19659e-14 & 1850.6 \\ \bottomrule
    \end{tabular}
    \caption{Results varying $\tilde r$ for $f_2$ on the ellipse}
    \label{tab:results-f2-ellipse-r}
\end{table}

\begin{table}[!ht]
    \centering
    \begin{tabular}{ccccccc}
    \toprule
        $\tilde r$ & $\tilde R$ & \textbf{MSE} & \textbf{MaxAE} & \textbf{MRE} & \textbf{MaxRE} & \textbf{ExTime} \\ \toprule
        25 & 351 & 0.00459835 & 0.189267 & 0.00173096 & 0.0370683 & 102.431 \\ \hline
        30 & 496 & 0.00409551 & 0.187603 & 0.00160015 & 0.0364001 & 195.6 \\ \hline
        35 & 666 & 0.00396018 & 0.306979 & 0.00139148 & 0.153327 & 342.204 \\ \hline
        40 & 861 & 0.0033322 & 0.213863 & 0.00119205 & 0.0526855 & 559.824 \\ \hline
        45 & 1081 & 0.00640472 & 0.470402 & 0.00131135 & 0.359608 & 868.598 \\ \hline
        50 & 1326 & 0.00380542 & 0.334224 & 0.00112627 & 0.174184 & 1291.73 \\ \hline
        55 & 1596 & 0.0272258 & 1.10313 & 0.00263095 & 1.93311 & 1855.01 \\ \bottomrule
    \end{tabular}
    \caption{Results varying $\tilde r$ for $f_2$ on the annulus}
    \label{tab:results-f2-annulus-r}
\end{table}

\begin{table}[!ht]
    \centering
    \begin{tabular}{ccccccc}
    \toprule
        $\tilde r$ & $\tilde R$ & \textbf{MSE} & \textbf{MaxAE} & \textbf{MRE} & \textbf{MaxRE} & \textbf{ExTime} \\ \toprule
        25 & 351 & 0.00162521 & 0.125222 & 0.00107555 & 0.0256307 & 107.901 \\ \hline
        30 & 496 & 0.00105956 & 0.0880324 & 0.000706411 & 0.0108741 & 202.808 \\ \hline
        35 & 666 & 0.000854884 & 0.0848942 & 0.000559544 & 0.00937306 & 352.785 \\ \hline
        40 & 861 & 0.000718592 & 0.103829 & 0.000440516 & 0.0176728 & 571.705 \\ \hline
        45 & 1081 & 0.00064009 & 0.0870424 & 0.000383435 & 0.0112156 & 884.359 \\ \hline
        50 & 1326 & 0.000841451 & 0.146591 & 0.000350139 & 0.0330292 & 1308.41 \\ \hline
        55 & 1596 & 0.00170031 & 0.299026 & 0.000369271 & 0.0677911 & 1869.45 \\ \bottomrule
    \end{tabular}
    \caption{Results varying $\tilde r$ for $f_2$ on the polygon}
    \label{tab:results-f2-polygon-r}
\end{table}

\clearpage

In conclusion, we observe that there is no universally optimal choice of $m$ or $\tilde r$, as different combinations may perform differently depending on the approximated function, the domain, and even on the parameters of the chosen domain (recall that we are considering an ellipse of semiaxes $1.5$ and $1$, an annulus of radii $1/4$ and $1$ and a regular polygon of $12$ sides). See Table \ref{tab:best-results}, which reports, for each function $f_1,\dots,f_6$ on the ellipse, the annulus and the polygon, the lowest achieved MSE together with the corresponding values of $m$ and $\tilde r$. 

Moreover, these choices may also depend on the specific problem being addressed. In particular, if computational time is not a limiting factor, then larger values of $m$ may be preferred, as they tend to improve the accuracy, although at the cost of increased complexity and higher computational time, both in constructing the operator and in evaluating it. Conversely, smaller values of $m$ and $\tilde r$ lead to a more efficient, but typically less accurate, approximation.

\begin{table}[h]
\centering
\begin{tabular}{ccccccc}\toprule
$(m,\tilde r)$ & $f_1$                                                          & $f_2$                                                           & $f_3$                                                           & $f_4$                                                           & $f_5$                                                          & $f_6$                                                           \\ \toprule
Ellipse        & \begin{tabular}[c]{@{}c@{}}5.0344e-16 \\ (20,24)\end{tabular}  & \begin{tabular}[c]{@{}c@{}}8.51465e-16 \\ (20, 24)\end{tabular} & \begin{tabular}[c]{@{}c@{}}3.80758e-16 \\ (20, 30)\end{tabular} & \begin{tabular}[c]{@{}c@{}}1.93682e-13 \\ (20, 50)\end{tabular} & \begin{tabular}[c]{@{}c@{}}2.8102e-16 \\ (15, 18)\end{tabular} & \begin{tabular}[c]{@{}c@{}}2.34182e-14 \\ (20, 45)\end{tabular} \\ \hline
Annulus        & \begin{tabular}[c]{@{}c@{}}0.00330492 \\ (20, 40)\end{tabular} & \begin{tabular}[c]{@{}c@{}}0.00319141\\ (30, 40)\end{tabular}   & \begin{tabular}[c]{@{}c@{}}0.000446446\\ (35, 40)\end{tabular}  & \begin{tabular}[c]{@{}c@{}}0.000431085\\ (20, 40)\end{tabular}  & \begin{tabular}[c]{@{}c@{}}0.0142244\\ (20, 40)\end{tabular}   & \begin{tabular}[c]{@{}c@{}}0.00044145\\ (20, 45)\end{tabular}   \\ \hline
Polygon        & \begin{tabular}[c]{@{}c@{}}0.000567284\\ (20,40)\end{tabular}  & \begin{tabular}[c]{@{}c@{}}0.000707322\\ (20, 40)\end{tabular}  & \begin{tabular}[c]{@{}c@{}}0.00130365\\ (20, 40)\end{tabular}   & \begin{tabular}[c]{@{}c@{}}0.000749607\\ (20, 50)\end{tabular}  & \begin{tabular}[c]{@{}c@{}}0.000707498\\ (20, 50)\end{tabular} & \begin{tabular}[c]{@{}c@{}}0.00108535\\ (20, 40)\end{tabular} \\ \bottomrule
\end{tabular} 
\caption{Best MSE achieved for each function and each domain, together with the corresponding pair $(m,\tilde r)$.}
\label{tab:best-results}
\end{table}

After evaluating the performance of the interpolation–regression operators for function approximation under different parameter configurations, we now turn to their application to numerical integration via Gaussian cubature rules, highlighting their potential beyond approximation.

\subsection{Cubature formulas on planar surfaces}\label{sec:num-exp-cubature}
\begin{figure}
    \centering
    \begin{tabular}{cc}
        \includegraphics[height=3.5cm]{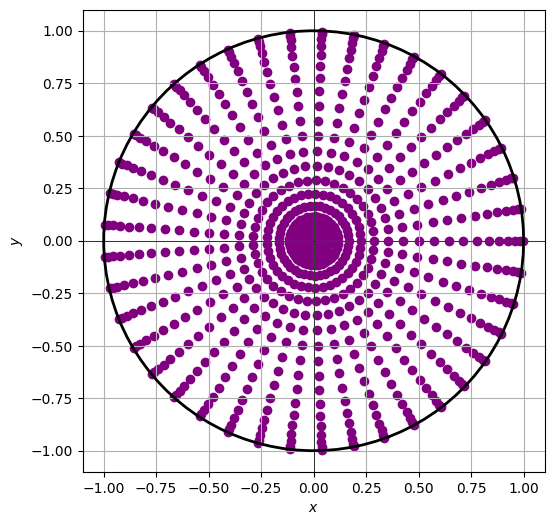} & \includegraphics[height=3.5cm]{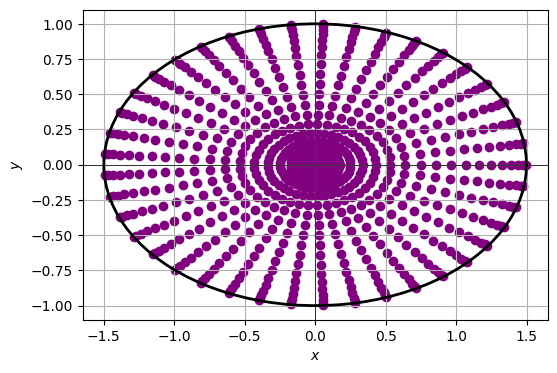} \\
        (a) Unit disk & (b) Ellipse with semiaxes $A=1.5$, $B=1$ \\
        \includegraphics[height=3.5cm]{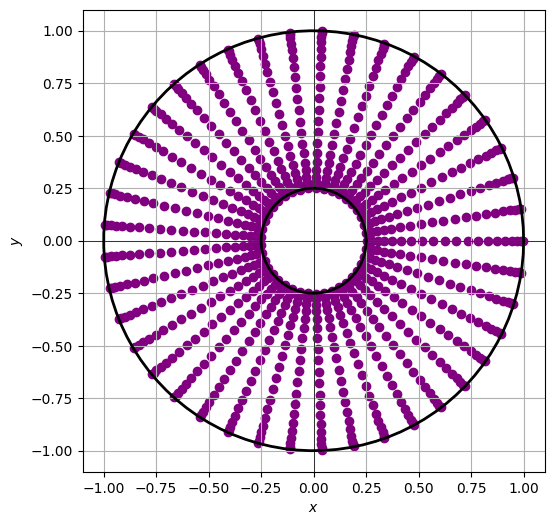} & \includegraphics[height=3.5cm]{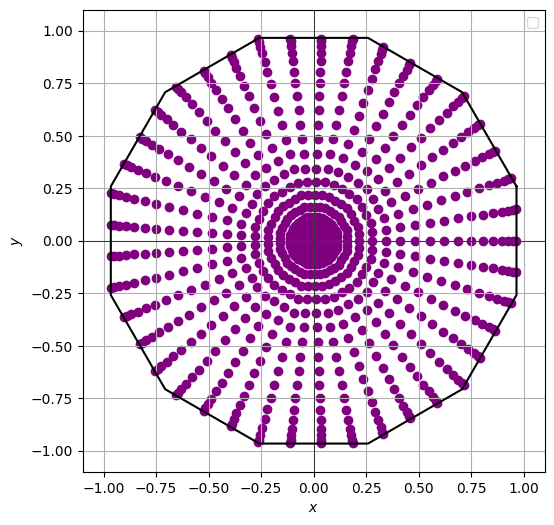} \\
        (c) Annulus with radii $a=1/4$, $A=1$ & (d) Regular polygon with $p=12$ sides
    \end{tabular}
    \caption{Cubature nodes on the disk, the ellipse, the annulus, and the polygon.}
    \label{fig:cubature-nodes}
\end{figure}

Following the methodology introduced in Section~\ref{sec:cubature}, we compute Gaussian nodes on the unit disk with an algebraic degree of exactness equal to $40$, as shown in Figure~\ref{fig:cubature-nodes}(a). These computations are performed by adapting a \href{https://www.math.unipd.it/~alvise/POINTSETS/DISK/cub_disk_prodrule.m}{MATLAB script} from \cite{Som26}. As detailed in Section~\ref{sec:cubature}, the cubature nodes on the ellipse, annulus, and polygon are obtained by applying the corresponding diffeomorphism $\varphi$ between the disk and each target domain; see Figure~\ref{fig:cubature-nodes} (b)-(d). 

By using the cubature nodes and their associated weights defined in \eqref{eq:cubature-ellipse}, \eqref{eq:cubature-annulus}, and \eqref{eq:cubature-polygon}, we approximate the integrals of the functions $f_1, \dots, f_6$ over these domains via the cubature rule \eqref{eq:cubature}. In addition, we consider the constant function $f_0(x,y)=1$, whose integral yields the area of $\Omega$. For the underlying interpolation-regression operator, the parameters are set to $n=100$, $m=20$, and $\tilde r = 24$. Although alternative parameter configurations can be considered, it is important to note that increasing $m$ and $\tilde r$ increases the structural complexity of the operator. 

For each function $f_0, \dots, f_6$, Tables~\ref{tab:cubature-ellipse}, \ref{tab:cubature-annulus}, and \ref{tab:cubature-polygon} summarize the results obtained for the ellipse, annulus, and polygon, respectively, including:
\begin{itemize}
    \item The reference value of the integral $\displaystyle\int_\Omega f_i(\mathbf x) \, d\mathbf x $,  computed numerically using the function \href{https://docs.scipy.org/doc/scipy/reference/generated/scipy.integrate.dblquad.html}{\texttt{dblquad}}.

    \item Since \texttt{dblquad} performs a numerical computation, it also provides an estimate of the associated error, which is included in the tables.
    
    \item The approximation obtained using the cubature formulas based on the interpolation-regression operator \eqref{eq:cubature}.
    
    \item The execution time required to evaluate the cubature formula \eqref{eq:cubature}.
    
    \item The squared error between the reference value of the integral and its approximation.
    
    \item The relative error between the reference value of the integral and its approximation. When the reference value is exactly $0$, a hyphen (`-') is displayed.
\end{itemize}

\begin{table}[!ht]
    \centering
    \begin{tabular}{ccccccc}
    \toprule
        $f$ & \textbf{Actual Int.} & \textbf{Estimated error} & \textbf{Cubature} & \textbf{ExTime} & \textbf{Sq. Error} & \textbf{Rel. Error} \\ \toprule
        $f_0$ & 4.71239 & 3.48787e-14 & 4.71239 & 7.84301 & 2.84779e-28 & 3.58107e-15 \\ \hline
        $f_1$ & -1.13583e-17 & 8.03213e-15 & -3.65701e-16 & 7.75366 & 1.25559e-31 & - \\ \hline
        $f_2$ & 4.93799 & 1.18395e-12 & 4.93799 & 8.37978 & 6.38977e-29 & 1.6188e-15 \\ \hline
        $f_3$ & 2.37736 & 1.75964e-14 & 2.37736 & 9.54007 & 4.93038e-30 & 9.33996e-16 \\ \hline
        $f_4$ & 2.83274 & 2.8456e-13 & 2.83274 & 8.39214 & 1.07934e-20 & 3.66751e-11 \\ \hline
        $f_5$ & 2.03608e-17 & 1.08611e-14 & 1.92527e-16 & 9.60892 & 2.96413e-32 & - \\ \hline
        $f_6$ & 2.59994 & 1.9244e-14 & 2.59994 & 9.803 & 4.7276e-22 & 8.3629e-12 \\ \bottomrule
    \end{tabular}
    \caption{Integral approximations obtained using cubature formulas and interpolation-regression operators on the ellipse}
    \label{tab:cubature-ellipse}
\end{table}

\begin{table}[!ht]
    \centering
    \begin{tabular}{ccccccc}
    \toprule
        $f$ & \textbf{Actual Int.} & \textbf{Estimated error} & \textbf{Cubature} & \textbf{ExTime} & \textbf{Sq. Error} & \textbf{Rel. Error} \\ \toprule
        $f_0$ & 2.94524 & 4.35984e-14 & 2.94524 & 7.07569 & 1.59744e-29 & 1.35704e-15 \\ \hline
        $f_1$ & -1.90716e-17 & 7.23667e-15 & -0.000100917 & 7.14614 & 1.01843e-08 & - \\ \hline
        $f_2$ & 3.01129 & 7.37973e-14 & 3.01136 & 7.41202 & 4.591e-09 & 2.25009e-05 \\ \hline
        $f_3$ & 1.79553 & 2.65791e-14 & 1.79616 & 7.8618 & 4.08461e-07 & 0.000355946 \\ \hline
        $f_4$ & 1.98713 & 2.94154e-14 & 1.98773 & 7.95165 & 3.6419e-07 & 0.000303695 \\ \hline
        $f_5$ & 4.77142e-17 & 1.62111e-14 & 0.000450718 & 9.63817 & 2.03147e-07 & - \\ \hline
        $f_6$ & 1.20757 & 1.78756e-14 & 1.20693 & 8.19103 & 4.09854e-07 & 0.000530155 \\ \bottomrule
    \end{tabular}
    \caption{Integral approximations obtained using cubature formulas and interpolation-regression operators on the annulus.}
    \label{tab:cubature-annulus}
\end{table}

\begin{table}[!ht]
    \centering
    \begin{tabular}{ccccccc}
    \toprule
        $f$ & \textbf{Actual Int.} & \textbf{Estimated error} & \textbf{Cubature} & \textbf{ExTime} & \textbf{Sq. Error} & \textbf{Rel. Error} \\ \toprule
        $f_0$ & 3 & 3.79433e-09 & 2.99996 & 8.68664 & 1.90953e-09 & 1.45661e-05 \\ \hline
        $f_1$ & 2.99948e-18 & 4.92229e-15 & 3.82355e-05 & 8.68612 & 1.46195e-09 & - \\ \hline
        $f_2$ & 3.05756 & 6.98765e-07 & 3.05745 & 9.95199 & 1.22847e-08 & 3.625e-05 \\ \hline
        $f_3$ & 1.93235 & 2.92238e-10 & 1.93233 & 11.2689 & 4.45636e-10 & 1.09246e-05 \\ \hline
        $f_4$ & 2.10582 & 1.14705e-09 & 2.10578 & 15.0048 & 1.19337e-09 & 1.64046e-05 \\ \hline
        $f_5$ & 5.06915e-17 & 1.16798e-14 & 2.76049e-05 & 13.3739 & 7.62032e-10 & - \\ \hline
        $f_6$ & 1.11736 & 4.16171e-09 & 1.11732 & 12.6922 & 1.3748e-09 & 3.31839e-05 \\ \bottomrule
    \end{tabular}
    \caption{Integral approximations obtained using cubature formulas and interpolation-regression operators on the polygon}
    \label{tab:cubature-polygon}
\end{table}

The previous numerical findings naturally lead to the following conclusions regarding the behavior and applicability of the proposed methods.

\section{Concluding remarks and open problems}

In this paper, we have analyzed interpolation–regression methods for functions defined on planar domains. A key aspect of the approach is the existence of a diffeomorphism between these domains and the unit disk, which allows the construction of orthogonal bases by transporting Zernike polynomials to the domain under consideration.

Numerical experiments illustrate the performance of the proposed methods, focusing on the accuracy of the operator. This is assessed through an error analysis based on several metrics and different choices of the dimensions of the interpolation and regression parts.

The analysis of the above methods for Hermite interpolation in these domains constitutes a natural extension of the results obtained in the framework of the constrained mock-Chebyshev least-squares approximation on the real line; see \cite{DMN26}.

Further developments of these interpolation-regression methods for functions defined on domains in $d$-dimensional Euclidean space, such as the $d$-ball, the $d$-ellipse, the $d$-circular annulus, the $d$-dimensional unit cube $[0,1]^{d}$ and the $d$-simplex, respectively, constitute the subject of ongoing work.

\section{Acknowledgments}

The authors appreciate the computing time provided by the ``Servicio de Supercomputación de la Universidad de Granada'' (\url{https://supercomputacion.ugr.es}).

The work of L. Fernández and J. A. Villegas is partially supported by grants PID2023-149117NB-I00, PID2024-155133NB-I00 and CEX
2020-001105-M, all funded by ``Ministerio de Ciencia, Innovación y Universidades'' (MICIU/AEI/10.13039/501100011033 and ERDF/EU), Spain.

The work of F. Marcellán has been supported by the research project PID2024-155133NB-I00, \emph{Ortogonalidad, Aproximaci\'on e Integrabilidad: Aplicaciones en Procesos Estoc\'asticos Cl\'asicos y Cu\'anticos} funded by MCIU/AEI, Spain.

\bibliographystyle{plain}
\bibliography{references}{}

\end{document}